\documentclass{article}
\usepackage{amsmath}
\usepackage{amsopn}
\usepackage{amssymb}
\usepackage{latexsym}
\usepackage{amsfonts}
\usepackage[all]{xy}
\usepackage{amsthm}
\usepackage{graphicx}

\newtheorem{thm}{Theorem}[section]

\newtheorem{pro}[thm]{Proposition}
\newtheorem{lem}[thm]{Lemma}
\newtheorem{cla}[thm]{Claim}

\newtheorem{cor}[thm]{Corollary}

\theoremstyle{definition}
\newtheorem{obs}[thm]{Observation}

\newtheorem{exa}[thm]{Example}
\newtheorem{defn}[thm]{Definition}

\newtheorem{conj}[thm]{Conjecture}

\newcommand{\een}{\end{enumerate}}
\newcommand{\blem}{\begin{lem}}
\newcommand{\elem}{\end{lem}}
\newcommand{\bcl}{\begin{cla}}
\newcommand{\ecl}{\end{cla}}
\newcommand{\ethm}{\end{thm}}
\newcommand{\bpr}{\begin{pro}}
\newcommand{\epr}{\end{pro}}
\newcommand{\bco}{\begin{cor}}
\newcommand{\eco}{\end{cor}}
\newcommand{\bcon}{\begin{conj}}
\newcommand{\econ}{\end{conj}}
\newcommand{\bde}{\begin{defn}}
\newcommand{\ede}{\end{defn}}
\newcommand{\bex}{\begin{exa}}
\newcommand{\eexa}{\end{exa}}
\newcommand{\bobs}{\begin{obs}}
\newcommand{\eobs}{\end{obs}}
\newcommand{\bexe}{\begin{exe}}
\newcommand{\eexe}{\end{exe}}

\begin{document}
\title{Generalization of Menger's Edge Theorem to Four Vertices}
\author{Avraham Goldstein\\ BMCC-CUNY, New York, USA; avraham.goldstein.nyc@gmail.com}
\maketitle
\begin{abstract}
Menger's Edge Theorem asserts that there exist $k$ pairwise edge-disjoint paths between two vertices in an undirected graph $G$ if and only if a deletion of any $k-1$ or less edges does not disconnect these two vertices. Alternatively, there exist $k$ pairwise summand-disjoint formal sums of edges with coefficients in $\mathbb{F}_2$, which are mapped by the boundary map to the sum of two vertices in $G$, if and only if after a deletion of any $k-1$ or less edges there still exist a formal sum of edges with coefficients in $\mathbb{F}_2$, which is mapped by the boundary map to the sum of these two vertices.
\\ \\
We extend this result to four vertices $A,B,C,D$. We prove that in an undirected graph $G$, in which all the vertices different from the vertices $A,B,C,D$ have even degrees, the following two statements are equivalent: There exist $k$ pairwise summand-disjoint formal sums of edges with coefficients in $\mathbb{F}_2$, which are mapped by the boundary map to $A+B+C+D$; After a deletion of any $k-1$ or less edges there still exists a formal sum of edges with coefficients in $\mathbb{F}_2$, which is mapped by the boundary map to $A+B+C+D$.
\\ \\
Equivalently, if after a deletion of any $k-1$ or less edges, the four vertices $A,B,C,D$ can be split into two pairs of vertices and the two vertices in each pair then can be connected by a path in such a way that these two paths are edge-disjoint, then the four vertices $A,B,C,D$ can be split $k$ times into two pairs of vertices and the two vertices in each one of these $2k$ pairs then connected by a path in such a way that these $2k$ paths are pairwise edge-disjoint.
\end{abstract}
\section{Introduction}
Menger's Edge Theorem is a fundamental and a classical result in the Graph Theory, published by Carl Menger in 1927, which has important applications in Flow Networks and other areas. Menger's Edge Theorem is usually viewed in the Graph Theory as stating that for every partition $V_1,V_2$ of the vertex set $V$ of a graph $G$ into $2$ non-empty disjoint sets, such that $V_1$ contains some fixed vertex $A$ and $V_2$ contains some fixed vertex $B$, there are at least $k$ crossing edges between $V_1$ and $V_2$ if and only if $G$ contains $k$ pairwise edge-disjoint trees, each one of which contains both vertices $A$ and $B$. This viewpoint of the Menger's Edge Theorem led to a generalization, which we will now describe.
\\ \\
In 1961 Crispin St. John Alvah Nash-Williams and William Thomas Tutte independently published a theorem, asserting that a graph $G$ has $k$ pairwise edge-disjoint spanning trees, which are the trees composed of all the vertices and of some of the edges of the graph, if and only if for every partition $V_1,\ldots,V_t$ of the vertex set $V$ of $G$ into $t$ non-empty pairwise disjoint sets, there are at least $k(t-1)$ crossing edges, which are the edges connecting between any two vertices belonging to the different sets in the partition. A corollary of this theorem is that if after a deletion of any $2k-1$ or less edges from the graph $G$, the graph $G$ still remains connected, then $G$ contains $k$ pairwise edge-disjoint spanning trees. These are classical results in the Graph Theory, and we refer the reader to \cite{Diestel}, \cite{Harary}, \cite{West}.
\\ \\
In 1999 Matthias Kriesell published \cite{Kriesell1}, in which he conjectured, that if a subset $U$ of the vertex set $V$ of a graph $G$ remains connected in $G$ after any deletion of $2k-1$ or less edges from $G$, then $G$ contains $k$ pairwise edge-disjoint trees, each of which contains all the vertices belonging to $U$. In that work Kriesell proved his conjecture for the case when all the vertices from $V-U$ have even degrees in $G$. In 2010 Matthias Kriesell published \cite{Kriesell2}, in which he proved that if a subset $U$ of the vertex set $V$ of a graph $G$ consisting of four vertices, and if after a deletion of any $\lceil\frac{3k}{2}\rceil-1$ or less edges from $G$, $U$ still remains connected in $G$, then $G$ contains $k$ pairwise edge-disjoint trees, each of which contains all the four vertices belonging to $U$. Here $\lceil\frac{3k}{2}\rceil$ is the least integer, which is not smaller than $\frac{3k}{2}$. In that work Kriesell proved, that $\lceil\frac{3k}{2}\rceil$ is the best possible bound.
\\ \\
However, we suggest to view Menger's Edge Theorem as a statement about the chains of edges in $G$, chains of vertices in $G$, and the boundary map from the former to the latter. Our motivation behind this approach is that Menger's Edge Theorem is really a theorem about the topological space which is a geometric realization of the graph $G$, and it relates between certain topological properties of that space. As such, it is independent of the particular choice of the graph $G$.
\\ \\
To explain our motivation we in the Appendix of this work give a brief summary of the topological viewpoint of the Graph Theory. In it we define the topological space obtained from the geometric realization of a graph, provide a rigorous treatment of these topological spaces and their properties, and state and prove the topological versions of the Menger's Edge Theorem and of our Theorem \ref{Menger.homol.5}.
\section{Preliminaries: Paths, Chains, Cycles and Boundaries in Graphs}
In this work all the algebra is done over the field $\mathbb{F}_2$ and all our graphs are finite and undirected. A graph $G$ is a triple, consisting of a set $E$, which is called the edge set of $G$, a nonempty set $V$, which is called the vertex set of $G$, and a function $\partial$ from the set $E$ to the set of all the subsets of $V$ consisting of exactly two different elements. The elements of $V$ are called vertices in $G$, the elements of $E$ are called edges in $G$, and $\partial$ is called the boundary function in $G$. In our graph we can have some edges $e_1\ne e_2$, such that $\partial(e_1)=\partial(e_2)$, but we cannot have a loop, which is an edge whose boundary is two copies of the same vertex.
\bde The degree $\deg(v)$ of a vertex $v$ in a graph $G$ is the total number of different edges $e$ in $G$ such that the set $\partial(e)$ contains the vertex $v$.
\ede
The degree of a vertex in a graph, clearly, always is a non-negative integer. If we wanted to permit loops in our graphs, we would define the boundary function $\partial$ in $G$ as going from $E$ to $S$, where $S$ is a set of all the sets consisting of two elements, each one of which is a copy of an element in $V$. In that case, the degree of a vertex $v$ in a graph $G$ is the total number of copies of $v$ which appear as elements in all the boundaries of all the edges in $G$. However, deleting all the loops from such a graph has no bearing on it with respect to the requirements and the statements of our Theorems \ref{Menger.homol.4}, \ref{Menger.homol.5} and \ref{Menger.homol.6}. That implies that these theorems are correct for the graphs with loops.
\bde Let $S$ be a set. An $\mathbb{F}_2$-linear space of all the formal $\mathbb{F}_2$-linear combinations of elements of $S$ is called the chain space of $S$ and is denoted by $A(S)$. If $S=\emptyset$ then $A(S)$ is the $0$-dimensional $\mathbb{F}_2$-linear space consisting of the zero element only. The elements of $A(S)$ are called chains of elements of $S$ or chains of elements in $S$.
\ede
There is a natural one-to-one correspondence between the subsets of $S$ and the elements of $A(S)$, which corresponds to each subset of $S$ the sum of all the elements of that subset. Under that correspondence the empty subset of $S$ corresponds to $0\in A(S)$. Moreover, the $XOR$ operation on the subsets of $S$, which we denote by $\oplus$, corresponds to the addition operation in $A(S)$. Abusing the notation, we will speak of the subsets of $S$ and of the corresponding elements of $A(S)$ interchangeably.
\\ \\
For any graph $G$, the function $\partial$ extends to a function from subsets of $E$ to subsets of $V$ by mapping the empty set to the empty set and mapping the set $\{e_1,\ldots,e_m\}$ to the set $\{\partial(e_1)\oplus\cdots\oplus\partial(e_m)\}$. Abusing the notation, we call this function $\partial$. Additionally, the function $\partial$, which is regarded as mapping each edge $e$ to the sum in $A(V)$ of the two elements of the set $\partial(e)$, extends by $\mathbb{F}_2$-linearity to a function from $A(E)$ to $A(V)$, and again we call this extension $\partial$.
\\ \\
Let $G$ be an undirected graph and $u,v$ be any two vertices in $G$. A path between $u$ and $v$ in $G$ is an alternating sequence $$P=(w_0,e_1,w_1,\ldots,w_{t-1},e_t,w_t)$$ of vertices and edges in $G$, with $w_0=u$ and $w_t=v$, or $w_0=v$ and $w_t=u$, such that for each edge $e_i$, $\partial(e_i)=\{w_{i-1},w_i\}$. We also require that $e_i\ne e_j$ for $i\ne j$. However, we do permit $w_i=w_j$ for $i\ne j$. We say that the length of the path $P$ is the number $t$ of edges in it. If $u=v$ then that path is called a cycle. We permit trivial cycles $C=(u)$ of the length $0$ between a vertex $u$ and itself.
\bde
For $k\ge 0$, we call vertices $u$ and $v$ $k$-path-connected in $G$ if there are $k$ paths $P_1,\ldots,P_k$ in $G$ between $u$ and $v$, such that no two of these paths have any common edges.
\ede
Since any trivial cycle $C=(u)$ has no common edges with itself, every vertex $u$ is $k$-path-connected to itself for all $k$. If vertices $u$ and $v$ are $1$-path-connected they are also called connected, and if not they are called disconnected.
\bde
A homological path between vertices $u$ and $v$ in $G$ is a chain of edges $p\in A(E)$ such that $\partial(p)=u+v$. A homological cycle is a chain of edges $c\in A(E)$ such that $\partial(c)=0$. We permit the trivial homological cycle $0\in A(E)$. Every homological cycle is regarded as a homological path between any vertex $u$ in $G$ and itself.
\ede
To each cycle $C=(u)$ we correspond the trivial homological cycle $c=0$, and to each path $P=(w_0,e_1,w_1,\ldots,w_{t-1},e_t,w_t)$ with $t\ge 1$ we correspond the homological path $p=e_1+\cdots+e_t$, for which $\partial(p)=w_0+w_t$.
\\ \\
In the other direction such a correspondence becomes vague. If we are given a homological path $p$ such that $\partial(p)=u+v$ with $u\ne v$, we can always find some edges $e_1,\ldots,e_t$ which appear as summands in $p$ such that there is a path $$P=(u=w_0,e_1,w_1,\ldots,w_{t-1},e_t,w_t=v)$$
This statement is trivial if $p$ consists of only one edge $e_1$, and is proved by induction on the number $m$ of summands in $p$. Thus, assume that for any homological path $p$ with $\partial(p)=u+v$, where $u\ne v$, such that $p$ contains $m$ or less summands in it, we know how to find a path $$P=(u=w_0,e_1,w_1,\ldots,w_{t-1},e_t,w_t=v)$$ such that the edges $e_1,\ldots,e_t$ are summands in $p$. Next, for any two different vertices $u$ and $v$, let $p$ be a homological path such that $\partial(p)=u+v$, and such that $p$ contains $m+1$ summands in it. We find a summand in $p$, which we call $e$, for which $\partial(e)=w+v$, where $w$ is any vertex in the graph. The homological path $p'=p+e$ contains only $m$ summands in it, and we get $\partial(p')=(u+v)+(w+v)=u+w$. If $u=w$ then we take the path $P=(u=w_0,e_1=e,w_1=v)$. If $u\ne w$ then we can find some path $$P'=(u=w_0,e_1,w_1,\ldots,w_{t-1},e_t,w_t=w)$$ such that all the edges $e_1,\ldots,e_t$ appear as summands in $p'$. We take the path $P$ to be $$P=(u=w_0,e_1,w_1,\ldots,w_{t-1},e_t,w_t=w,e_{t+1}=e,w_{t+1}=v)$$
Notice, that there can be many different choices of edges $e_1,\ldots,e_t$ which will be used in the path $P$, and for each such a choice there can be different possible orderings of these edges in $P$. For example, for the homological path $p=h_1+h_2+h_3$ such that $\partial(h_1)=\partial(h_2)=\partial(h_3)=u+v$ with $u\ne v$ we can take the paths $P'_1=(u=w_0,e_1=h_1,w_1=v)$, $P'_2=(u=w_0,e_1=h_2,w_1=v)$, $P'_3=(u=w_0,e_1=h_3,w_1=v)$, $P_1=(u=w_0,e_1=h_1,w_1=v,e_2=h_2,w_2=u,e_3=h_3,w_3=v)$, $P_2=(u=w_0,e_1=h_1,w_1=v,e_2=h_3,w_2=u,e_3=h_2,w_3=v)$, or any one of the other four paths which are obtained from the permutations of the order of the three edges $h_1,h_2,h_3$.
\\ \\
In general, there can be some summands in $p$ which cannot appear in any path $P$ between the two vertices $u$ and $v$. Also, there can be some summands in $p$ which cannot appear in a path $P$ between the vertices $u$ and $v$ if $P$ is composed of the edges which appear as summands in $p$. For example, if $p=h_1+h_2+h_3$ such that $\partial(h_1)=u+v$ and $\partial(h_2)=\partial(h_3)=u'+v'$, where all the vertices $u,u',v,v'$ are pairwise different, then a path $P$, which goes between the vertices $u$ and $v$ and which is composed of the edges which appear as summands in $p$, cannot contain the edges $h_2$ and $h_3$ in it. Furthermore, if the graph $G$ consists of four vertices $u,u',v,v'$ and of four edges $h_1,h_2,h_3,h_4$, where $\partial(h_1),\partial(h_2),\partial(h_3)$ are as above, and $\partial(h_4)=v+v'$, then any path $P$ between the vertices $u$ and $v$ in $G$ cannot contain the edges $h_2$ and $h_3$ in it.
\\ \\
Clearly, if we have a homological path $p$ such that $\partial(p)=u+v$ with $v\ne u$, and we take any path $$P=(u=w_0,e_1,w_1,\ldots,w_{t-1},e_t,w_t=v)$$ and create a new homological path $c=p+e_1+\cdots+e_t$ then we get a homological cycle. If all the edges $e_1,\ldots,e_t$ were summands in $p$ then this homological cycle $c$ is obtained from $p$ by removing these summands from $p$.
\\ \\
Now we discuss the correspondence between the homological cycles and the regular cycles. To the trivial homological cycle $c=0$ we can correspond any cycle $C=(u)$ where $u$ is any vertex in the graph. Since we do not permit loops, a homological cycle cannot contain only one summand and a cycle cannot contain only one edge.
\\ \\
If we are given a homological cycle $c=h_1+\cdots+h_m$ with $m\ge 2$ then we can always find some cycles $C_1,\cdots,C_r$, each one of them containing two or more edges, such that no two of these cycles have any common edges, and such that $c$ is the sum of all the edges which appear in $C_1,\cdots,C_r$. This is proved by induction on the number $m$ of summands in $c$. Indeed, if $c=h_1+h_2$ contains two summands in it, then we must have $\partial(h_1)=\partial(h_2)=u+w$, and we can take $C_1=(u=w_0,e_1=h_1,w_1=w,e_2=h_2,w_2=u)$. If $c$ has more than two summands, then we select any summand $h$ in $c$, and we call one of the vertices in $\partial(h)$ by $u$ and the other vertex in $\partial(h)$ by $v$. Now, $c+h$ is a homological path between $u$ and $v$, and we can find some path $$P=(u=w_0,e_1,w_1,\ldots,w_{t-1},e_t,w_t=v)$$ such that all the edges $e_1,\ldots,e_t$ are summands in $c+h$. Then $$C_1=(u=w_0,e_1,w_1,\ldots,w_{t-1},e_t,w_t=v,e_{t+1}=h,w_{t+1}=u)$$ is a cycle, and all the edges $e_1,\ldots,e_t,e_{t+1}$ which appear in $C_1$ are summands in $c$. Finally, $c'=c+e_1+\cdots+e_t+e_{t+1}$ is a homological cycle with $m-t-1$ summands in it. If $m-t-1=0$ then $c$ is the sum of all the edges which appear in $C_1$. Otherwise, by induction, we can find some cycles $C_2,\cdots,C_r$, such that no two of these cycles have any common edges, and such that $c'$ is the sum of all the edges which appear in $C_2,\cdots,C_r$. Then $c$ is the sum of all the edges which appear in the cycles $C_1,\cdots,C_r$, and no two of these cycles have any common edges.
\bde The image of $A(E)$ in $A(V)$ under $\partial$ is called the boundary space of the graph $G$ and is denoted by $B(G)$. The kernel of $\partial$ in $A(E)$ is called the cycle space of the graph $G$ and is denoted by $Cycle(G)$.
\ede
The quotient $\mathbb{F}_2$-linear space $A(V)/B(G)$ has a dimension $n>0$. Indeed, the image of every edge $e\in E$ under $\partial$ has an even number of vertex summands in it, and cancelations in summations in $A(V)$ always remove an even number of vertex summands. This implies that every element of $B(G)$ contains an even number of vertex summands in it. Hence, $v\notin B(G)$ for any vertex $v\in V$. Thus, one can find some $v_1,\ldots,v_n\in V$ such that $A(V)$ is the union of pairwise disjoint sets $B(G),v_1+B(G),\ldots,v_n+B(G)$. Every $v\in V$ belongs to some unique set $v_i+B(G)$, and $u,v\in V$ are connected in $G$ if and only if they belong to the same set $v_i+B(G)$.
\\ \\
The sets $V_1=v_1+B(G),\ldots,V_n=v_n+B(G)$ are called the connected components of $V$ in $G$. For every $e\in E$, $\partial(e)$ must be a subset of one of the connected components of $V$ in $G$. Thus, $E$ is the union of $n$ pairwise disjoint sets $E_1,\ldots,E_n$, such that $\partial$ takes each $E_i$ to subsets of $V_i$. We obtain a decomposition of the graph $G$ into the disjoint union of the graphs $G_1,\ldots,G_n$, which are called the connected components of $G$. It is easy to see that $B(G)=B(G_1)\oplus\cdots\oplus B(G_n)$. Thus, for $z_i\in A(E_i)$ and $z_j\in A(E_j)$ with $i\ne j$, $\partial(z_i)=\partial(z_j)$ if and only if $\partial(z_i)=\partial(z_j)=0$. It is also clear that all the summands in a homological path $p$, such that $\partial(p)=u+v$, which belong to any connected component of the graph, different from the connected component containing the two vertices $u$ and $v$, must add-up to a homological cycle.
\\ \\
A graph $G$, which consists of only one connected component, and such that $Cycle(G)$ is the trivial group $\{0\}$, is called \textit{a tree}. In a tree there always is a unique path between any two vertices. In a tree which contains one or more edges, a vertex $v$, such that $v$ appears as a summand in $\partial(e)$ of only one unique edge $e$, is called \textit{a leaf}, and this edge $e$ is called \textit{the petiole of the leaf} $v$. Clearly, $v$ is a leaf in a tree if and only if $\deg(v)=1$. It can be shown by induction on the number of edges in a tree, that every tree must have at least two leafs in it. Indeed, this is obvious when a tree contains only one edge or two edges. For a tree with three or more edges in it, if we can find an edge, which is not a petiole of any leaf, then deleting that edge will break the tree into two disjoint trees, each one of which has at least one edge in it. Thus, each one of these two trees must have two or more leafs in it by the induction hypothesis - four or more leafs in total. After restoring the deleted edge, the tree must still have at least two of these four or more leafs.
\\ \\
The notion of homological path $p$ between two different vertices $u$ and $v$ can be generalized to any even collection of different vertices. Thus, we can think of elements $p\in A(E)$ such that $\partial(p)=v_1+v_2+\cdots+v_{2n-1}+v_{2n}$, where $n$ is any positive integer, as generalized homological paths between $2n$ pairwise different vertices $v_1,v_2,\ldots,v_{2n-1},v_{2n}$.
\begin{lem}\label{basic} Let $p\in A(E)$ be such that $\partial(p)=v_1+v_2+\cdots+v_{2n-1}+v_{2n}$. There exist $n$ paths $P_1,\ldots,P_n$, each path $P_i$ going between some vertices $v_{i,1}$ and $v_{i,2}$, where the set $\{v_{1,1},v_{1,2},v_{2,1},v_{2,2},\ldots,v_{n,1},v_{n,2}\}$ is equal to the set $\{v_1,v_2,\ldots,v_{2n-1},v_{2n}\}$, such that all the edges which appear in the paths $P_1,\ldots,P_n$ are summands in $p$, and that no two of the paths $P_1,\ldots,P_n$ have any common edges in them.
\end{lem}
\begin{proof} We prove this lemma by induction on $n$, and for each fixed $n$ by induction on the number $m$ of summands in $p$. For $n=1$ the lemma was already established by us above. Assume that the lemma is correct when the number of summands in $\partial(p)$ is between $2$ and $2n$. We will prove it for the case $\partial(p)=v_1+v_2+\cdots+v_{2n+1}+v_{2n+2}$. The chain of edges $p$ must have at least $n+1$ summands in it, since the boundary of any edge has two summands in it. If $p$ has $n+1$ summands then our lemma is trivial, since each summand in $p$ is an edge whose boundary is a sum of some two unique vertices from the set $\{v_1,v_2,\ldots,v_{2n+1},v_{2n+2}\}$. Assume that the lemma is correct for all $p\in A(E)$ consisting of $m$ or less summands. Let $p\in A(E)$ have $m+1$ summands in it. Since $\partial(p)=v_1+v_2+\cdots+v_{2n+1}+v_{2n+2}$, there must be some summand $e$ in $p$ such that $\partial(e)=w+v_{2n+2}$, where $w\ne v_{2n+2}$ is a vertex in the graph.
\\ \\
If $w=v_j$ for some $v_j\in \{v_1,v_2,\ldots,v_{2n},v_{2n+1}\}$ then we take the path $P_{n+1}$ to be $P_{n+1}=(v_j=w_0,e_1=e,w_1=v_{2n+2})$ and we take $v_{n+1,1}=v_j$ and $v_{n+1,2}=v_{2n+2}$. Next, since $\partial(p+e)$ is a sum of $2n$ vertices, we can find $n$ paths $P_1,\ldots,P_n$, with each path $P_i$ going between some vertices $v_{i,1}$ and $v_{i_2}$, where the set $\{v_{1,1},v_{1,2},v_{2,1},v_{2,2},\ldots,v_{n,1},v_{n,2}\}$ is equal to the set $\{v_1,v_2,\ldots,v_{2n-1},v_{2n},v_{2n+1}\}-\{v_j\}$, such that all the edges which appear in the paths $P_1,\ldots,P_n$ are summands in $p+e$. Finally, paths $P_1,\ldots,P_n,P_{n+1}$ satisfy the requirements of the lemma.
\\ \\
If $w\notin \{v_1,v_2,\ldots,v_{2n},v_{2n+1}\}$ then $\partial(p+e)$ is a sum of $2n+2$ vertices, and the chain of edges $p+e$ contains in it $m$ summands. By our induction hypothesis, we can find $n+1$ paths $P_1,\ldots,P_{n+1}$, with each path $P_i$ going between some vertices $v_{i,1}$ and $v_{i,2}$, where the set $\{v_{1,1},v_{1,2},v_{2,1},v_{2,2},\ldots,v_{n,1},v_{n,2}\}$ is equal to the set $\{v_1,v_2,\ldots,v_{2n+1},w\}$, such that all the edges which appear in the paths $P_1,\ldots,P_{n+1}$ are summands in $p+e$. Then some path $P_j$ among these $n+1$ paths $P_1,\ldots,P_{n+1}$ must go between some vertex $v_j\ne w$ and the vertex $w$. We append the edge $e$ to the path $P_j$ and obtain $n+1$ paths which satisfy the requirements of the lemma.
\end{proof}
\section{Statement and Proof of the Extended Menger's Edge Theorem}
In this section we will state the Extended Menger's Edge Theorem first in the language of generalized homological paths, and prove it using that language, and then in the language of the paths. Before doing all that, we need the following five lemmas and one theorem:
\begin{lem}\label{tech}
Let $G$ be a graph consisting of one connected component. For any positive integer $n$ let $v_1,v_2,\ldots,v_{2n-1},v_{2n}$ be any $2n$ different vertices in $G$. There exists an element $p\in A(E)$ such that $\partial(p)=v_1+v_2+\cdots+v_{2n-1}+v_{2n}$.
\end{lem}
\begin{proof}
If $G$ is not a tree then we find a non-trivial cycle $c$ in $G$, select any edge, which appears as a summand in $c$, and delete this edge from our graph $G$. After that deletion $G$ will still consist of one connected component, and if $G$ at that point did not become a tree, then we repeat this process of finding a cycle and deleting an edge again. After a finite number of repetitions, $G$ will become a tree. Since the tree $G$ contains at least the vertices $v_1$ and $v_2$, it must contain at least one edge. Thus, $G$ contains at least two leafs. Next, if some of the leafs in $G$ do not belong to the set $\{v_1,v_2,\ldots,v_{2n-1},v_{2n}\}$ then we select any leaf $v$ in $G$ not belonging to the set $\{v_1,v_2,\ldots,v_{2n-1},v_{2n}\}$, and we delete this leaf $v$ and its petiole edge from the tree $G$. After this deletion, $G$ will still consist of one connected component and contain all the vertices $\{v_1,v_2,\ldots,v_{2n-1},v_{2n}\}$ and at least one edge in it. If now all the leafs in $G$ belong to the set $\{v_1,v_2,\ldots,v_{2n-1},v_{2n}\}$ then we are done, and if not, we repeat this process of deleting a leaf and its petiole again. After a finite number of repetitions, all the leafs in $G$ will be from the set $\{v_1,v_2,\ldots,v_{2n-1},v_{2n}\}$.
\\ \\
Next, if our lemma is wrong, then there exist some counter-examples to it, and we can take any counter-example, and by the above-described deletions turn it into a tree, such that all its leafs belong to the set $\{v_1,v_2,\ldots,v_{2n-1},v_{2n}\}$. The tree, which we will obtain in such a way, will still be a counter-example to our lemma. Let a tree $G$ with all its leafs belonging to the set $\{v_1,v_2,\ldots,v_{2n-1},v_{2n}\}$ be such a counter-example to our lemma, that any other counter-example to our lemma has the same number of edges as $G$ or more edges than $G$. Thus, $G$ is a counter-example with the minimal number of edges in it. After a possible re-indexing of the vertices $v_1,v_2,\ldots,v_{2n-1},v_{2n}$ we will have that $v_{2n}$ is a leaf in $G$. Let $e$ be the petiole edge of $v_{2n}$. Then $\partial(e)=v+v_{2n}$ for some vertex $v$ in $G$.
\\ \\
If $v$ belongs to the set $\{v_1,v_2,\ldots,v_{2n-1},v_{2n}\}$, then we can assume that $v=v_{2n-1}$. In that case, if $n=1$ then the element $e\in A(E)$ will satisfy $\partial(e)=v_1+v_2+\cdots+v_{2n-1}+v_{2n}$, and if $n>1$ then we take a new tree $T$, which is $G$ with the edge $e$ and the vertex $v_{2n}$ deleted from it. From the minimality of the counter-example $G$ it follows, that for this new tree $T$ we can find an element $p'\in A(E_T)$, where $E_T$ is the set of edges of $T$, such that $\partial(p')=v_1+v_2+\cdots+v_{2n-3}+v_{2n-2}$. This implies $\partial(p'+e)=v_1+v_2+\cdots+v_{2n-3}+v_{2n-2}+v_{2n-1}+v_{2n}$, which shows that $G$ is not a counter-example to our lemma.
\\ \\
If $v$ does not belong to the set $\{v_1,v_2,\ldots,v_{2n-1},v_{2n}\}$, then we take a new tree $T$, which is $G$ with the edge $e$ and vertex $v_{2n}$ deleted from it. From the minimality of the counter-example $G$ it follows, that for this new tree $T$ we can find an element $p'\in A(E_T)$ such that $\partial(p')=v_1+v_2+\cdots+v_{2n-2}+v_{2n-1}+v$. This implies $\partial(p'+e)=v_1+v_2+\cdots+v_{2n-3}+v_{2n-2}+v_{2n-1}+v_{2n}$, which again shows that $G$ is not a counter-example to our lemma. Thus, a counter-example to our Lemma does not exist.
\end{proof}
A direct corollary of Lemma \ref{tech} is that if $G$ is a graph and $v_1,v_2,\ldots,v_{2n-1},v_{2n}$, where $n$ is a positive integer, are any $2n$ different vertices in $G$, then there does not exist an element $p\in A(E)$ such that $\partial(p)=v_1+v_2+\cdots+v_{2n-1}+v_{2n}$ if and only if there is a connected component of $G$ which contains an odd number of vertices from the set $\{v_1,v_2,\ldots,v_{2n-1},v_{2n}\}$.
\begin{lem}\label{minimal1} Let $k$ be a positive integer. Let $G$ be a graph with the vertex set $V=\{v_1,v_2,v_{3},v_{4}\}$ and the edge set $E$ such that the degree of each one of the vertices $v_1,v_2,v_{3},v_{4}$ is at least $k$. Then there exist $k$ elements $p_1,\ldots,p_k$ in $A(E)$, such that no two of $p_1,\ldots,p_k$ have any common nontrivial summands, satisfying $$\partial(p_1)=\cdots=\partial(p_k)=v_1+v_2+v_3+v_4$$
\end{lem}
\begin{proof}
For $i=1,2,3$ and $j=i+1,\ldots,4$, let $\alpha(i,j)$ be the number of edges $e$ such that $\partial(e)=v_i+v_j$. We denote these edges by $e_{i,j,1},\ldots,e_{i,j,\alpha(i,j)}$. Of course, $\alpha(i,j)$ can be equal to $0$, in which case the list $e_{i,j,1},\ldots,e_{i,j,\alpha(i,j)}$ does not actually list any edges in it.
\\ \\
\includegraphics[scale=0.5]{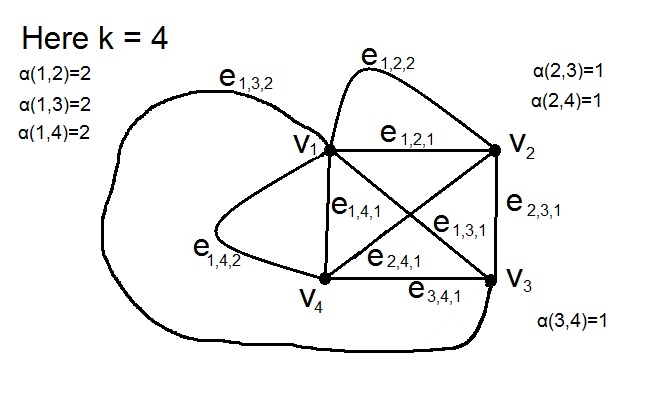}
\\
Let $$m_1=\min(\alpha(1,2),\alpha(3,4))$$ $$m_2=\min(\alpha(1,3),\alpha(2,4))$$ $$m_3=\min(\alpha(1,4),\alpha(2,3))$$ We define $p_t$ as follows:
\\
For $1\le t\le m_1$ we define $p_t=e_{1,2,t}+e_{3,4,t}$;
\\
For $m_1+1\le t\le m_1+m_2$ we define $p_t=e_{1,3,(t-m_1)}+e_{2,4,(t-m_1)}$;
\\
For $m_1+m_2+1\le t\le m_1+m_2+m_3$ we define $p_t=e_{1,4,(t-m_1-m_2)}+e_{2,3,(t-m_1-m_2)}$.
\\
It is easy to verify that $\partial(p_1)=\cdots=\partial(p_{m_1+m_2+m_3})=v_1+v_2+v_3+v_4$. Thus, if $m_1+m_2+m_3\ge k$ then the statement of lemma is true. Assume that $m_1+m_2+m_3<k$.
\\ \\
Since we can re-index the vertices $v_1,v_2,v_{3},v_{4}$, we can always assume that $\alpha(3,4)\le \alpha(1,2)$ and $\alpha(2,4)\le \alpha(1,3)$. If $\alpha(1,4)\le \alpha(2,3)$ we would get $m_1+m_2+m_3=\alpha(3,4)+\alpha(2,4)+\alpha(1,4)<k$, which implies that the degree of $v_4$ is less than $k$. Thus, $\alpha(1,4)>\alpha(2,3)$. Next, if $\alpha(3,4)=\alpha(1,2)$ then $m_1+m_2+m_3=\alpha(1,2)+\alpha(2,4)+\alpha(2,3)<k$, which implies that the degree of $v_2$ is less than $k$. Thus, $\alpha(3,4)<\alpha(1,2)$. Similarly, if $\alpha(2,4)=\alpha(1,3)$ then $m_1+m_2+m_3=\alpha(3,4)+\alpha(1,3)+\alpha(2,3)<k$, which implies that the degree of $v_3$ is less than $k$. Thus, $\alpha(2,4)<\alpha(1,3)$.
\\ \\
Next, let $$m=\min(\alpha(1,4)-\alpha(2,3),\alpha(1,2)-\alpha(3,4),\alpha(1,3)-\alpha(2,4))$$ If $m_1+m_2+m_3+m<k$ then either $$\alpha(3,4)+\alpha(2,4)+\alpha(2,3)+(\alpha(1,4)-\alpha(2,3))<k$$ or $$\alpha(3,4)+\alpha(2,4)+\alpha(2,3)+(\alpha(1,2)-\alpha(3,4))<k$$ or $$\alpha(3,4)+\alpha(2,4)+\alpha(2,3)+(\alpha(1,3)-\alpha(2,4))<k$$ But this implies that either the degree of $v_4$ or the degree of $v_2$ or the degree of $v_3$ is less than $k$. Thus, $m_1+m_2+m_3+m\ge k$.
\\ \\
Finally, for all $t=1,\ldots,m$ we define $$p_{m_1+m_2+m_3+t}=e_{1,4,(m_3+t)}+e_{1,3,(m_2+t)}+e_{1,2,(m_1+t)}$$ It is easy to verify that $\partial(p_{m_1+m_2+m_3+t})=v_1+v_2+v_3+v_4$. Thus, we produced $k$ or more elements in $A(E)$, such that no two of them have any common nontrivial summands, which are mapped by $\partial$ to $v_1+v_2+v_3+v_4$.
\end{proof}
\begin{lem}\label{minimal2} Let $G$ be a graph with the vertex set $V$ such that $V$ is a union of disjoint sets $\{v_1,v_2,v_{3},v_{4}\}$ and some empty or nonempty set $Y$, and the edge set $E$ such that for every $e\in E$, at least one of the two summands in boundary the $\partial(e)$ of $e$ belongs to the set $\{v_1,v_2,v_{3},v_{4}\}$.
\\ \\
If the degree of each one of the vertices $v_1,v_2,v_{3},v_{4}$ is at least $k$, and the degree of every vertex belonging to $Y$ is even, and for each $y\in Y$ and each $1\le i\le 4$, the number of the edges $e$ for which $\partial(e)=y+v_i$ is not greater than the total number of all the edges $a$ for which $\partial(a)=y+w$ with any $w\ne v_i$, then there exist $k$ chains of edges $p_1,\ldots,p_k$ in $A(E)$, such that no two of $p_1,\ldots,p_k$ have any common nontrivial summands, satisfying $$\partial(p_1)=\cdots=\partial(p_k)=v_1+v_2+v_3+v_4$$
\end{lem}
\begin{proof}
If $Y$ is empty then this lemma follows from Lemma \ref{minimal1}. We proof the lemma by induction on the number of vertices in $Y$. Assume that our lemma is correct when $Y$ contains $m$ or less vertices. Suppose that $Y$ contains $m+1$ vertices. We select any vertex $y\in Y$. We will modify the graph $G$ by deleting some edges and the vertex $y$ from it and by drawing some new edges in it in such a way, that the degrees of the vertices $v_1,v_2,v_{3},v_{4}$ will not change and that whenever we find $k$ chains of edges $p'_1,\ldots,p'_k$ in the modified graph $G$, such that no two of $p'_1,\ldots,p'_k$ have any common nontrivial summands, satisfying $\partial(p'_1)=\cdots=\partial(p'_k)=v_1+v_2+v_3+v_4$ in the modified $G$, then we can produce $k$ chains of edges $p_1,\ldots,p_k$ in the original graph $G$, such that no two of $p_1,\ldots,p_k$ have any common nontrivial summands, satisfying $\partial(p_1)=\cdots=\partial(p_k)=v_1+v_2+v_3+v_4$ in the original $G$.
\\ \\
Theoretically, for our vertex $y$ there could exist five different cases:
\\ \\
\textit{Case 1} is when there are no edges $e\in E$ such that $\partial(e)$ contains $y$ as a summand. In that case we delete the vertex $y$ from our graph $G$, and apply the induction hypothesis on the number of vertices in $Y$;
\\ \\
\textit{Case 2} is when for all the edges $e\in E$, such that $\partial(e)$ contains $y$ as a summand, we have $\partial(e)=y+v_i$ for some unique index $i$ between $1$ and $4$. In that case the number of the edges $e$ for which $\partial(e)=y+v_i$, is greater than the total number of all the edges $a$ for which $\partial(a)$ is a sum of $y$ and of any vertex $w$ different from $v_i$. This contradicts the requirements of our lemma, which implies that this case cannot happen;
\\ \\
\textit{Case 3} is when for all the edges $e\in E$, such that $\partial(e)$ contains $y$ as a summand, we have either $\partial(e)=y+v_i$ or $\partial(e)=y+v_j$ for some two fixed indices $i<j$ between $1$ and $4$. In that case the number of the edges $a$ such that $\partial(a)=y+v_i$ must be not greater than the number of the edges $b$ such that $\partial(b)=y+v_j$, and the number of the edges $b$ such that $\partial(b)=y+v_j$ must be not greater than the number of the edges $a$ such that $\partial(a)=y+v_i$. Thus, we will have $n>0$ edges $a_1,\ldots,a_n$ such that $\partial(a_1)=\cdots=\partial(a_n)=y+v_i$, and $n>0$ edges $b_1,\ldots,b_n$ such that $\partial(b_1)=\cdots=\partial(b_n)=y+v_j$. We pair these edges $a_z$ and $b_z$ for all $z=1,\ldots,n$ in pairs, delete each one of these pairs from $G$, and for each deleted pair $a_z,b_z$ we draw a new edge $\gamma_z$ in $G$ such that $\partial(\gamma_z)=v_i+v_j$. After we do that for all the $2n$ edges whose boundaries contain $y$ as a summand, we delete the vertex $y$ from $G$.
\\ \\
We claim that if after this modification of the graph $G$, there exist $k$ chains of edges $p'_1,\ldots,p'_k$ in the modified $G$, such that no two of $p'_1,\ldots,p'_k$ have any common nontrivial summands, satisfying $\partial(p'_1)=\cdots=\partial(p'_k)=v_1+v_2+v_3+v_4$ in the modified $G$, then there exist $k$ chains of edges $p_1,\ldots,p_k$ in the original $G$, such that no two of $p_1,\ldots,p_k$ have any common nontrivial summands, satisfying $\partial(p_1)=\cdots=\partial(p_k)=v_1+v_2+v_3+v_4$ in the original $G$. Indeed, in each place where the chains of edges $p'_1,\ldots,p'_k$ in the modified $G$ contain some new edge $\gamma_z$ as a summand, we substitute that $\gamma_z$ by the sum $a_z+b_z$, and we obtain the required chains of edges $p_1,\ldots,p_k$ in the original graph $G$. Notice, that our modification of the graph $G$ did not change the degrees of the vertices $v_1,v_2,v_{3},v_{4}$, and so the existence of the required $k$ chains of edges $p'_1,\ldots,p'_k$ in the modified graph $G$ follows from our induction hypothesis on the number of vertices in $Y$;
\\ \\
\emph{Case 4} is when for all the edges $e\in E$, such that $\partial(e)$ contains $y$ as a summand, we have $\partial(e)=y+v_i$ or $\partial(e)=y+v_j$ or $\partial(e)=y+v_t$ for some three fixed indices $i<j<t$ between $1$ and $4$. In that case, after a possible re-indexing of the vertices $v_1,v_2,v_{3},v_{4}$, we will have:\\
Zero edges $e$ in $G$ such that $\partial(e)=y+v_4$;\\
$n_3>0$ edges $e_{3,1},\ldots,e_{3,n_3}$ in $G$ such that $\partial(e_{3,1})=\cdots=\partial(e_{3,n_3})=y+v_3$;\\
$n_2\ge n_3$ edges $e_{2,1},\ldots,e_{2,n_2}$ in $G$ such that $\partial(e_{2,1})=\cdots=\partial(e_{2,n_2})=y+v_2$;\\
$n_1\ge n_2$ edges $e_{1,1},\ldots,e_{1,n_1}$ in $G$ such that $\partial(e_{1,1})=\cdots=\partial(e_{1,n_1})=y+v_1$.
\\ \\
Since the degree of the vertex $y$ is even, $n_1+n_2+n_3$ must be even. The requirement that for each $1\le i\le 4$, the number of the edges $e$ for which $\partial(e)=y+v_i$ is not greater than the total number of all the edges $a$ for which $\partial(a)=y+w$ with any $w\ne v_i$, stated in our lemma, is equivalent to the requirement that $n_1$ is not greater than $n_2+n_3$. We get two subcases:
\\ \\
\textit{Subcase 1 of Case 4} is when $n_1<n_2+n_3$. In this subcase we must have $n_1\le n_2+n_3-2$ because $n_1+n_2+n_3$ must be even. Notice that $\varsigma=\frac{n_2+n_3-n_1}{2}=\frac{n_2+n_3+n_1-2n_1}{2}$ is an integer, which is not greater than $\frac{n_3}{2}$ because $n_2\le n_1$. Thus, $n_3-\varsigma>0$ because $n_3\ge 1$. In this subcase we delete $\varsigma$ edges $e_{2,(n_2+1-\varsigma)},\ldots,e_{2,n_2}$ and $\varsigma$ edges $e_{3,(n_3+1-\varsigma)},\ldots,e_{3,n_3}$ from the graph $G$, and draw in $G$ $\varsigma$ new edges $h_{1},\ldots,h_{\varsigma}$ such that $\partial(h_{1})=\cdots=\partial(h_{\varsigma})=v_2+v_3$. In the modified graph $G$, the number of the edges $e$ for which $\partial(e)=y+v_1$ is $n'_1=n_1$, and the number of the edges $e$ for which $\partial(e)=y+v_2$ is $n'_2=n_2-\varsigma$, and the number of the edges $e$ for which $\partial(e)=y+v_3$ is $n'_3=n_3-\varsigma$. Thus, in the modified graph $G$ we get $n'_1=n'_2+n'_3$ with $n'_2\ge n'_3$ and $n'_3>0$. The degrees of the vertices $v_1,v_2,v_{3},v_{4}$ did not change during this modification.
\\
\includegraphics[scale=0.5]{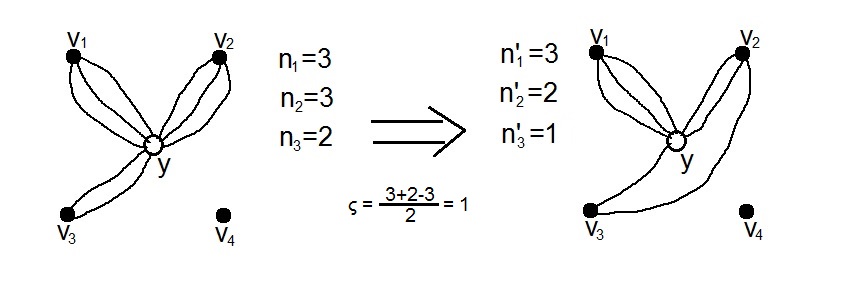}
\\
If in the modified graph $G$ there exist $k$ chains of edges $p'_1,\ldots,p'_k$, such that no two of $p'_1,\ldots,p'_k$ have any common nontrivial summands, satisfying $\partial(p'_1)=\cdots=\partial(p'_k)=v_1+v_2+v_3+v_4$ in the modified $G$, then substituting for each $h_{z}$ which appears as a summand in $p'_1,\ldots,p'_k$ the corresponding sum $e_{2,(n_2+z-\varsigma)}+e_{3,(n_3+z-\varsigma)}$ produces $k$ chains of edges $p_1,\ldots,p_k$ in the original graph $G$, such that no two of $p_1,\ldots,p_k$ have any common nontrivial summands, satisfying $\partial(p_1)=\cdots=\partial(p_k)=v_1+v_2+v_3+v_4$ in the original $G$. Thus, we need to show that we can find the required $k$ chains of edges $p'_1,\ldots,p'_k$ in the modified graph $G$. To do so, we plug the modified graph $G$ as the original graph $G$ in our Subcase 2 of Case 4;
\\ \\
\textit{Subcase 2 of Case 4} is when $n_1=n_2+n_3$. In this subcase we delete $n_2$ edges $e_{1,(n_1+1-n_2)},\ldots,e_{1,n_1}$ and $n_2$ edges $e_{2,1},\ldots,e_{2,n_2}$ from the graph $G$, and draw in $G$ $n_2$ new edges $h_{2,1},\ldots,h_{2,n_2}$ such that $\partial(h_{2,1})=\cdots=\partial(h_{2,n_2})=v_1+v_2$, and we delete $n_3=n_1-n_2$ edges $e_{1,1},\ldots,e_{1,n_3}$ and $n_3$ edges $e_{3,1},\ldots,e_{2,n_3}$ from the graph $G$, and draw in $G$ $n_3$ new edges $h_{3,1},\ldots,h_{3,n_3}$ such that $\partial(h_{3,1})=\cdots=\partial(h_{3,n_3})=v_1+v_3$. In the modified graph $G$ there are no edges which contain the vertex $y$ as a summand in their boundaries. We delete the vertex $y$ from our modified graph $G$. The degrees of the vertices $v_1,v_2,v_{3},v_{4}$ did not change during this entire modification.
\\
\includegraphics[scale=0.5]{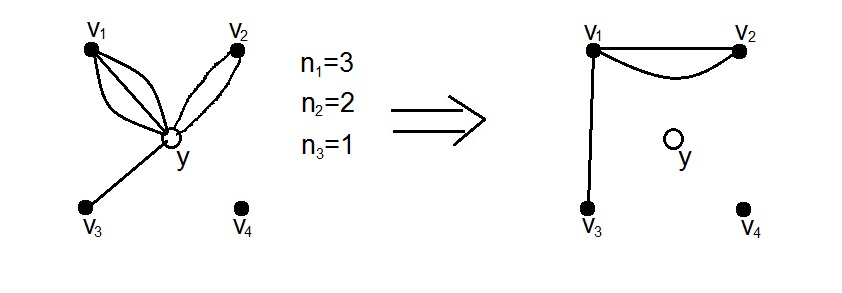}
\\
If in the modified graph $G$ there exist $k$ chains of edges $p'_1,\ldots,p'_k$, such that no two of $p'_1,\ldots,p'_k$ have any common nontrivial summands, satisfying $\partial(p'_1)=\cdots=\partial(p'_k)=v_1+v_2+v_3+v_4$ in the modified $G$, then substituting $e_{1,(n_3+z)}+e_{2,z}$ for each $h_{2,z}$ which appears as a summand in $p'_1,\ldots,p'_k$, and substituting $e_{1,z}+e_{3,z}$ for each $h_{3,z}$ which appears as a summand in $p'_1,\ldots,p'_k$ produces $k$ chains of edges $p_1,\ldots,p_k$ in the original graph $G$, such that no two of $p_1,\ldots,p_k$ have any common nontrivial summands, satisfying $\partial(p_1)=\cdots=\partial(p_k)=v_1+v_2+v_3+v_4$ in the original $G$. Finally, by our induction hypothesis on the number of vertices in $Y$, we can find the required $k$ chains of edges $p'_1,\ldots,p'_k$ in the modified graph $G$. This concludes our Case 4;
\\ \\
\textit{Case 5} is when in $G$ there are:\\
$n_1>0$ edges $e_{1,1},\ldots,e_{1,n_1}$ such that $\partial(e_{1,1})=\cdots=\partial(e_{1,n_1})=y+v_1$;\\
$n_2>0$ edges $e_{2,1},\ldots,e_{2,n_2}$ such that $\partial(e_{2,1})=\cdots=\partial(e_{2,n_2})=y+v_2$;\\
$n_3>0$ edges $e_{3,1},\ldots,e_{3,n_3}$ such that $\partial(e_{3,1})=\cdots=\partial(e_{3,n_3})=y+v_3$;\\
$n_4>0$ edges $e_{4,1},\ldots,e_{4,n_4}$ such that $\partial(e_{4,1})=\cdots=\partial(e_{4,n_4})=y+v_4$.\\
We can arrange our indices in such a way, that $n_1\ge n_2\ge n_3\ge n_4>0$. Since the degree of the vertex $y$ is even, $n_1+n_2+n_3+n_4$ must be even. The requirement that for each $1\le i\le 4$, the number of the edges $e$ for which $\partial(e)=y+v_i$ is not greater than the total number of all the edges $a$ for which $\partial(a)=y+w$ with any $w\ne v_i$, stated in our lemma, is equivalent to the requirement that $n_1$ is not greater than $n_2+n_3+n_4$. Since the degree of the vertex $y$ is even, either $n_1=n_2+n_3+n_4$ or $n_1\le n_2+n_3+n_4-2$.
\\ \\
\textit{Subcase 1 of Case 4} is when $n_1\le n_2+n_3+n_4-2$. Notice that $\frac{n_2+n_3+n_4-n_1}{2}=\frac{n_2+n_3+n_4+n_1-2n_1}{2}$ is a positive integer. Let $\varsigma=\min(\frac{n_2+n_3+n_4-n_1}{2},n_4)$. In this subcase we delete $\varsigma$ edges $e_{3,(n_3+1-\varsigma)},\ldots,e_{3,n_3}$ and $\varsigma$ edges $e_{4,(n_4+1-\varsigma)},\ldots,e_{4,n_4}$ from the graph $G$, and draw in $G$ $\varsigma$ new edges $h_{1},\ldots,h_{\varsigma}$ such that $\partial(h_{1})=\cdots=\partial(h_{\varsigma})=v_3+v_4$. In the modified graph $G$, the number of the edges $e$ for which $\partial(e)=y+v_1$ is $n'_1=n_1$, and the number of the edges $e$ for which $\partial(e)=y+v_2$ is $n'_2=n_2$, and the number of the edges $e$ for which $\partial(e)=y+v_3$ is $n'_3=n_3-\varsigma$, and the number of the edges $e$ for which $\partial(e)=y+v_4$ is $n'_4=n_4-\varsigma$. Thus, in the modified graph $G$ we either get $n'_1=n'_2+n'_3+n'_4$ with $n'_1\ge n'_2> n'_3\ge n'_4\ge 0$, or get $n'_1<n'_2+n'_3$ and $n'_4=0$. The degrees of the vertices $v_1,v_2,v_{3},v_{4}$ did not change during this modification.
\\ \\
If in the modified graph $G$ there exist $k$ chains of edges $p'_1,\ldots,p'_k$, such that no two of $p'_1,\ldots,p'_k$ have any common nontrivial summands, satisfying $\partial(p'_1)=\cdots=\partial(p'_k)=v_1+v_2+v_3+v_4$ in the modified $G$, then substituting for each $h_{z}$ which appears as a summand in $p'_1,\ldots,p'_k$ the corresponding sum $e_{3,(n_3+z-\varsigma)}+e_{4,(n_4+z-\varsigma)}$ produces $k$ chains of edges $p_1,\ldots,p_k$ in the original graph $G$, such that no two of $p_1,\ldots,p_k$ have any common nontrivial summands, satisfying $\partial(p_1)=\cdots=\partial(p_k)=v_1+v_2+v_3+v_4$ in the original $G$. Thus, we need to show that we can find the required $k$ chains of edges $p'_1,\ldots,p'_k$ in the modified graph $G$. If $n'_3=n'_4=0$ we plug the modified graph $G$ as the original graph $G$ in our Case 3. If $n'_3>n'_4=0$ we plug the modified graph $G$ as the original graph $G$ in our Case 4. If $n'_3\ge n'_4>0$  we plug the modified graph $G$ as the original graph $G$ in our Subcase 2 of Case 5;
\\ \\
\textit{Subcase 2 of Case 5} is when $n_1=n_2+n_3+n_4$. In this subcase we delete $n_2$ edges $e_{1,(n_1+1-n_2)},\ldots,e_{1,n_1}$ and $n_2$ edges $e_{2,1},\ldots,e_{2,n_2}$ from the graph $G$, and draw in $G$ $n_2$ new edges $h_{2,1},\ldots,h_{2,n_2}$ such that $\partial(h_{2,1})=\cdots=\partial(h_{2,n_2})=v_1+v_2$, and we delete $n_3$ edges $e_{1,(n_1+1-n_3-n_2)},\ldots,e_{1,(n_1-n_2)}$ and $n_3$ edges $e_{3,1},\ldots,e_{3,n_3}$ from the graph $G$, and draw in $G$ $n_3$ new edges $h_{3,1},\ldots,h_{3,n_3}$ such that $\partial(h_{3,1})=\cdots=\partial(h_{3,n_3})=v_1+v_3$, and we delete $n_4=n_1-n_2-n_3$ edges $e_{1,1},\ldots,e_{1,n_4}$ and $n_4$ edges $e_{4,1},\ldots,e_{4,n_4}$ from the graph $G$, and draw in $G$ $n_4$ new edges $h_{4,1},\ldots,h_{4,n_4}$ such that $\partial(h_{4,1})=\cdots=\partial(h_{4,n_4})=v_1+v_4$. In the modified graph $G$ there are no edges which contain the vertex $y$ as a summand in their boundaries. We delete the vertex $y$ from our modified graph $G$. The degrees of the vertices $v_1,v_2,v_{3},v_{4}$ did not change during this entire modification.
\\ \\
If in the modified graph $G$ there exist $k$ chains of edges $p'_1,\ldots,p'_k$, such that no two of $p'_1,\ldots,p'_k$ have any common nontrivial summands, satisfying $\partial(p'_1)=\cdots=\partial(p'_k)=v_1+v_2+v_3+v_4$ in the modified $G$, then substituting $e_{1,(n_4+n_3+z)}+e_{2,z}$ for each $h_{2,z}$ which appears as a summand in $p'_1,\ldots,p'_k$, and substituting $e_{1,(n_4+z)}+e_{3,z}$ for each $h_{3,z}$ which appears as a summand in $p'_1,\ldots,p'_k$, and substituting $e_{1,z}+e_{4,z}$ for each $h_{4,z}$ which appears as a summand in $p'_1,\ldots,p'_k$ produces $k$ chains of edges $p_1,\ldots,p_k$ in the original graph $G$, such that no two of $p_1,\ldots,p_k$ have any common nontrivial summands, satisfying $\partial(p_1)=\cdots=\partial(p_k)=v_1+v_2+v_3+v_4$ in the original $G$. Finally, by our induction hypothesis on the number of vertices in $Y$, we can find the required $k$ chains of edges $p'_1,\ldots,p'_k$ in the modified graph $G$. This concludes our Case 5, and concludes the proof of the lemma.
\end{proof}
The requirement that the degrees of all the vertices, except $v_1,v_2,v_3,v_4$, must always be even is instrumental. Without this requirement, Lemma \ref{minimal2} would be wrong. Consider the following counter-example.
\begin{exa}\label{counter}
Consider a graph $G$ with the vertex set $V$ and the edge set $E$ $$V=\{v_1,v_2,v_3,v_4,y_1,y_2\}\;\;,\;\;E=\{e_1,e_2,e_3,h_1,h_2,h_3,g\}$$ such that $$\partial(e_1)=\{v_1,y_1\}, \partial(e_2)=\{v_2,y_1\}, \partial(e_3)=\{v_3,y_1\},$$ $$\partial(h_1)=\{v_2,y_2\}, \partial(h_2)=\{v_3,y_2\}, \partial(h_3)=\{v_4,y_2\},$$ and $\partial(g)=\{v_1,v_4\}$. We use $\partial(x)=\{a,b\}$ and $\partial(x)=a+b$ interchangeably.\\
\includegraphics[scale=0.5]{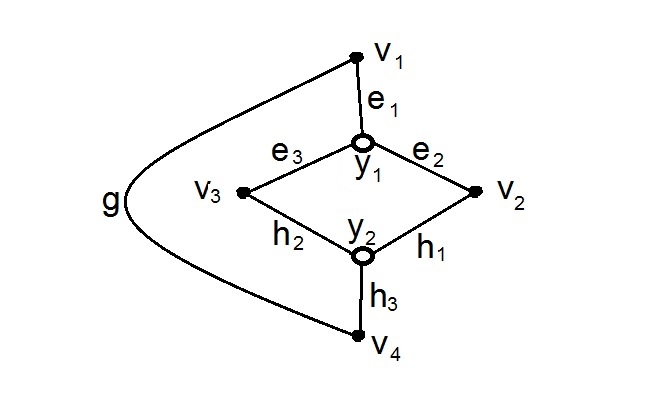}
\\
All the requirements of Lemma \ref{minimal2}, except for the requirement that the degrees of $y_1$ and $y_2$ must be even, are satisfied. Indeed, for every edge $x\in E$, at least one of the summands in $\partial(x)$ belongs to the set $\{v_1,v_2,v_{3},v_{4}\}$. The degree of each one of the vertices $v_1,v_2,v_{3},v_{4}$ is $k=2$. The number of edges $e$ for which $\partial(e)=y_1+v_1$ is not greater than the total number of all the edges $a$ for which $\partial(a)=y_1+v_2$ or $\partial(a)=y_1+v_3$ or $\partial(a)=y_1+v_4$. The number of edges $e$ for which $\partial(e)=y_1+v_2$ is not greater than the total number of all the edges $a$ for which $\partial(a)=y_1+v_1$ or $\partial(a)=y_1+v_3$ or $\partial(a)=y_1+v_4$. The number of edges $e$ for which $\partial(e)=y_1+v_3$ is not greater than the total number of all the edges $a$ for which $\partial(a)=y_1+v_1$ or $\partial(a)=y_1+v_2$ or $\partial(a)=y_1+v_4$. The number of edges $e$ for which $\partial(e)=y_1+v_4$ is not greater than the total number of all the edges $a$ for which $\partial(a)=y_1+v_1$ or $\partial(a)=y_1+v_2$ or $\partial(a)=y_1+v_3$. The same is also true for $y_2$.
\\ \\
After a deletion of any $k-1=1$ edge, the graph $G$ still consists of one connected component and hence there still exists some chain of edges $p$ such that $\partial(p)=v_1+v_2+v_3+v_4$. However, since the degree of $y_1$ is $3$ and the degree of $y_2$ is $3$, there do not exist $k=2$ chains of edges $p_1,p_2$ in $A(E)$, such that $p_1$ and $p_2$ do not have any common nontrivial summands, satisfying $\partial(p_1)=\partial(p_2)=v_1+v_2+v_3+v_4$. Indeed, any path between vertices $v_i$ and $v_j$ with $i\ne j$, except for the path $P=(v_1=w_0,e_1=g,w_1=v_4)$, must contain two or more edges in it. Thus, if such two chains of edges $p_1$ and $p_2$ would exist, then every edge in the set $E$, which consists of seven edges, would be a summand in $p_1$ or in $p_2$. Thus, $\partial(p_1+p_2)=y_1+y_2$ instead of being $0$.
\end{exa}
Let now $G$ be any undirected graph and let $v_1,v_2,v_3,v_4$ be any four different vertices in $G$. Let $y$ be a vertex different from $v_1,v_2,v_3,v_4$, such that there is an odd number of edges $\gamma$ in $G$ for which $\partial(\gamma)=y+w$ with $w\in\{v_1,v_2,v_3,v_4\}$, and such that the number of edges $e$ in $G$ for which $\partial(e)=y+v_1$ is not greater than the total number of all the edges $a$ in $G$ for which $\partial(a)=y+v_2$ or $\partial(a)=y+v_3$ or $\partial(a)=y+v_4$, and the number of edges $e$ in $G$ for which $\partial(e)=y+v_2$ is not greater than the total number of all the edges $a$ in $G$ for which $\partial(a)=y+v_1$ or $\partial(a)=y+v_3$ or $\partial(a)=y+v_4$, and the number of edges $e$ in $G$ for which $\partial(e)=y+v_3$ is not greater than the total number of all the edges $a$ in $G$ for which $\partial(a)=y+v_1$ or $\partial(a)=y+v_2$ or $\partial(a)=y+v_4$, and the number of edges $e$ in $G$ for which $\partial(e)=y+v_4$ is not greater than the total number of all the edges $a$ in $G$ for which $\partial(a)=y+v_1$ or $\partial(a)=y+v_2$ or $\partial(a)=y+v_3$.
\\ \\
Our next lemma asserts, that if we create a new graph $G'$ either by deleting from $G$ an edge $\beta$ such that $\partial(\beta)=y+w$ with $w\in\{v_1,v_2,v_3,v_4\}$, or by drawing in $G$ a new edge $\alpha$ such that $\partial(\alpha)=y+w$ with $w\in\{v_1,v_2,v_3,v_4\}$, then the number of edges $e$ in $G'$ for which $\partial(e)=y+v_1$ is not greater than the total number of all the edges $a$ in $G'$ for which $\partial(a)=y+v_2$ or $\partial(a)=y+v_3$ or $\partial(a)=y+v_4$, and the number of edges $e$ in $G'$ for which $\partial(e)=y+v_2$ is not greater than the total number of all the edges $a$ in $G'$ for which $\partial(a)=y+v_1$ or $\partial(a)=y+v_3$ or $\partial(a)=y+v_4$, and the number of edges $e$ in $G'$ for which $\partial(e)=y+v_3$ is not greater than the total number of all the edges $a$ in $G'$ for which $\partial(a)=y+v_1$ or $\partial(a)=y+v_2$ or $\partial(a)=y+v_4$, and the number of edges $e$ in $G'$ for which $\partial(e)=y+v_4$ is not greater than the total number of all the edges $a$ in $G'$ for which $\partial(a)=y+v_1$ or $\partial(a)=y+v_2$ or $\partial(a)=y+v_3$.
\begin{lem}\label{odd.vertex} Let $G$ be a graph in which any four vertices $v_1,v_2,v_3,v_4$ were selected. Let $y$ be a vertex in $G$ different from $v_1,v_2,v_3,v_4$, such that there is an odd number of edges $\gamma$ in $G$ for which $\partial(\gamma)=y+w$ with $w\in\{v_1,v_2,v_3,v_4\}$. If for all $i=1,2,3,4$, the number of edges $e$ in $G$ for which $\partial(e)=y+v_i$ is not greater than the total number of all the edges $a$ in $G'$ for which $\partial(a)=y+w$ with $w\in\{v_1,v_2,v_3,v_4\}-\{v_i\}$, then if we create a new graph $G'$ either by deleting from $G$ an edge $\beta$ such that $\partial(\beta)=y+w$ with $w\in\{v_1,v_2,v_3,v_4\}$, or by drawing in $G$ a new edge $\alpha$ such that $\partial(\alpha)=y+w$ with $w\in\{v_1,v_2,v_3,v_4\}$, then for all $i=1,2,3,4$, the number of edges $e$ in $G'$ for which $\partial(e)=y+v_i$ is not greater than the total number of all the edges $a$ in $G'$ for which $\partial(a)=y+w$ with $w\in\{v_1,v_2,v_3,v_4\}-\{v_i\}$.
\end{lem}
\begin{proof} Let the vertices $$v(G,y)_1,v(G,y)_2,v(G,y)_3,v(G,y)_4$$ be the four vertices $v_1,v_2,v_3,v_4$ indexed relative to the vertex $y$ and the graph $G$ in any such a way, that the number $\vartheta(G,y,1)$ of the edges in $G$, which have the boundary $y+v(G,y)_1$, is greater than or equal to the number $\vartheta(G,y,2)$ of the edges in $G$, which have the boundary $y+v(G,y)_2$, and the number $\vartheta(G,y,2)$ of the edges in $G$, which have the boundary $y+v(G,y)_2$, is greater than or equal to the number $\vartheta(G,y,3)$ of the edges in $G$, which have the boundary $y+v(G,y)_3$, and the number $\vartheta(G,y,3)$ of the edges in $G$, which have the boundary $y+v(G,y)_3$, is greater than or equal to the number $\vartheta(G,y,4)$ of the edges in $G$, which have the boundary $y+v(G,y)_4$.
\\
\includegraphics[scale=0.5]{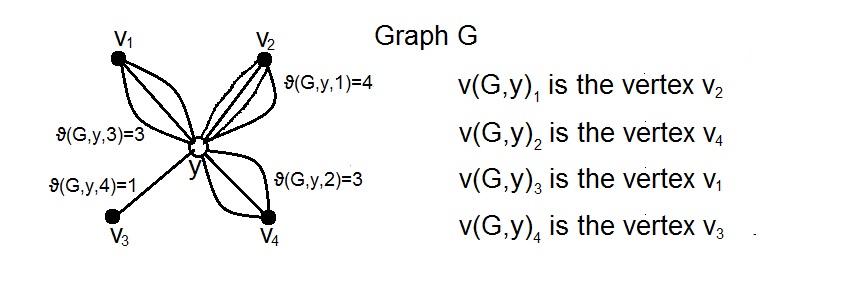}
\\
We have the following inequalities:
$$\vartheta(G,y,1)\ge\vartheta(G,y,2)\ge\vartheta(G,y,3)\ge\vartheta(G,y,4)\ge0$$
$$\vartheta(G,y,1)\le\vartheta(G,y,2)+\vartheta(G,y,3)+\vartheta(G,y,4)$$
Since there is an odd number of edges $\gamma$ in $G$ for which $\partial(\gamma)=y+w$ with $w\in\{v_1,v_2,v_3,v_4\}$, the second inequality must be strict
$$\vartheta(G,y,1)<\vartheta(G,y,2)+\vartheta(G,y,3)+\vartheta(G,y,4)$$
and we also get $\vartheta(G,y,3)\ge 1$, since if $\vartheta(G,y,3)=0$ then we must have $\vartheta(G,y,4)=0$ and $\vartheta(G,y,1)=\vartheta(G,y,2)$.
\\ \\
When we create a new graph $G'$ by deleting from $G$ an edge $\beta$ such that $\partial(\beta)=y+w$ with $w\in\{v_1,v_2,v_3,v_4\}$, then if $w$ is one of vertices $v(G,y)_2,v(G,y)_3,v(G,y)_4$ then we can index the four vertices $v_1,v_2,v_3,v_4$ relative to the vertex $y$ and the graph $G'$ in such a way, that the vertex $v(G',y)_1$ is the vertex $v(G,y)_1$. In that case we get
$$\vartheta(G',y,1)=\vartheta(G,y,1)\le\vartheta(G,y,2)+\vartheta(G,y,3)+\vartheta(G,y,4)-1=$$ $$=\vartheta(G',y,2)+\vartheta(G',y,3)+\vartheta(G',y,4)$$
If $w$ is the vertex $v(G,y)_1$ and $\vartheta(G,y,1)=\vartheta(G,y,2)$ then we can index the four vertices $v_1,v_2,v_3,v_4$ relative to the vertex $y$ and the graph $G'$ in such a way, that the vertex $v(G',y)_1$ is the vertex $v(G,y)_2$. In that case we get
$$\vartheta(G',y,1)=\vartheta(G,y,2)=\vartheta(G,y,1)\le\vartheta(G,y,2)+\vartheta(G,y,3)+\vartheta(G,y,4)-1=$$ $$=\vartheta(G,y,1)+\vartheta(G,y,3)+\vartheta(G,y,4)-1=\vartheta(G',y,2)+\vartheta(G',y,3)+\vartheta(G',y,4)$$
If $w$ is the vertex $v(G,y)_1$ and $\vartheta(G,y,1)>\vartheta(G,y,2)$ then we can index the four vertices $v_1,v_2,v_3,v_4$ relative to the vertex $y$ and the graph $G'$ in such a way, that the vertex $v(G',y)_1$ is the vertex $v(G,y)_1$. In that case we get
$$\vartheta(G',y,1)=\vartheta(G,y,1)-1<\vartheta(G,y,1)<\vartheta(G,y,2)+\vartheta(G,y,3)+\vartheta(G,y,4)=$$ $$=\vartheta(G',y,2)+\vartheta(G',y,3)+\vartheta(G',y,4)$$\\
When we create a new graph $G'$ by drawing in $G$ a new edge $\alpha$ such that $\partial(\alpha)=y+w$ with $w\in\{v_1,v_2,v_3,v_4\}$, then if $w$ is $v(G,y)_r$ with the index $r$ being equal to $2$ or $3$ or $4$, then if $\vartheta(G,y,r)<\vartheta(G,y,1)$ then we can index the four vertices $v_1,v_2,v_3,v_4$ relative to the vertex $y$ and the graph $G'$ in such a way, that the vertex $v(G',y)_1$ is the vertex $v(G,y)_1$. In that case we get
$$\vartheta(G',y,1)=\vartheta(G,y,1)<\vartheta(G,y,2)+\vartheta(G,y,3)+\vartheta(G,y,4)<$$ $$<\vartheta(G,y,2)+\vartheta(G,y,3)+\vartheta(G,y,4)+1=\vartheta(G',y,2)+\vartheta(G',y,3)+\vartheta(G',y,4)$$
If $w$ is $v(G,y)_r$ with the index $r$ being equal to $2$ or $3$ or $4$, and $\vartheta(G,y,r)=\vartheta(G,y,1)$ then we can index the four vertices $v_1,v_2,v_3,v_4$ relative to the vertex $y$ and the graph $G'$ in such a way, that the vertex $v(G',y)_1$ is the vertex $v(G,y)_r$. In that case we get
$$\vartheta(G',y,1)=\vartheta(G,y,r)+1=\vartheta(G,y,1)+1\le\vartheta(G,y,2)+\vartheta(G,y,3)+\vartheta(G,y,4)=$$ $$=\sum\limits_{t\in \{1,2,3,4\}-\{r\}}\vartheta(G,y,t)=\vartheta(G',y,2)+\vartheta(G',y,3)+\vartheta(G',y,4)$$
If $w$ is the vertex $v(G,y)_1$ then we can index the four vertices $v_1,v_2,v_3,v_4$ relative to the vertex $y$ and the graph $G'$ in such a way, that the vertex $v(G',y)_1$ is the vertex $v(G,y)_1$. In that case we get
$$\vartheta(G',y,1)=\vartheta(G,y,1)+1\le\vartheta(G,y,2)+\vartheta(G,y,3)+\vartheta(G,y,4)=$$ $$=\vartheta(G',y,2)+\vartheta(G',y,3)+\vartheta(G',y,4)$$
\end{proof}
Now we state and prove our fifth lemma. The idea of that lemma is as follows: Assume that in a graph $G$ all the edges go either between $v_i$ and $v_j$ or between $v_i$ and vertices from $Y$. Thus, the boundary of every edge contains among its summands one or two of the vertices $v_1,v_2,v_3,v_4$. Assume that after a deletion of any $k-1$ or less edges from $G$, there still exists some chain of edges $p$ such that $\partial(p)=v_1+v_2+v_3+v_4$. We look at all the vertices in $G$ which are different from the four vertices $v_1,v_2,v_3,v_4$, and mark-down all the vertices $y_1,\ldots,y_m$ with odd degrees in $G$. Next, we draw new edges in the graph $G$ in any such a way, that every vertex $y_1,\ldots,y_m$ appears as a summand in the boundary of exactly one new edge. Each one of these new edges we can choose to draw either between two of the vertices $y_1,\ldots,y_m$ or between one of the vertices $y_1,\ldots,y_m$ and one of the vertices $v_1,v_2,v_3,v_4$. After we are done drawing all these new edges, the degree of each one of the vertices $y_1,\ldots,y_m$ increases exactly by one, and becomes even.
\\ \\
The lemma claims that regardless of how we draw these new edges, after we draw them there will be $k$ chains of edges $p_1,\ldots,p_k$ in the modified graph $G$, such that no two of $p_1,\ldots,p_k$ have any common nontrivial summands and that $\partial(p_1)=\cdots=\partial(p_k)=v_1+v_2+v_{3}+v_{4}$ in the modified $G$.
\\ \\
For example, in our Example \ref{counter} one can draw a new edge $\epsilon_1$ with boundary $\partial(\epsilon_1)=\{y_1,y_2\}$ or draw two new edges $\epsilon_1,\epsilon_2$ with boundaries $\partial(\epsilon_1)=\{y_1,v_i\}$ and $\partial(\epsilon_2)=\{y_2,v_j\}$, where $i$ and $j$ can be any numbers between $1$ and $4$, and in the new graph one will find $k=2$ chains of edges $p_1,p_2$ such that they do not have any common summands and that $\partial(p_1)=\partial(p_2)=v_1+v_2+v_{3}+v_{4}$. The cases when one draws two edges $\epsilon_1,\epsilon_2$ with boundaries $\partial(\epsilon_1)=\{y_1,v_i\}$ and $\partial(\epsilon_2)=\{y_2,v_j\}$ are covered by Lemma \ref{minimal2}. If one drew one new edge $\epsilon_1$ with boundary $\partial(\epsilon_1)=\{y_1,y_2\}$ then we can take $p_1=g+e_2+e_3$ and $p_2=e_1+\epsilon_1+h_3+h_1+h_2$, or take $p_1=g+h_1+h_2$ and $p_2=e_2+\epsilon_1+h_3+e_1+e_3$, or take $p_1=g+e_2+\epsilon_1+h_2$ and $p_2=e_1+e_3+h_1+h_3$.\\
\includegraphics[scale=0.5]{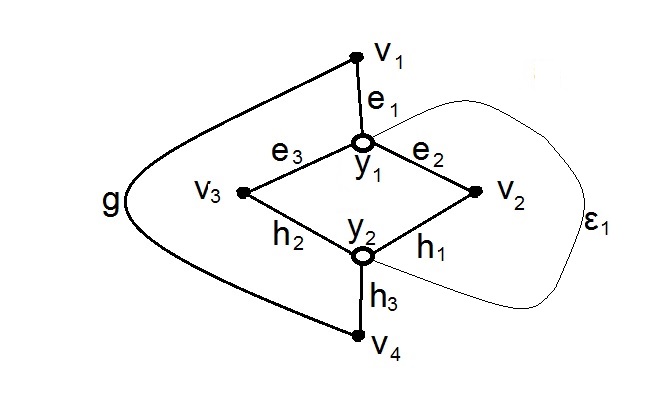}
\begin{lem}\label{minimal.final} Let $G$ be a graph with the vertex set $V$, which is a union of disjoint sets $\{v_1,v_2,v_{3},v_{4}\}$ and some empty or nonempty set $Y$, and the edge set $E$ such that for every $e\in E$, at least one of the summands in $\partial(e)$ belongs to the set $\{v_1,v_2,v_{3},v_{4}\}$. Let $y_1,\ldots,y_m$ be all the vertices in $Y$ such that their degrees in $G$ are odd. Here $m$ is $0$ if no $y\in Y$ has an odd degree in $G$.
\\ \\
If after a deletion of any $k-1$ or less edges from $G$ there still remains some chain of edges $p$ such that $\partial(p)=v_1+v_2+v_{3}+v_{4}$ then if one creates a new graph $GN$ by drawing in $G$ new edges $\epsilon_1,\ldots,\epsilon_n$ in any such a way that the boundary $\partial(\epsilon_i)$ of each new edge $\epsilon_i$ is a sum of two different vertices, one of them from the set $\{y_1,\ldots,y_m\}$ and the other one from the set $\{v_1,v_2,v_3,v_4,y_1,\ldots,y_m\}$, and that each vertex $y_1,\ldots,y_m$ is a summand in the boundary of exactly one of the new edges $\epsilon_1,\ldots,\epsilon_n$, then in the graph $GN$ there exist $k$ chains of edges $p_1,\ldots,p_k$, such that no two of $p_1,\ldots,p_k$ have any common nontrivial summands and that $\partial(p_1)=\cdots=\partial(p_k)=v_1+v_2+v_{3}+v_{4}$ in $GN$. Here $n$ is $0$ if $m=0$, which is the case when no new edges can be drawn and $GN=G$.
\end{lem}
\begin{proof} For each edge $y\in Y$ let the vertices $$v(G,y)_1,v(G,y)_2,v(G,y)_3,v(G,y)_4$$ denote the four vertices $v_1,v_2,v_3,v_4$ indexed relative to the vertex $y$ and the graph $G$ in any such a way, that the number $\vartheta(G,y,1)$ of the edges in $G$, which have the boundary $y+v(G,y)_1$, is greater than or equal to the number $\vartheta(G,y,2)$ of the edges in $G$, which have the boundary $y+v(G,y)_2$, and the number $\vartheta(G,y,2)$ of the edges in $G$, which have the boundary $y+v(G,y)_2$, is greater than or equal to the number $\vartheta(G,y,3)$ of the edges in $G$, which have the boundary $y+v(G,y)_3$, and the number $\vartheta(G,y,3)$ of the edges in $G$, which have the boundary $y+v(G,y)_3$, is greater than or equal to the number $\vartheta(G,y,4)$ of the edges in $G$, which have the boundary $y+v(G,y)_4$.
\\ \\
Assume that after a deletion of any $k-1$ or less edges from $G$ there still exists some chain of edges $p$ in $G$ such that $\partial(p)=v_1+v_2+v_{3}+v_{4}$. We claim that if for some vertex $y\in Y$ we have $$\vartheta(G,y,1)>\vartheta(G,y,2)+\vartheta(G,y,3)+\vartheta(G,y,4)$$ then deleting from $G$ any $r=\vartheta(G,y,1)-[\vartheta(G,y,2)+\vartheta(G,y,3)+\vartheta(G,y,4)]$ edges $\gamma_1,\ldots,\gamma_r$ such that $\partial(\gamma_1)=\cdots=\partial(\gamma_r)=y+v(G,y)_1$, produces a new graph $G'$ such that after a deletion of any $k-1$ or less edges from $G'$ there still exists some chain of edges $p'$ in $G'$, such that $\partial(p')=v_1+v_2+v_{3}+v_{4}$.
\\
\includegraphics[scale=0.5]{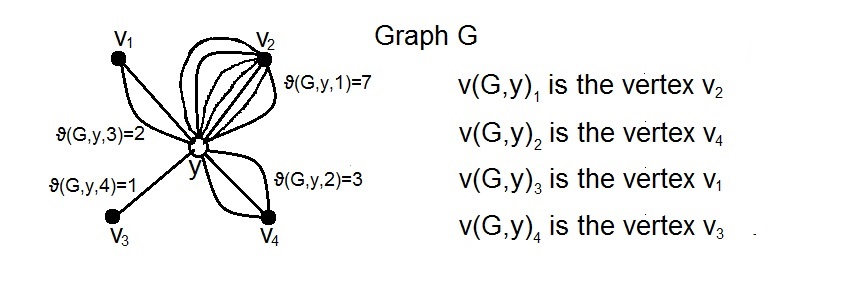}
\\
Indeed, if a after deletion of some $k-1$ or less edges from $G'$, one of the vertices $v_1,v_2,v_3,v_4$ becomes disconnected in $G'$ from the other three vertices, then after that deletion the vertex $y$ will either be in the same connected component as the vertex $\vartheta(G,y,1)$ or in a different connected component. If $y$ is in the same connected component as $\vartheta(G,y,1)$ then we can assume that there is no edge $\gamma$ such that $\partial(\gamma)=y+v(G,y)_1$ among the $k-1$ or less edges deleted from $G'$. Thus, deleting the same $k-1$ or less edges from $G$ will disconnect one of the vertices $v_1,v_2,v_3,v_4$ in $G$ from the other three vertices. If $y$ is in the different connected component than $\vartheta(G,y,1)$ then all the $\vartheta(G,y,2)+\vartheta(G,y,3)+\vartheta(G,y,4)$ edges in $G'$ whose boundaries are equal to $y+v(G,y)_1$ are among the $k-1$ or less edges deleted from $G'$. That implies, that selecting the same $k-1$ or less edges in $G$, and then un-selecting these $\vartheta(G,y,2)+\vartheta(G,y,3)+\vartheta(G,y,4)$ edges whose boundaries in $G'$ are $y+v(G,y)_1$, and instead of them selecting all the $\vartheta(G,y,2)+\vartheta(G,y,3)+\vartheta(G,y,4)$ or less edges whose boundaries in $G$ are either $y+v(G,y)_2$ or $y+v(G,y)_3$ or $y+v(G,y)_4$ and which were not already selected, and deleting these $k-1$ or less selected edges from $G$ will disconnect one of the vertices $v_1,v_2,v_3,v_4$ in $G$ from the other three vertices.
\\ \\
However, such a deletion of $\vartheta(G,y,1)-[\vartheta(G,y,2)+\vartheta(G,y,3)+\vartheta(G,y,4)]$ edges from $G$ can change the parity of the degree of the vertex $y$, which, in turn, tampers with our requirement to draw a new edge whose boundary contains $y$ as a summand in it. Hence, we have to be more careful. We proceed as follows:
\\
\includegraphics[scale=0.5]{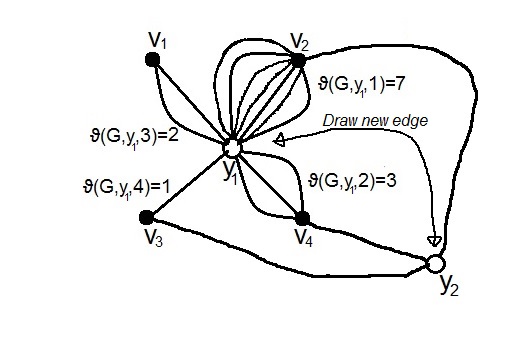}
\\
If for some vertex $y\in Y$ with an even degree in $G$ we have $$\vartheta(G,y,1)>\vartheta(G,y,2)+\vartheta(G,y,3)+\vartheta(G,y,4)$$ then we delete from $G$ any $$\vartheta(G,y,1)-[\vartheta(G,y,2)+\vartheta(G,y,3)+\vartheta(G,y,4)]$$ edges whose boundary is $y+v(G,y)_1$. If for some $y\in Y$ with an odd degree in $G$ we have $$\vartheta(G,y,1)>\vartheta(G,y,2)+\vartheta(G,y,3)+\vartheta(G,y,4)$$ then we delete from $G$ any $$\vartheta(G,y,1)-1-[\vartheta(G,y,2)+\vartheta(G,y,3)+\vartheta(G,y,4)]$$ edges, whose boundary is $y+v(G,y)_1$. Hence, we can assume that $$\vartheta(G,y,1)\le\vartheta(G,y,2)+\vartheta(G,y,3)+\vartheta(G,y,4)$$ for every vertex $y\in Y$ of an even degree in $G$ and $$\vartheta(G,y,1)\le 1+\vartheta(G,y,2)+\vartheta(G,y,3)+\vartheta(G,y,4)$$ for every vertex $y\in Y$ of an odd degree in $G$. Now, if $$\vartheta(G,y,1)=\vartheta(G,y,2)+\vartheta(G,y,3)+\vartheta(G,y,4)$$ then the degree of $y$ in $G$ is even. Thus, for all the vertices $y_1,\ldots,y_m$ we will have either $$\vartheta(G,y_i,1)=1+\vartheta(G,y_i,2)+\vartheta(G,y_i,3)+\vartheta(G,y_i,4)$$ or $$\vartheta(G,y_i,1)<\vartheta(G,y_i,2)+\vartheta(G,y_i,3)+\vartheta(G,y_i,4)$$\\
Assume that we have some vertices $y\in \{y_1,\ldots,y_m\}$ for which
$$\vartheta(G,y,1)>\vartheta(G,y,2)+\vartheta(G,y,3)+\vartheta(G,y,4)$$ We are going to show that for any given construction of a new graph $GN$ by drawing the new edges $\epsilon_1,\ldots,\epsilon_n$ in $G$ in any such a way as required by our lemma, we can reduce the graph $G$, by deleting some edges from $G$, to some graph $G'$ such that $G'$ satisfies the requirements of our lemma, and such that the number of vertices $y\in Y$ for which
$$\vartheta(G',y,1)>\vartheta(G',y,2)+\vartheta(G',y,3)+\vartheta(G',y,4)$$ is by one or two smaller than the number of vertices $y\in Y$ for which
$$\vartheta(G,y,1)>\vartheta(G,y,2)+\vartheta(G,y,3)+\vartheta(G,y,4)$$
and such that if our lemma is correct for the graph $G'$ then we can find $k$ chains of edges $p_1,\ldots,p_k$ in the graph $GN$, which was given to us, such that no two of $p_1,\ldots,p_k$ have any common nontrivial summands and that $\partial(p_1)=\cdots=\partial(p_k)=v_1+v_2+v_{3}+v_{4}$ in $GN$. Thus, assume that someone constructed the graph $GN$ by drawing some new edges $\epsilon_1,\ldots,\epsilon_n$ in $G$ in any such a way as required by our lemma, and gave this $GN$ to us.
\\ \\
We take any vertex $y_i$ for which $$\vartheta(G,y_i,1)=1+\vartheta(G,y_i,2)+\vartheta(G,y_i,3)+\vartheta(G,y_i,4)$$ Exactly one of the $n$ new edges $\epsilon_1,\ldots,\epsilon_n$ which were drawn in $G$ to create $GN$ has $y_i$ appearing as a summand in its boundary. We can assume that $\partial(\epsilon_1)=y_i+w$, where $w$ can be one of the vertices $v_1,v_2,v_3,v_4$ or can be one of the vertices $y_1,\ldots,y_{i-1},y_{i+1},\ldots,y_m$. We consider three possibilities for $w$.
\\ \\
\textit{Case 1:} The vertex $w$ is one of the vertices $v_1,v_2,v_3,v_4$. In that case we create our graph $G'$ by deleting from $G$ any one edge $\gamma$ such that $\partial(\gamma)=y_i+v(G,y_i)_1$. After a deletion of any $k-1$ or less edges from $G'$, we can still find some chain of edges $p'$ in $G'$ such that $\partial(p')=v_1+v_2+v_{3}+v_{4}$ in $G'$. Next we create the graph $GN'$ by drawing in $G'$ the $n-1$ new edges $\epsilon_2,\ldots,\epsilon_n$. The degree of the vertex $v_i$ in $G'$ is even, so this construction of $GN'$ from $G'$ satisfies the requirements of our lemma. The graph $GN'$ can be obtained from the graph $GN$ by deleting from the graph $GN$ the edge $\epsilon_1$ and the edge $\gamma$. Thus, if in the graph $GN'$ there exist $k$ chains of edges $p_1,\ldots,p_k$, such that no two of $p_1,\ldots,p_k$ have any common nontrivial summands and that $\partial(p_1)=\cdots=\partial(p_k)=v_1+v_2+v_{3}+v_{4}$ in $GN'$, then these same $k$ chains of edges $p_1,\ldots,p_k$ are chains of edges in the graph $GN$, which satisfy $\partial(p_1)=\cdots=\partial(p_k)=v_1+v_2+v_{3}+v_{4}$ in $GN$.
\\ \\
\textit{Case 2:} The vertex $w=y_j$, $j\ne i$, is one of the vertices $y_1,\ldots,y_{i-1},y_{i+1},\ldots,y_m$. In this case we have two subcases:
\\ \\
\textit{Subcase 1} is when $\vartheta(G,y_j,1)=1+\vartheta(G,y_j,2)+\vartheta(G,y_j,3)+\vartheta(G,y_j,4)$. In that subcase we create our graph $G'$ by deleting from $G$ any one edge $\gamma$ such that $\partial(\gamma)=y_i+v(G,y_i)_1$ and any one edge $\alpha$ such that $\partial(\alpha)=y_j+v(G,y_j)_1$. After a deletion of any $k-1$ or less edges from $G'$, we can still find some chain of edges $p'$ in $G'$ such that $\partial(p')=v_1+v_2+v_{3}+v_{4}$ in $G'$. Next we create the graph $GN'$ by drawing in $G'$ the $n-1$ new edges $\epsilon_2,\ldots,\epsilon_n$. The degrees of the vertices $v_i$ and $v_j$ in $G'$ are even, so this construction of $GN'$ from $G'$ satisfies the requirements of our lemma. The graph $GN'$ can be obtained from the graph $GN$ by deleting from the graph $GN$ the edge $\epsilon_1$ and the edges $\gamma$ and $\alpha$. Thus, if in the graph $GN'$ there exist $k$ chains of edges $p_1,\ldots,p_k$, such that no two of $p_1,\ldots,p_k$ have any common nontrivial summands and that $\partial(p_1)=\cdots=\partial(p_k)=v_1+v_2+v_{3}+v_{4}$ in $GN'$, then these same $k$ chains of edges $p_1,\ldots,p_k$ are chains of edges in the graph $GN$, which satisfy $\partial(p_1)=\cdots=\partial(p_k)=v_1+v_2+v_{3}+v_{4}$ in $GN$;
\\ \\
\textit{Subcase 2} is when $\vartheta(G,y_j,1)<\vartheta(G,y_j,2)+\vartheta(G,y_j,3)+\vartheta(G,y_j,4)$. In that subcase we create our graph $G'$ by deleting from $G$ any one edge $\gamma$ such that $\partial(\gamma)=y_i+v(G,y_i)_1$. After a deletion of any $k-1$ or less edges from $G'$, we can still find some chain of edges $p'$ in $G'$ such that $\partial(p')=v_1+v_2+v_{3}+v_{4}$ in $G'$. Next we create the graph $GN'$ by drawing in $G'$ the $n$ new edges $\epsilon'_1,\epsilon_2,\ldots,\epsilon_n$, where the boundary $\partial(\epsilon'_1)$ of the edge $\epsilon'_1$ is $v_j+v(G,y_i)_1$. Thus, $\partial(\epsilon'_1)$ in $GN'$ is the sum of the same vertices as $\partial(\epsilon_1+\gamma)$ is in $GN$. The degree of the vertex $v_i$ in $G'$ is even, so this construction of $GN'$ from $G'$ satisfies the requirements of our lemma. The graph $GN'$ can be obtained from the graph $GN$ by deleting from the graph $GN$ the edge $\epsilon_1$ and the edge $\gamma$, and then drawing the edge $\epsilon'_1$. Thus, if in the graph $GN'$ there exist $k$ chains of edges $p'_1,\ldots,p'_k$, such that no two of $p'_1,\ldots,p'_k$ have any common nontrivial summands and that $\partial(p'_1)=\cdots=\partial(p'_k)=v_1+v_2+v_{3}+v_{4}$ in $GN'$, then substituting in these $k$ chains of edges $p'_1,\ldots,p'_k$ for the edge $\epsilon'_{1}$, whenever $\epsilon'_{1}$ appears as a summand in one of these chains of edges, the sum $\epsilon_1+\gamma$ produces $k$ chains of edges $p_1,\ldots,p_k$ in the graph $GN$, with no two of them having any common nontrivial summands, which satisfy $\partial(p_1)=\cdots=\partial(p_k)=v_1+v_2+v_{3}+v_{4}$ in $GN$.
\\
\includegraphics[scale=0.5]{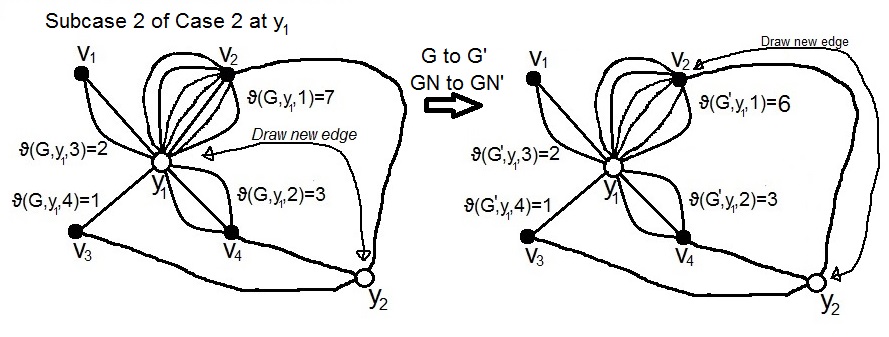}
\\
Repeating this reduction of the graph $G$ to the graph $G'$, and after each reduction substituting $G'$ for $G$ and $GN'$ for $GN$, will eventually produce a graph $G'$ in which $\vartheta(G',y,1)\le\vartheta(G',y,2)+\vartheta(G',y,3)+\vartheta(G',y,4)$ for every vertex $y\in Y$, such that this $G'$ satisfies the requirements of our lemma, and such that if our lemma is correct for that $G'$ then there exist $k$ chains of edges $p_1,\ldots,p_k$ in the graph $GN$, which was created from $G$ by drawing the new edges $\epsilon_1,\ldots,\epsilon_n$ in $G$ and given to us, with no two of these $k$ chains of edges having any common nontrivial summands, which satisfy $\partial(p_1)=\cdots=\partial(p_k)=v_1+v_2+v_{3}+v_{4}$ in $GN$. Thus, we only need to prove our lemma for graphs $G$ such that $\vartheta(G,y,1)\le\vartheta(G,y,2)+\vartheta(G,y,3)+\vartheta(G,y,4)$ for every vertex $y\in Y$.
\\ \\
Now, we start with some graph $G$ which satisfies the requirements of our lemma, such that $\vartheta(G,y,1)\le\vartheta(G,y,2)+\vartheta(G,y,3)+\vartheta(G,y,4)$ for every vertex $y\in Y$. Assume that someone constructed the graph $GN$ by drawing some new edges $\epsilon_1,\ldots,\epsilon_n$ in $G$ in any such a way as required by our lemma, and gave this graph $GN$ to us. The vertices $v_1,v_2,v_3,v_4$ in the graph $G$ have degrees not smaller than $k$. Hence, the degrees of these four vertices in the graph $GN$ are also not smaller than $k$. Every vertex $y\in Y$ has an even degree in $GN$. Also, every vertex $y\in Y$ which has an even degree in $G$ is not contained in the boundary of any new edge $\epsilon_1,\ldots,\epsilon_n$ in $GN$, hence at such $y$ we have $\vartheta(GN,y,1)\le\vartheta(GN,y,2)+\vartheta(GN,y,3)+\vartheta(GN,y,4)$. Additionally, for every vertex $y\in Y$ which has an odd degree in $G$, such that some new edge $\epsilon_i$ has in $GN$ the boundary $\partial(\epsilon_i)=y+w_i$ with $w_i\in \{v_1,v_2,v_3,v_4\}$, Lemma \ref{odd.vertex} asserts that $\vartheta(GN,y,1)\le\vartheta(GN,y,2)+\vartheta(GN,y,3)+\vartheta(GN,y,4)$.
\\ \\
However, we may have some $2z>0$ vertices in the set $\{y_1,\ldots,y_m\}$, we can assume that they are $y_1,\ldots,y_{2z}$, which in the graph $GN$ are connected between themselves by some $z>0$ new edges from the set $\{\epsilon_1,\ldots,\epsilon_n\}$, we can assume that these are the edges $\epsilon_1,\ldots,\epsilon_{z}$ for which $\partial(\epsilon_i)=y_{2i-1}+y_{2i}$, and this prevents us from applying Lemma \ref{minimal2} to the graph $GN$. We are going to construct a graph $GN'$ from the graph $GN$ by deleting the $z$ edges $\epsilon_1,\ldots,\epsilon_{z}$ and another $z$ edges which we will specify below, and then drawing new edges $\epsilon'_1,\ldots,\epsilon'_{z}$ with $\partial(\epsilon'_i)=u_i+y_{2i}$ where $u_i\in \{v_1,v_2,v_3,v_4\}$ for all $i=1,\ldots,z$.
\\ \\
We will perform this construction of $GN'$ in such a way, that the degrees of the vertices $v_1,v_2,v_3,v_4$ in the graph $GN'$ will remain the same as their degrees in the graph $GN$, and that every vertex $y\in Y$ will have an even degree in $GN'$, and that $$\vartheta(GN',y,1)\le\vartheta(GN',y,2)+\vartheta(GN',y,3)+\vartheta(GN',y,4)$$ will be satisfied at every vertex $y\in Y$. Thus, Lemma \ref{minimal2} will assert the existence in the graph $GN'$ of $k$ chains of edges $p'_1,\ldots,p'_k$, such that no two of $p'_1,\ldots,p'_k$ have any common nontrivial summands and that $\partial(p'_1)=\cdots=\partial(p'_k)=v_1+v_2+v_{3}+v_{4}$ in $GN'$. Finally, we will use these $k$ chains of edges $p'_1,\ldots,p'_k$ in $GN'$ to construct $k$ chains of edges $p_1,\ldots,p_k$ in the graph $GN$, with no two of $p_1,\ldots,p_k$ having any common nontrivial summands, which satisfy $\partial(p_1)=\cdots=\partial(p_k)=v_1+v_2+v_{3}+v_{4}$ in $GN$.
\\ \\
The graph $GN'$ is constructed from the graph $GN$ by deleting from $GN$ any $z$ edges $\gamma_1,\ldots,\gamma_z$ such that $\partial(\gamma_i)=\vartheta(G,y_{2i-1},1)+y_{2i-1}$, and deleting the $z$ edges $\epsilon_1,\ldots,\epsilon_{z}$, and then drawing $z$ new edges $\epsilon'_1,\ldots,\epsilon'_{z}$ such that the boundary $\partial(\epsilon'_i)$ of $\epsilon'_i$ in $GN'$ is $\vartheta(G,y_{2i-1},1)+y_{2i}$. Thus, we delete some edge $\gamma_1$ such that $\partial(\gamma_1)=\vartheta(G,y_{1},1)+y_{1}$, and we delete the edge $\epsilon_1$, and we draw an edge $\epsilon'_1$ such that $\partial(\epsilon'_1)=\vartheta(G,y_{1},1)+y_{2}$, and so on.
\\
\includegraphics[scale=0.6]{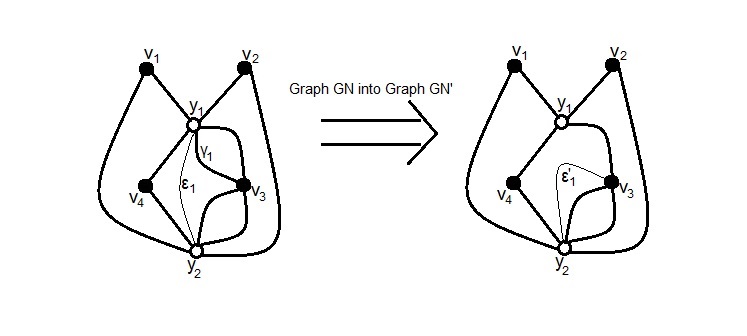}
\\
The boundary of every edge in the graph $GN'$ contains at least one of the four vertices $v_1,v_2,v_3,v_4$ among its two summands. At each vertex $y_{2i-1}$ the graph $GN'$ differs from the graph $G$ by one edge $\gamma_i$ which is deleted from $G$. At each vertex $y_{2i}$ the graph $GN'$ differs from the graph $G$ by one new edge $\epsilon'_1$ which is drawn in $G$. Thus, by Lemma \ref{odd.vertex}, we have $$\vartheta(GN',y_i,1)\le\vartheta(GN',y_i,2)+\vartheta(GN',y_i,3)+\vartheta(GN',y_i,4)$$
at all the $2z$ vertices $y_1,y_2\ldots,y_{2z-1},y_{2z}$ in the graph $GN'$. Thus, $$\vartheta(GN',y,1)\le\vartheta(GN',y,2)+\vartheta(GN',y,3)+\vartheta(GN',y,4)$$ at every vertex $y\in Y$ in $GN'$. The degrees of the vertices $y_1,y_3\ldots,y_{2z-1}$ in the the graph $GN'$ are by $1$ smaller than the degrees of these vertices in the graph $G$, and by $2$ smaller than the degrees of these vertices in the graph $GN$. The degrees of the vertices $y_2,y_4\ldots,y_{2z}$ in the the graph $GN'$ are the same as the degrees of these vertices in the graph $GN$. Thus, the degree of every vertex $y\in Y$ in $GN'$ is always even, and is equal to $$\vartheta(GN',y,1)+\vartheta(GN',y,2)+\vartheta(GN',y,3)+\vartheta(GN',y,4)$$ because every edge in $GN'$ contains at least one of the four vertices $v_1,v_2,v_3,v_4$ as a summand in its boundary. The degrees of the four vertices $v_1,v_2,v_3,v_4$ in the graph $GN'$ remained the same as the degrees of these four vertices in the graph $GN$, which implies that these degrees are not smaller than $k$. Lemma \ref{minimal2} now asserts the existence of $k$ chains of edges $p'_1,\ldots,p'_k$ in the graph $GN'$, such that no two of $p'_1,\ldots,p'_k$ have any common nontrivial summands, and that $\partial(p'_1)=\cdots=\partial(p'_k)=v_1+v_2+v_{3}+v_{4}$ in $GN'$.
\\ \\
Finally, we take these $k$ chains of edges $p'_1,\ldots,p'_k$, and for each edge $\epsilon'_i$, which we find as a summand in one of these $k$ chains of edges, we substitute $\epsilon_i+\gamma_i$. We obtain $k$ chains of edges $p_1,\ldots,p_k$ in the graph $GN$, such that no two of $p_1,\ldots,p_k$ have any common nontrivial summands and that $\partial(p'_1)=\cdots=\partial(p'_k)=v_1+v_2+v_{3}+v_{4}$ in $GN$. Thus, in any graph $GN$, which was constructed in any way as required by our lemma, and given to us, we can find $k$ chains of edges $p_1,\ldots,p_k$, with no two of $p_1,\ldots,p_k$ having any common nontrivial summands, which satisfy $\partial(p_1)=\cdots=\partial(p_k)=v_1+v_2+v_{3}+v_{4}$ in $GN$.
\end{proof}
Now we state and prove the edge version of the classical Menger's Theorem using the language of homological paths. We will use this theorem and the proof, which we provide here, in the proof of our Extended Menger's Edge Theorem.
\begin{thm}\label{Menger.homol.2} For any two vertices $v_1$ and $v_2$ in a graph $G$ the following two statements are equivalent: After a deletion of any $k-1$ or less edges from $G$ there still exists some chain of edges $p$ in $G$ such that $\partial(p)=v_1+v_2$; There exist $k$ chains of edges $p_1,\ldots,p_k$ in $G$, such that no two of $p_1,\ldots,p_k$ have any common nontrivial summands and that $\partial(p_1)=\cdots=\partial(p_k)=v_1+v_2$.
\end{thm}
\begin{proof}
For any $k$ chains of edges $p_1,\ldots,p_k$ in $G$ such that no two of $p_1,\ldots,p_k$ have any common nontrivial summands, a deletion of any $k-1$ or less edges from $G$ will not effect one or more of these $k$ chains of edges. Thus, if $\partial(p_1)=\cdots=\partial(p_k)=v_1+v_2$, then after such a deletion there exists some chain of edges $p$ in $G$ such that $\partial(p)=v_1+v_2$. Now we prove the other direction in the theorem.
\\ \\
We prove the other direction by induction on $k$. Our theorem is trivial for $k=1$. Assume that for some integer $k>1$, our theorem is correct for all the integers from $1$ to $k-1$. We will prove it for $k$. Assume that the theorem is wrong for $k$. Then there exist some counter-examples to it. Among these counter-examples we select a counter-example with the minimal number of edges in its graph $G$. Thus, in all the other counter-examples the graphs have the same number of edges as $G$ or more edges than $G$.
\\ \\
Now, for every edge $x_1$ in the graph $G$ there must exist some $k-1$ edges $x_2,\ldots,x_{k}$ in $G$ such that after the deletion of all the $k$ edges $x_1,x_2,\ldots,x_{k}$ from $G$ there does not exist any chain of edges $p$ in $G$ such that $\partial(p)=v_1+v_2$. Indeed, if for some edge $x_1$ in the graph $G$ no such $k-1$ edges $x_2,\ldots,x_{k}$ exist in $G$ then we can delete the edge $x_1$ from our graph $G$ and obtain a new counter-example to our theorem, in which the graph has less edges than $G$. This contradicts the minimality of the number of edges in $G$.
\\ \\
Assume that there is an edge $x_1$ in $G$ such that $\partial(x_1)=u_{1,1}+u_{1,2}$ does not contain the vertices $v_1$ and $v_2$ among its two summands $u_{1,1}$ and $u_{1,2}$. We take any $k-1$ edges $x_2,\ldots,x_{k}$ in $G$ such that after the deletion of all the $k$ edges $x_1,x_2,\ldots,x_{k}$ from the graph $G$, there remains no chain of edges $p$ in $G$ such that $\partial(p)=v_1+v_2$. For each $i=2,\ldots,k$ the boundary $\partial(x_i)$ is $u_{i,1}+u_{i,2}$, where $u_{i,1}$ and $u_{i,2}$ are two vertices in $G$. Let the graph $\Gamma$ be equal to the graph $G$ without the edges $x_1,x_2,\ldots,x_{k}$. The vertices $v_1$ and $v_2$ must belong to the different connected components of the graph $\Gamma$. Let $K_1$ be the connected component of the graph $\Gamma$ which contains the vertex $v_1$ and let $K_2$ be the connected component of the graph $\Gamma$ which contains the vertex $v_2$. We label the two vertices in the boundaries $\partial(x_2),\ldots,\partial(x_{k})$ in such a way, that all the vertices $u_{1,1},\ldots,u_{k,1}$ belong to $K_1$ and all the vertices $u_{1,2},\ldots,u_{k,2}$ belong to $K_2$.
\\
\includegraphics[scale=0.5]{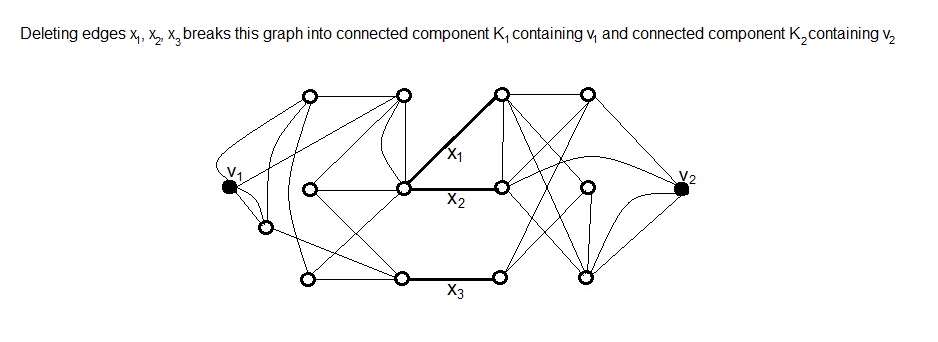}
\\
We construct a new graph $T$ from the graph $G$ by performing the following three steps:
\\
\textit{Step 1:} Draw a new vertex $v'_1$;
\\
\textit{Step 2:} For all the edges $x_i$, where $i=1,\ldots,k$, define the boundary operation in $T$ by $\partial(x_i)=u_{i,2}+v'_1$;
\\
\textit{Step 3:} Delete all the vertices and all the edges belonging to $K_1$.
\\
\includegraphics[scale=0.5]{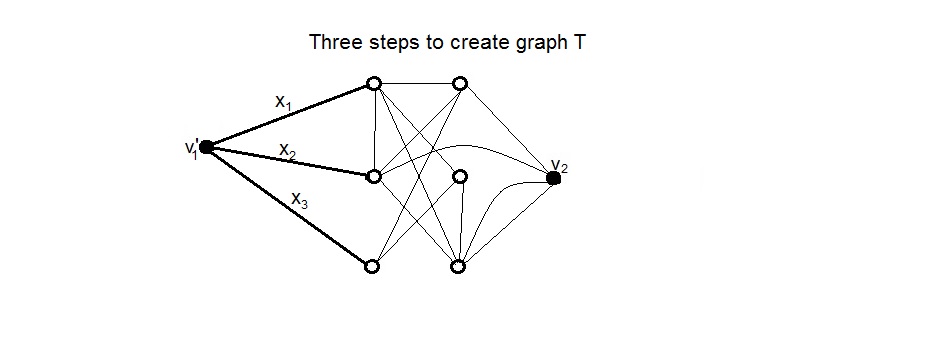}
\\
The boundary operation for all the other edges in $T$ remains the same as it was in $G$. Since $u_{1,1}\ne v_1$, the connected component $K_1$ of $\Gamma$ has at least one edge in it. Hence, the graph $T$ has less edges than the graph $G$, which implies that $T$ cannot be a counter-example to our theorem.  Now, if there exist some $k-1$ or less edges in $T$, such that after the deletion of these edges from $T$ there remains no chain of edges $p'$ in $T$ such that $\partial(p')=v'_1+v_2$ in $T$, then after the deletion of these same $k-1$ or less edges from $G$ there is no chain of edges $p$ in $G$ such that $\partial(p)=v_1+v_2$ in $G$. Indeed, if such a chain of edges $p$ would exist in $G$ then by removing from that $p$ all the edges, which appear as summands in $p$ and which belong to $K_1$, we would obtain a chain of edges $p'$ in $T$ such that $\partial(p')=v'_1+v_2$ in $T$. Thus, we can find $k$ chains of edges $p_{T,1},\ldots,p_{T,k}$ in $T$, such that no two of $p_{T,1},\ldots,p_{T,k}$ have any common nontrivial summands and that $\partial(p_{T,1})=\cdots=\partial(p_{T,k})=v'_1+v_2$ in $T$.
\\ \\
Since the only edges which contain the vertex $v'_1$ as a summand in their boundary in the graph $T$ are the $k$ edges $x_1,\ldots,x_k$, we see that each one of the $k$ chains of edges $p_{T,1},\ldots,p_{T,k}$ must contain as a summand in it one edge from among $x_1,\ldots,x_k$. Thus, we can assume that the chain of edges $p_{T,i}$, for all $i=1,\ldots,k$, contains in at as a summand the edge $x_i$.
\\
\includegraphics[scale=0.5]{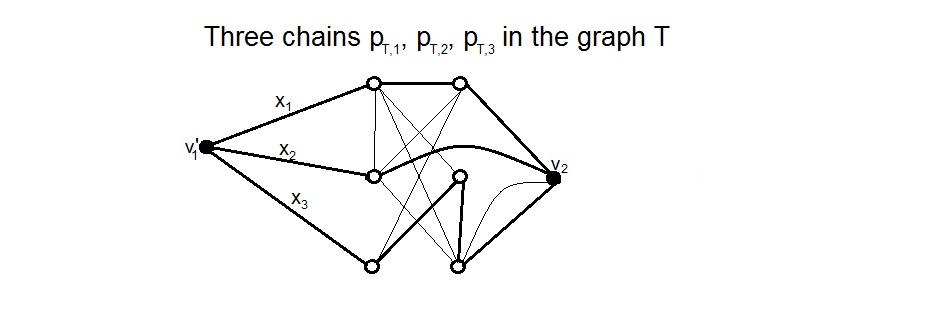}
\\
Next, we construct a new graph $L$ from the graph $G$ by performing the following three steps:
\\
\textit{Step 1:} Draw a new vertex $v'_2$;
\\
\textit{Step 2:} For all the edges $x_i$, where $i=1,\ldots,k$, define the boundary operation in $L$ by $\partial(x_i)=u_{i,1}+v'_2$;
\\
\textit{Step 3:} Delete all the vertices and all the edges belonging to $K_2$.
\\
\includegraphics[scale=0.5]{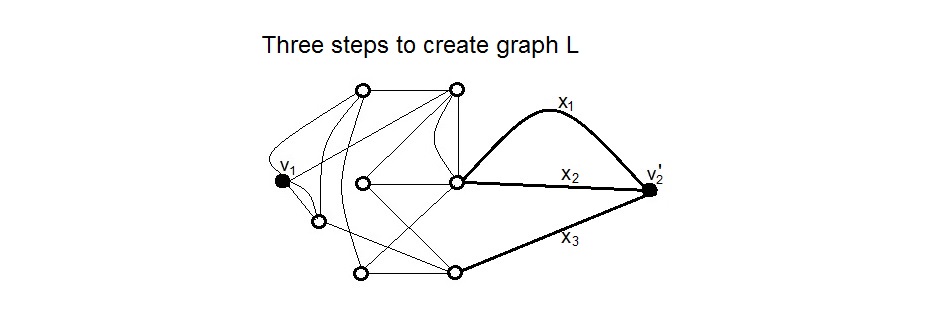}
\\
The boundary operation for all the other edges in $L$ remains the same as it was in $G$. Since $u_{1,2}\ne v_2$, the connected component $K_2$ of $\Gamma$ has at least one edge in it. Thus, the graph $L$ has less edges than the graph $G$, which implies that $L$ cannot be a counter-example to our theorem.  Now, if there exist some $k-1$ or less edges in $L$, such that after the deletion of these edges from $L$ there remains no chain of edges $p'$ in $L$ such that $\partial(p')=v_1+v'_2$ in $L$, then after the deletion of these same $k-1$ or less edges from $G$ there is no chain of edges $p$ in $G$ such that $\partial(p)=v_1+v_2$ in $G$. Indeed, if such a chain of edges $p$ would exist in $G$ then by removing from that $p$ all the edges, which appear as summands in $p$ and which belong to $K_2$, we would obtain a chain of edges $p'$ in $L$ such that $\partial(p')=v_1+v'_2$ in $L$. Thus, we can find $k$ chains of edges $p_{L,1},\ldots,p_{L,k}$ in $L$, such that no two of $p_{L,1},\ldots,p_{L,k}$ have any common nontrivial summands and that $\partial(p_{L,1})=\cdots=\partial(p_{L,k})=v_1+v'_2$ in $L$.
\\ \\
Since the only edges which contain the vertex $v'_2$ as a summand in their boundary in the graph $L$ are the $k$ edges $x_1,\ldots,x_k$, we see that each one of the $k$ chains of edges $p_{L,1},\ldots,p_{L,k}$ must contain as a summand in it one edge from among $x_1,\ldots,x_k$. Thus, we can assume that the chain of edges $p_{L,i}$, for all $i=1,\ldots,k$, contains in at as a summand the edge $x_i$.
\\
\includegraphics[scale=0.5]{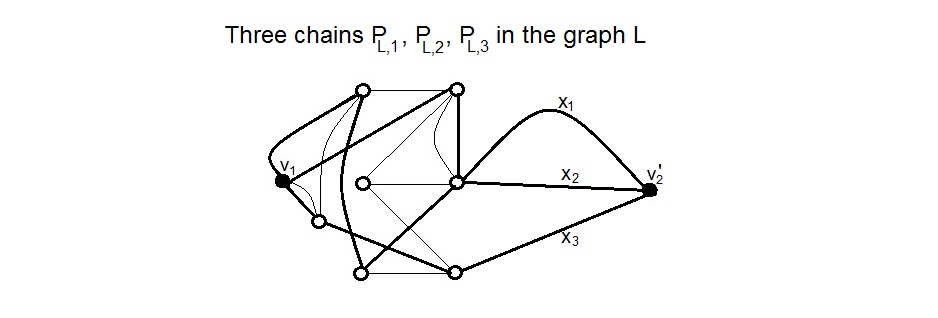}
\\
Finally, we construct $k$ chains of edges $p_{1},\ldots,p_{k}$ in $G$ by taking the chains of edges $p_{T,1},\ldots,p_{T,k}$, and in each chain of edges $p_{T,i}$ substituting the summand $x_i$ by the whole chain of edges $p_{L,i}$. The boundary operation in graph $G$ takes each $p_{T,i}$ to $u_{i,1}+v_2$ and takes each $p_{L,i}$ to $\partial(x_i)=v_1+u_{i,2}$. Since in $G$ the boundary of each $x_i$, which in $p_{T,i}$ was substituted by $p_{L,i}$, is $u_{i,1}+u_{i,2}$, we have $\partial(p_1)=\cdots=\partial(p_k)=v_1+v_2$ in the graph $G$. We see that $G$ is not a counter-example to our theorem, which implies that there is no edge $x_1$ in $G$ such that its boundary $\partial(x_1)=u_{1,1}+u_{1,2}$ does not contain either $v_1$ or $v_2$ or both of them among its two summands $u_{1,1}$ and $u_{1,2}$. Thus, for every edge $e$ in $G$ we must have $\partial(e)=v_1+w$ or $\partial(e)=v_2+w$ or $\partial(e)=v_1+v_2$.
\\ \\
Now, for any vertex $w$ in $G$, such that $w$ is different from $v_1$ and $v_2$, if the number $n_1$ of edges the $a$ such that $\partial(a)=v_1+w$ is greater than the number $n_2$ of the edges $b$ such that $\partial(b)=v_2+w$, then deleting from the graph $G$ any $(n_1-n_2)$ edges whose boundaries are $v_1+w$ does not effect the property that after a deletion of any $k-1$ or less edges from $G$ there still exists some chain of edges $p$ such that $\partial(p)=v_1+v_2$. In the other direction, if the number $n_2$ of the edges $b$ such that $\partial(b)=v_2+w$ is greater than the number $n_1$ of the edges $a$ such that $\partial(a)=v_1+w$, then deleting from the graph $G$ any $(n_2-n_1$) edges whose boundaries are $v_2+w$ does not effect the property that after a deletion of any $k-1$ or less edges from $G$ there still exists some chain of edges $p$ such that $\partial(p)=v_1+v_2$.
\\ \\
Thus, we can assume that for each vertex $w$ not belonging to the set $\{v_1,v_2\}$, the number of the edges $a$ such that $\partial(a)=v_1+w$ is equal to the number of the edges $b$ such that $\partial(b)=v_2+w$. This means that we can pair each edge $a$ such that $\partial(a)=v_1+w$ with some edge $b$ such that $\partial(b)=v_2+w$. We create a new graph $G'$ from $G$ by deleting all such pairs of edges $a,b$ for every vertex $w\notin\{v_1,v_2\}$, and for each deleted pair $a,b$ drawing one new edge $\gamma$ with $\partial(\gamma)=v_1+v_2$. Clearly, $\partial(\gamma)=\partial(a+b)$. The degrees of the vertices $v_1$ and $v_2$ in the new graph $G'$ are the same as the degrees of these two vertices in the original graph $G$, which implies that the degrees of the vertices $v_1$ and $v_2$ in the new graph $G'$ are greater than $k-1$. Thus, there exist $k$ chains of edges $p'_1,\ldots,p'_k$ in $G'$, such that no two of $p'_1,\ldots,p'_k$ have any common nontrivial summands and that $\partial(p'_1)=\cdots=\partial(p'_k)=v_1+v_2$ in $G'$. We take these $k$ chains of edges $p'_1,\ldots,p'_k$ and we substitute the sum of the pair $a+b$ for the corresponding edge $\gamma$ whenever we find $\gamma$ appearing as a summand in one of the $p'_1,\ldots,p'_k$. This produces $k$ chains of edges $p_1,\ldots,p_k$ in $G$, such that no two of $p_1,\ldots,p_k$ have any common nontrivial summands and that $\partial(p_1)=\cdots=\partial(p_k)=v_1+v_2$ in $G$. Thus, there cannot be a counter-example to our theorem.
\end{proof}
The structure of the proof of Theorem \ref{Menger.homol.2} was as follows: First we consider a hypothetical counter-example with the minimal number of edges in it. We prove that such a counter-example cannot have any edges not containing at least one of the two vertices $v_1,v_2$ in its boundary. This leaves us with the graphs, in which every edge has at least one of the two vertices $v_1,v_2$ in its boundary, and we reduce these graphs to the graphs, in which every edge has the boundary $v_1+v_2$.
\\ \\
The proof of our extension of Menger's Edge Theorem to the case of four vertices $v_1,v_2,v_3,v_4$ has a similar structure. We start by considering a hypothetical counter-example with the minimal number of edges in its graph $G$. We show, that if this example contains an edge $x_1$ whose boundary does not contain any of the vertices $v_1,v_2,v_3,v_4$, then we can break $G$ into two graphs $K_1$ and $K'$, each of them containing at least one edge, with $K_1$ containing only one of the four vertices $v_1,v_2,v_3,v_4$ and $K'$ containing the other three vertices, so that $K_1$ and $K'$ are connected by the edge $x_1$ and some $k-1$ edges $x_2,\ldots\,x_k$. This step is not trivial, because the deletion of the edge $x_1$ changes the parities of the degrees of the two vertices $u_{1,1}$ and $u_{1,2}$ which are the summands in the boundary of $x_1$, and that tampers with certain requirements in our theorem.
\\ \\
Next we show, that if we take the graph $K'$, add to it a new vertex $\psi$, and then draw the edges $x_1,\ldots\,x_k$, substituting in their boundaries the vertex $\psi$ for the vertices not belonging to $K'$, we will obtain a graph $T$ which satisfies the requirements of our theorem for the vertex $\psi$ and the three of the four vertices $v_1,v_2,v_3,v_4$ which are contained in $K'$. Thus, from the minimality of the number of edges in $G$ it follows, that our theorem will be true for $T$ relative to the vertex $\psi$ and the three of the four vertices $v_1,v_2,v_3,v_4$ which are contained in $K'$.
\\ \\
Next we show, that if the take the graph $K_1$, add to it a new vertex $\omega$, and then draw the edges $x_1,\ldots\,x_k$, substituting in their boundaries the vertex $\omega$ for the vertices not belonging to $K_1$, we will obtain a graph $L$ which satisfies the requirements of Theorem \ref{Menger.homol.2} for the vertex $\omega$ and the one of the four vertices $v_1,v_2,v_3,v_4$, which is contained in $K_1$.
\\ \\
Next we show, that whenever one produces a new graph $GN$ by drawing in $G$ any new edges $\epsilon_1,\ldots,\epsilon_n$ in any such a way as required by our theorem, we can construct $k$ chains of edges $p_1,\ldots,p_k$ in $GN$, such that no two of $p_1,\ldots,p_k$ have any common nontrivial summands and that their boundaries in $GN$ are all equal to $v_1+v_2+v_3+v_4$, in the following way: Use Theorem \ref{Menger.homol.2} to create $k$ chains of edges $p_{L,1},\ldots,p_{L,k}$ in $L$, such that no two of $p_{L,1},\ldots,p_{L,k}$ have any common nontrivial summands and that their boundaries in $L$ are all equal to the sum of the two vertices in $L$; Construct a new graph $TN$ by drawing new edges in the graph $T$ in a certain way, which satisfies the requirement of our theorem, and in the process of drawing these new edges in $T$ construct certain chains of edges $\gamma_{*}$ in $GN$; Find $k$ chains of edges $p_{T,1},\ldots,p_{T,k}$ in $TN$, such that no two of $p_{T,1},\ldots,p_{T,k}$ have any common nontrivial summands and that their boundaries in $TN$ are all equal to the sum of the four vertices in $TN$; In the chains of edges $p_{T,1},\ldots,p_{T,k}$ substitute each one of the new edges, which were drawn in $T$ to create $TN$, either by one of the edges $\epsilon_1,\ldots,\epsilon_n$ or by one of the chains of edges $\gamma_{*}$, and substitute each $x_i$ there by the chain of edges $p_{L,i}$.
\\ \\
Thus, we show, that our hypothetical counter-example with the minimal number of edges in it cannot have an edge not containing at least one of the four vertices $v_1,v_2,v_3,v_4$ as a summand in its boundary. This reduces our theorem to Lemma \ref{minimal.final}, which asserts that such a counter-example cannot exist.
\begin{thm}\label{Menger.homol.4} Let $G$ be a graph such that its vertex set $V$ is a union of disjoint sets $\{v_1,v_2,v_3,v_4\}$ and some empty or nonempty set $Y$. Let $y_1,\ldots,y_m$ be all the vertices in $Y$ such that their degrees in $G$ are odd. Here $m$ is $0$ if no such vertices exist in $Y$. If after a deletion of any $k-1$ or less edges from $G$ there still exists some chain of edges $p$ in $G$ such that $\partial(p)=v_1+v_2+v_{3}+v_{4}$, then if one creates a new graph $GN$ by drawing new edges $\epsilon_1,\ldots,\epsilon_n$ in $G$ in any such a way that the boundary $\partial(\epsilon_i)$ of each new edge $\epsilon_i$ is a sum of two different vertices, one of them from the set $\{y_1,\ldots,y_m\}$ and the other one from the set $\{v_1,v_2,v_3,v_4,y_1,\ldots,y_m\}$, and that each vertex $y_1,\ldots,y_m$ is a summand in the boundary of exactly one of the new edges $\epsilon_1,\ldots,\epsilon_n$, then in the graph $GN$ there exist $k$ chains of edges $p_1,\ldots,p_k$, such that no two of $p_1,\ldots,p_k$ have any common nontrivial summands and that $\partial(p_1)=\cdots=\partial(p_k)=v_1+v_2+v_{3}+v_{4}$ in $GN$. Here $n$ is $0$ if $m=0$, which is the case when no new edges can be drawn and $GN=G$.
\end{thm}
\begin{proof}
We prove our theorem by induction on $k$. For $k=1$ our theorem follows from Lemma \ref{tech}. Assume that for some integer $k>1$, our theorem is correct for all the integers from $1$ to $k-1$. We will prove it for $k$. Assume that the theorem is wrong for $k$. Then there exist some counter-examples to it. Among these counter-examples we select a counter-example with the minimal number of edges in its graph $G$. Thus, all the graphs in all the other counter-examples have the same number of edges as $G$ or more edges than $G$.
\\ \\
If the vertices $v_1,v_2,v_{3},v_{4}$ are split between two or more connected components of $G$ then two of these four vertices must belong to one of the connected components of $G$ and the other two vertices must belong to some other connected component of $G$, because if one of these four vertices belongs to its own connected component then there cannot exist any chain of edges $p$ in $G$ such that $\partial(p)=v_1+v_2+v_{3}+v_{4}$. Deleting any $k-1$ or less edges from any one of the two connected components which contain the vertices $v_1,v_2,v_{3},v_{4}$ does not disconnect the two vertices which belong to that connected component. The theorem now follows from Theorem \ref{Menger.homol.2}. Thus, we can assume that all four vertices $v_1,v_2,v_{3},v_{4}$ belong to one connected component of $G$.
\\ \\
From the minimality of the number of edges in $G$ it follows, that the connected components, which do not contain the vertices $v_1,v_2,v_{3},v_{4}$, cannot contain any edges in them. Thus, they all consist of single vertices, which implies that they can be deleted from the graph $G$ without ruining the counter-example. Thus, we assume that $G$ consists of only one connected component.
\\ \\
Since $G$ with the four vertices $v_1,v_2,v_3,v_4$ is a counter-example to our theorem, after a deletion of any $k-1$ or less edges from $G$, there still remains some chain of edges $p$ such that $\partial(p)=v_1+v_2+v_{3}+v_{4}$, but there is a way to construct the new graph $GN$ by drawing the new edges $\epsilon_1,\ldots,\epsilon_n$ in $G$ as prescribed in the theorem, so that there are no $k$ chains of edges $p_1,\ldots,p_k$ in $GN$, such that no two of $p_1,\ldots,p_k$ have any common nontrivial summands, and that $\partial(p_1)=\cdots=\partial(p_k)=v_1+v_2+v_{3}+v_{4}$ in $GN$. Assume that such a graph $GN$ was constructed and given to us.
\\
\includegraphics[scale=0.5]{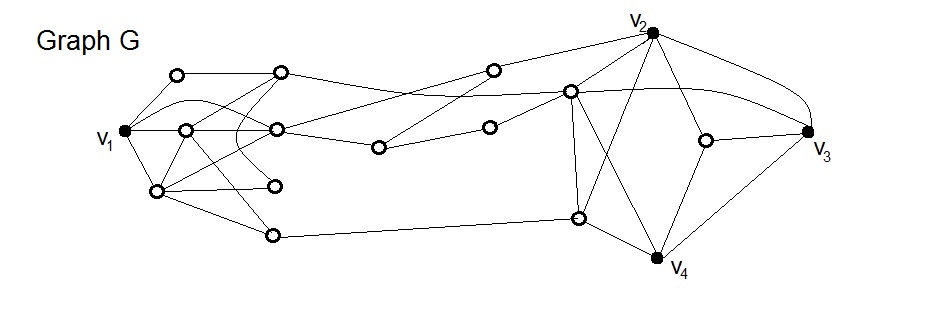}
\\
We claim that for every edge $x_1$ in the graph $G$ there must exist some $k-1$ edges $x_2,\ldots,x_{k}$ in $G$ such that after the deletion of all the $k$ edges $x_1,x_2,\ldots,x_{k}$ from $G$ there does not remain any chain of edges $p$ in $G$ such that $\partial(p)=v_1+v_2+v_{3}+v_{4}$. If our claim is false, then there exists some edge $x_1$ in the graph $G$ such that for any $k-1$ edges $x_2,\ldots,x_{k}$ in $G$, after the deletion of all the $k$ edges $x_1,x_2,\ldots,x_{k}$ from $G$ there still remains some chain of edges $p$ in $G$ such that $\partial(p)=v_1+v_2+v_{3}+v_{4}$. By the minimality of the number of edges is $G$, the graph $G'$, which is the graph $G$ without the edge $x_1$, cannot be a counter-example to our theorem. Let $\partial(x_1)=u_{1,1}+u_{1,2}$ in $G$. Let $y'_1,\ldots,y'_{\mu}$ be all the vertices in $G'$ which are different from the four vertices $v_1,v_2,v_3,v_4$ and which have odd degrees in $G'$.
\\ \\
Since after deleting any $k-1$ edges from $G'$ there still exists some chain of edges $p$ in $G'$ such that $\partial(p)=v_1+v_2+v_{3}+v_{4}$ in $G'$, if one creates a new graph $GN'$ from $G'$ by drawing new edges $\epsilon'_1,\ldots,\epsilon'_{\nu}$ in $G'$ in any such a way that the boundary $\partial(\epsilon'_i)$ of each new edge $\epsilon'_i$ is a sum a vertex from the set $\{y'_1,\ldots,y'_{\mu}\}$ and a vertex from the set $\{v_1,v_2,v_3,v_4,y'_1,\ldots,y'_{\mu}\}$, and that each vertex $y'_1,\ldots,y_{\mu}$ is a summand in the boundary of exactly one of the new edges $\epsilon'_1,\ldots,\epsilon'_{\nu}$, then in the graph $GN'$ one can find $k$ chains of edges $p'_1,\ldots,p'_k$, such that no two of $p'_1,\ldots,p'_k$ have any common nontrivial summands and that $\partial(p'_1)=\cdots=\partial(p'_k)=v_1+v_2+v_{3}+v_{4}$ in $GN'$. We will now show, that this implies that in the $GN$ we can also find $k$ chains of edges $p_1,\ldots,p_k$, such that no two of $p_1,\ldots,p_k$ have any common nontrivial summands and that $\partial(p_1)=\cdots=\partial(p_k)=v_1+v_2+v_{3}+v_{4}$ in $GN$. Thus, the graph $G$ cannot be a counter-example, which means that our claim that there must exist some $k-1$ edges $x_2,\ldots,x_{k}$ in $G$ such that after the deletion of all the $k$ edges $x_1,x_2,\ldots,x_{k}$ from $G$ there does not remain any chain of edges $p$ in $G$ such that $\partial(p)=v_1+v_2+v_{3}+v_{4}$, cannot be not false. We consider three cases:
\\ \\
\textit{Case 1} is when both vertices $u_{1,1}$ and $u_{1,2}$ do not belong to the set $\{y'_1,\ldots,y'_{\mu}\}$. In this case we can assume that $y'_1=y_1,\ldots,y'_{\mu}=y_{m-2}$ and $u_{1,1}=y_{m-1},u_{1,2}=y_m$. We get four possible subcases:
\\
\includegraphics[scale=0.5]{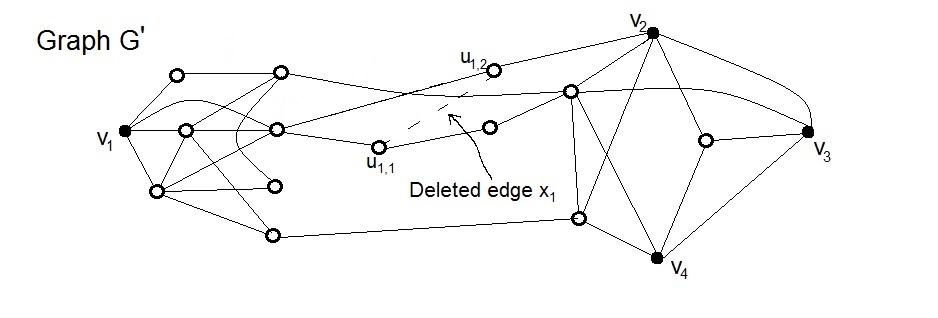}
\\
\textit{Subcase 1 of Case 1} is when for some edge $\epsilon_i$ in $GN$ we have $\partial(\epsilon_i)=u_{1,1}+u_{1,2}$. We can assume that $i=n$. In that subcase we construct the graph $GN'$ from the graph $G'$ by drawing the new edges $\epsilon'_1=\epsilon_1,\ldots,\epsilon'_{\nu}=\epsilon_{n-1}$. It is easy to verify, that this construction satisfies the requirements of our theorem. Thus, we can find $k$ chains of edges $p'_1,\ldots,p'_k$ in the graph $GN'$, such that no two of $p'_1,\ldots,p'_k$ have any common nontrivial summands and that $\partial(p'_1)=\cdots=\partial(p'_k)=v_1+v_2+v_{3}+v_{4}$ in $GN'$. These chains of edges $p'_1,\ldots,p'_k$ can be regarded as chains of edges $p_1,\ldots,p_k$ in the graph $GN$, and $\partial(p_1)=\cdots=\partial(p_k)=v_1+v_2+v_{3}+v_{4}$ in $GN$. Thus, in this subcase the graph $GN$ does not provide a counter-example to our theorem;
\\ \\
\textit{Subcase 2 of Case 1} is when for some two edges $\epsilon_i$ and $\epsilon_j$ in $GN$, where $i\ne j$, we have $\partial(\epsilon_i)=u_{1,1}+w_1$ and $\partial(\epsilon_j)=u_{1,2}+w_2$, where $w_1$ and $w_2$ are vertices belonging to the set $\{v_1,v_2,v_3,v_4\}$. In that subcase $w_1$ and $w_2$ can be the same vertex or two different vertices. We can assume that $i=n-1$ and $j=n$. In that subcase we construct the graph $GN'$ from the graph $G'$ by drawing the new edges $\epsilon'_1=\epsilon_1,\ldots,\epsilon'_{\nu}=\epsilon_{n-2}$. It is easy to verify, that this construction satisfies the requirements of our theorem. Thus, we can find $k$ chains of edges $p'_1,\ldots,p'_k$ in the graph $GN'$, such that no two of $p'_1,\ldots,p'_k$ have any common nontrivial summands and that $\partial(p'_1)=\cdots=\partial(p'_k)=v_1+v_2+v_{3}+v_{4}$ in $GN'$. These chains of edges $p'_1,\ldots,p'_k$ can be regarded as chains of edges $p_1,\ldots,p_k$ in the graph $GN$, and $\partial(p_1)=\cdots=\partial(p_k)=v_1+v_2+v_{3}+v_{4}$ in $GN$. Thus, in this subcase the graph $GN$ does not provide a counter-example to our theorem;
\\ \\
\textit{Subcase 3 of Case 1} is when for some two different edges $\epsilon_i$ and $\epsilon_j$ in $GN$ we have $\partial(\epsilon_i)=u_{1,1}+w_1$ and $\partial(\epsilon_j)=u_{1,2}+w_2$, where one of the two vertices $w_1,w_2$ belongs to the set $\{v_1,v_2,v_3,v_4\}$ and the other one of these two vertices belongs to the set $\{y_1,\ldots,y_{m-2}\}$. We can assume that $w_1\in\{v_1,v_2,v_3,v_4\}$ and $w_2\in\{y_1,\ldots,y_{m-2}\}$, and that $i=n-1$ and $j=n$. In that subcase we construct the graph $GN'$ from the graph $G'$ by drawing the new edges $\epsilon'_1=\epsilon_1,\ldots,\epsilon'_{\nu-1}=\epsilon_{n-2}$ and $\epsilon'_{\nu}$, such that $\partial(\epsilon'_{\nu})=w_1+w_2$. It is easy to verify, that this construction satisfies the requirements of our theorem. Thus, we can find $k$ chains of edges $p'_1,\ldots,p'_k$ in the graph $GN'$, such that no two of $p'_1,\ldots,p'_k$ have any common nontrivial summands and that $\partial(p'_1)=\cdots=\partial(p'_k)=v_1+v_2+v_{3}+v_{4}$ in $GN'$. Substituting $\epsilon_{n-1}+x_1+\epsilon_n$ for $\epsilon'_{\nu}$ whenever we find the edge $\epsilon'_{\nu}$ as a summand in one of the chains of edges $p'_1,\ldots,p'_k$, produces $k$ chains of edges $p_1,\ldots,p_k$ in the graph $GN$ such that no two of these $k$ chains of edges have any common summands and that $\partial(p_1)=\cdots=\partial(p_k)=v_1+v_2+v_{3}+v_{4}$ in $GN$;
\\
\includegraphics[scale=0.5]{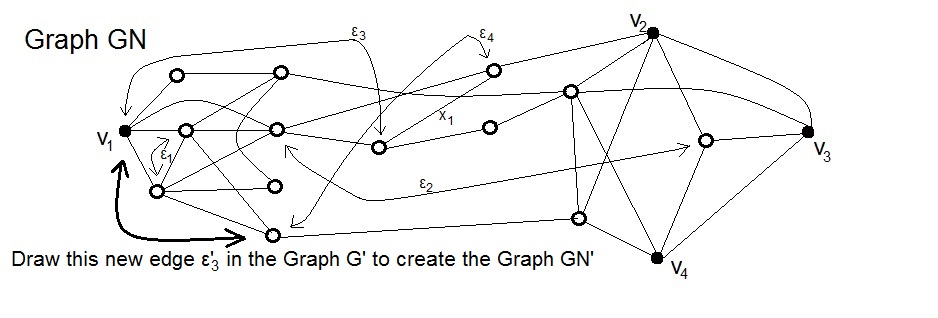}
\\
\textit{Subcase 4 of Case 1} is when for some two different edges $\epsilon_i$ and $\epsilon_j$ in $GN$ we have $\partial(\epsilon_i)=u_{1,1}+w_1$ and $\partial(\epsilon_j)=u_{1,2}+w_2$, where $w_1$ and $w_2$ are two different vertices belonging to the set $\{y_1,\ldots,y_{m-2}\}$. We can assume that $i=n-1$ and $j=n$. In that subcase, just like in the previous subcase, we construct the graph $GN'$ from the graph $G'$ by drawing the new edges $\epsilon'_1=\epsilon_1,\ldots,\epsilon'_{\nu-1}=\epsilon_{n-2}$ and $\epsilon'_{\nu}$, such that $\partial(\epsilon'_{\nu})=w_1+w_2$. It is easy to verify, that this construction satisfies the requirements of our theorem. Thus, we can find $k$ chains of edges $p'_1,\ldots,p'_k$ in the graph $GN'$, such that no two of $p'_1,\ldots,p'_k$ have any common nontrivial summands and that $\partial(p'_1)=\cdots=\partial(p'_k)=v_1+v_2+v_{3}+v_{4}$ in $GN'$. Substituting $\epsilon_{n-1}+x_1+\epsilon_n$ for $\epsilon'_{\nu}$ whenever we find the edge $\epsilon'_{\nu}$ as a summand in one of the chains of edges $p'_1,\ldots,p'_k$, produces $k$ chains of edges $p_1,\ldots,p_k$ in the graph $GN$ such that no two of these $k$ chains of edges have any common summands and that $\partial(p_1)=\cdots=\partial(p_k)=v_1+v_2+v_{3}+v_{4}$ in $GN$. This concludes our Case 1;
\\
\includegraphics[scale=0.5]{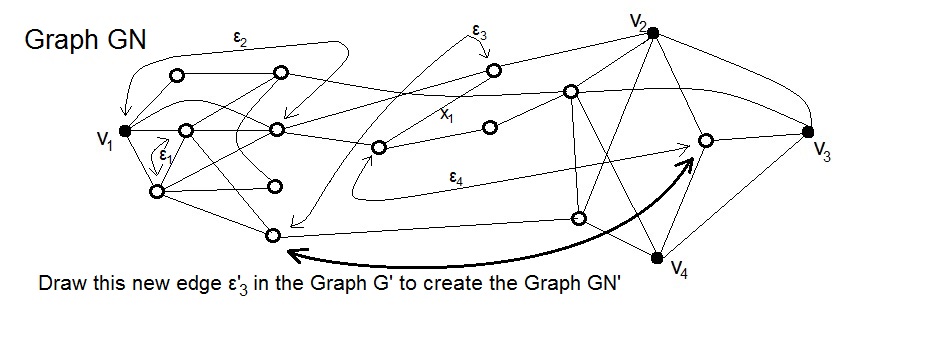}
\\
\textit{Case 2} is when one of the two vertices $u_{1,1}$ and $u_{1,2}$ belongs to the set $\{y'_1,\ldots,y'_{\mu}\}$ and the other one of these two vertices does not belong to the set $\{y'_1,\ldots,y'_{\mu}\}$. In this case we can assume that $y'_1=y_1,\ldots,y'_{\mu-1}=y_{m-1},y'_{\mu}=u_{1,2}$ and $u_{1,1}=y_m$.
\\
\includegraphics[scale=0.5]{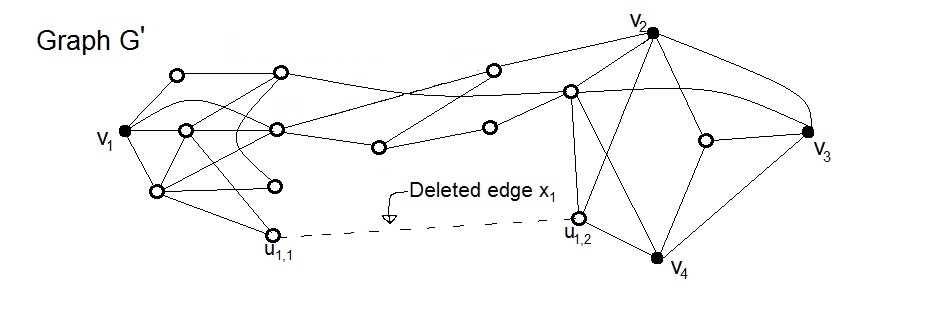}
\\
In this case, for some edge $\epsilon_i$ in $GN$ we have $\partial(\epsilon_i)=w+u_{1,1}$, where the vertex $w$ is either one of the four vertices $v_1,v_2,v_3,v_4$ or one of the $m-1$ vertices $y_1,\ldots,y_{m-1}$. We can assume that $i=n$. In that case we construct the graph $GN'$ from the graph $G'$ by drawing the new edges $\epsilon'_1=\epsilon_1,\ldots,\epsilon'_{\nu-1}=\epsilon_{n-1}$ and $\epsilon'_{\nu}$, such that $\partial(\epsilon'_{\nu})=w+u_{1,2}$. It is easy to verify, that this construction satisfies the requirements of our theorem. Thus, we can find $k$ chains of edges $p'_1,\ldots,p'_k$ in the graph $GN'$, such that no two of $p'_1,\ldots,p'_k$ have any common nontrivial summands and that $\partial(p'_1)=\cdots=\partial(p'_k)=v_1+v_2+v_{3}+v_{4}$ in $GN'$. Substituting $x_1+\epsilon_n$ for $\epsilon'_{\nu}$ whenever we find the edge $\epsilon'_{\nu}$ as a summand in one of the chains of edges $p'_1,\ldots,p'_k$, produces $k$ chains of edges $p_1,\ldots,p_k$ in the graph $GN$ such that no two of these $k$ chains of edges have any common summands and that $\partial(p_1)=\cdots=\partial(p_k)=v_1+v_2+v_{3}+v_{4}$ in $GN$. This concludes our Case 2;
\\ \\
\textit{Case 3} is when both vertices $u_{1,1}$ and $u_{1,2}$ belong to the set $\{y'_1,\ldots,y'_{\mu}\}$. In this case we can assume that $y'_1=y_1,\ldots,y'_{\mu-2}=y_{m},y'_{\mu-1}=u_{1,1},y'_{\mu}=u_{1,2}$. In that case we construct the graph $GN'$ from the graph $G'$ by drawing the new edges $\epsilon'_1=\epsilon_1,\ldots,\epsilon'_{\nu-1}=\epsilon_{n}$ and $\epsilon'_{\nu}$, such that $\partial(\epsilon'_{\nu})=u_{1,1}+u_{1,2}$. It is easy to verify, that this construction satisfies the requirements of our theorem. The edge $\epsilon'_{\nu}$ in $GN'$ can be regarded as the edge $x_1$ in $GN$, and the graphs $GN'$ and $GN$ are thus isomorphic. Hence, we can find $k$ chains of edges $p_1,\ldots,p_k$ in the graph $GN$, such that no two of $p_1,\ldots,p_k$ have any common nontrivial summands and that $\partial(p_1)=\cdots=\partial(p_k)=v_1+v_2+v_{3}+v_{4}$ in $GN$. This concludes our Case 3.
\\ \\
We showed that if for some edge $x_1$ in the graph $G$, for any $k-1$ edges $x_2,\ldots,x_{k}$ in $G$, the deletion of all the $k$ edges $x_1,x_2,\ldots,x_{k}$ from $G$ does not break all the chains of edges $p$ in $G$ such that $\partial(p)=v_1+v_2+v_{3}+v_{4}$, then regardless of how one constructs the graph $GN$ from $G$, as long as this construction satisfies the requirements of our theorem, we can find $k$ chains of edges $p_1,\ldots,p_k$ in the graph $GN$, such that no two of $p_1,\ldots,p_k$ have any common nontrivial summands and that $\partial(p_1)=\cdots=\partial(p_k)=v_1+v_2+v_{3}+v_{4}$ in $GN$. Thus, $G$ is not a counter-example to our theorem. This implies that for every edge $x_1$ in the graph $G$ there must exist some $k-1$ edges $x_2,\ldots,x_{k}$ in $G$ such that after the deletion of all the $k$ edges $x_1,x_2,\ldots,x_{k}$ from $G$ there does not remain any chain of edges $p$ in $G$ such that $\partial(p)=v_1+v_2+v_{3}+v_{4}$.
\\
\includegraphics[scale=0.5]{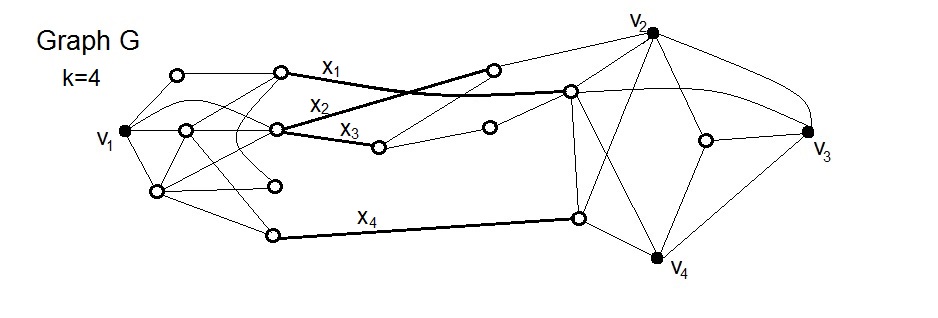}
\\
Next, assume that there is an edge $x_1$ in $G$ such that its boundary $\partial(x_1)=u_{1,1}+u_{1,2}$ does not contain any one of the four vertices $v_1,v_2,v_3,v_4$ among its two summands $u_{1,1}$ and $u_{1,2}$. We find some $k-1$ edges $x_2,\ldots,x_{k}$ in $G$ such that after the deletion of all the $k$ edges $x_1,x_2,\ldots,x_{k}$ from the graph $G$, there no longer exists any chain of edges $p$ in $G$ such that $\partial(p)=v_1+v_2+v_{3}+v_{4}$. We denote the boundary $\partial(x_i)$ of each edge $x_i$ by $u_{i,1}+u_{i,2}$, with $u_{i,1}$ and $u_{i,2}$ being two different vertices in $G$. We construct the graph $\Gamma'$ by deleting the edges $x_2,\ldots,x_{k}$ from the graph $G$. There exists some chain of edges $p$ in $\Gamma'$ such that $\partial(p)=v_1+v_2+v_{3}+v_{4}$ in $\Gamma'$, but if we delete the edge $x_1$ from the graph $\Gamma'$ then there will not remain any chain of edges in $\Gamma'$ such that its boundary is $v_1+v_2+v_{3}+v_{4}$. Let $K$ be the connected component of the graph $\Gamma'$ which contains the vertices $u_{1,1}$ and $u_{1,2}$. Since the deletion of the edge $x_1$ breaks in $\Gamma'$ some chain of edges whose boundary is $v_1+v_2+v_{3}+v_{4}$, $K$ must contain either two or four of the vertices from the set $\{v_1,v_2,v_{3},v_{4}\}$.
\\ \\
From Lemma \ref{tech} it follows that the deletion of the edge $x_1$ from the graph $\Gamma'$ must break the connected component $K$ of the graph $\Gamma'$ into two connected components $K_1$ and $K_2$ with the vertex $u_{1,1}$ contained in the connected component $K_1$ and the vertex $u_{1,2}$ contained in the connected component $K_2$. Indeed, the deletion of $x_1$ has no effect on the other connected components of the graph $\Gamma'$, and it follows from Lemma \ref{tech} that if the deletion of the edge $x_1$ does not break the connected component $K$ then the deletion of $x_1$ has no effect on the existence or the nonexistence of a chain of edges $p'$ in $K$ such that $\partial(p')$ is the sum of all the vertices from the set $\{v_1,v_2,v_{3},v_{4}\}$ which belong to $K$.
\\ \\
Furthermore, it follows from Lemma \ref{tech} that $K_1$ and $K_2$ must contain odd numbers of the vertices from the set $\{v_1,v_2,v_{3},v_{4}\}$. Thus, at least one of these two connected components must contain only one vertex from the set $\{v_1,v_2,v_{3},v_{4}\}$. We can assume that $K_1$ contains the vertex $v_1$. Now, hypothetically, unless we prove otherwise, the graph $\Gamma'$ may contain more connected components than just $K$. Moreover, hypothetically, one of these other connected components may contain two of the vertices $v_2,v_{3},v_{4}$. Since this matter is not essential for this proof, we will not provide an argument why this cannot happen in our minimal example, and why $K_2$ will always contain all the three vertices $v_2,v_{3},v_{4}$.
\\ \\
We construct a new graph $T$ from the graph $G$ by performing the following three steps:
\\
\textit{Step 1:} Draw a new vertex $v'_1$;
\\
\textit{Step 2:} For every edge $e$ in $G$ such that $\partial(e)=u+w$ with the vertex $u$ belonging to $K_1$ and the vertex $w$ not belonging to $K_1$, define the boundary $\partial(e)$ of the edge $e$ in $T$ to be $\partial(e)=v'_1+w$;
\\
\textit{Step 3:} Delete all the vertices and all the edges belonging to $K_1$.
\\ \\
The boundary operation for all the other edges in $T$ remains the same as it was in $G$. Since $u_{1,1}\ne v_1$, the connected component $K_1$ contains at least one edge. Hence, the graph $T$ has less edges than the graph $G$, which, by the minimality of edges in $G$, implies that $T$ cannot be a counter-example to our theorem.
\\
\includegraphics[scale=0.5]{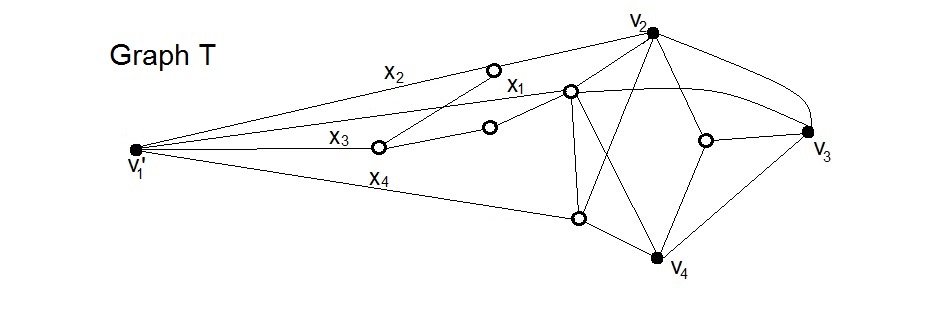}
\\
Next, we claim that if there exist some $k-1$ or less edges in $T$, such that after the deletion of these edges from $T$ there is no chain of edges $p'$ in $T$ such that $\partial(p')=v'_1+v_2+v_3+v_4$ in $T$, then after the deletion of these same $k-1$ or less edges from $G$ there is no chain of edges $p$ in $G$ such that $\partial(p)=v_1+v_2+v_3+v_4$ in $G$.
\\ \\
Assume that our claim is wrong and that after the deletion of these same $k-1$ or less edges from $G$ there still exists a chain of edges $p$ in $G$ such that $\partial(p)=v_1+v_2+v_3+v_4$ in $G$. We remove from this chain of edges $p$ every summand $e$ which belongs to $K_1$, and we obtain a new chain of edges $p'$, which can be regarded as a chain of edges in $T$. The vertex $v_1$ must appear as a summand in the boundary of an odd number of edges which are summands in $p$, and every other vertex belonging to $K_1$ must appear as a summand in the boundary of an even number of edges which are summands in $p$. When we remove an edge $e$ which belongs to $K_1$ from $p$, two vertices belonging to $K_1$ get removed from appearing as summands in $\partial(p)$ in $G$. Hence, the total number of vertices from $K_1$ which appear as summands in the boundary $\partial(p')$ of $p'$ in $G$ is odd. Thus, $\partial(p')$ in $T$ will contain a sum of an odd number of $v'_1$, which implies that $\partial(p')$ in $T$ contains $v'_1$ as a summand. With respect to all the vertices not belonging to $K_1$, they appear in $\partial(p')$ in $T$ in the same way as they appear in $\partial(p')$ in $G$, and they appear in $\partial(p')$ in $G$ in the same way as they appear in $\partial(p)$ in $G$. Thus, $\partial(p')$ in $T$ is $v'_1+v_2+v_{3}+v_{4}$.
\\ \\
This implies that our claim is right, which, in its turn, implies that after a deletion of any $k-1$ or less edges from $T$ there remains some chain of edges $p'$ in $T$ such that $\partial(p')=v'_1+v_2+v_3+v_4$ in $T$.
\\ \\
Now, for every vertex $w$ in $G$ which does not belong to $K_1$, the degree of $w$ in $T$ is the same as the degree of $w$ in $G$. Let $y_1,\ldots,y_r$ be all the vertices in the set $\{y_1,\ldots,y_m\}$ which do not belong to $K_1$. The vertices $y_1,\ldots,y_r$ are the only vertices in $T$, which do not belong to the set $\{v'_1,v_2,v_3,v_4\}$, and which have odd degrees in $T$. Since $T$ cannot be a counter-example to our theorem, if one creates a new graph $TN$ from the graph $T$ by drawing new edges $\epsilon_{T,1},\ldots,\epsilon_{T,t}$ in $T$ in any such a way that the boundary $\partial(\epsilon_{T,i})$ of each new edge $\epsilon_{T,i}$ is a sum of two vertices from the set $\{v'_1,v_2,v_3,v_4,y_1,\ldots,y_r\}$ with at least one of these two vertices belonging to the set $\{y_1,\ldots,y_r\}$, and that each vertex $y_1,\ldots,y_r$ is a summand in the boundary of exactly one of the new edges $\epsilon_{T,1},\ldots,\epsilon_{T,t}$, then in the graph $TN$ there exist $k$ chains of edges $p_{T,1},\ldots,p_{T,k}$, such that no two of $p_{T,1},\ldots,p_{T,k}$ have any common nontrivial summands and that $\partial(p_{T,1})=\cdots=\partial(p_{T,k})=v'_1+v_2+v_{3}+v_{4}$ in $TN$.
\\ \\
Next, we construct a new graph $L$ from the graph $G$ by performing the following three steps:
\\
\textit{Step 1:} Draw a new vertex $v'_2$;
\\
\textit{Step 2:} For every edge $e$ in $G$ such that $\partial(e)=u+w$ with the vertex $u$ belonging to $K_1$ and the vertex $w$ not belonging to $K_1$, define the boundary $\partial(e)$ of the edge $e$ in $L$ to be $\partial(e)=v'_2+w$;
\\
\textit{Step 3:} Delete all the vertices not belonging to $K_1$, and delete all the edges, whose boundaries are sums of two vertices not belonging to $K_1$.
\\ \\
The boundary operation for all the edges in $G$ which belong to $K_1$ remains in $L$ the same as it was in $G$. If there exist some $k-1$ or less edges in $L$, such that after the deletion of these edges from $L$ there is no chain of edges $p'$ in $L$ such that $\partial(p')=v_1+v'_2$ in $L$, then the deletion of these same $k-1$ or less edges from $G$ disconnects in $G$ the vertex $v_1$ from all the vertices not belonging to $K_1$. Thus, after this deletion, there would not remain a chain of edges $p$ in $G$ such that $\partial(p)=v_1+v_2+v_3+v_4$ in $G$, which would imply that $G$ is not a counter-example. Hence, after a deletion of any $k-1$ or less edges from $L$ there still remains some chain of edges in $L$ such that its boundary in $L$ is $v_1+v'_2$. Thus, by Theorem \ref{Menger.homol.2}, we can find $k$ chains of edges $p_{L,1},\ldots,p_{L,k}$ in $L$, such that no two of $p_{L,1},\ldots,p_{L,k}$ have any common nontrivial summands and that $\partial(p_{L,1})=\cdots=\partial(p_{L,k})=v_1+v'_2$ in $L$.
\\ \\
If for some edge $e$ in $G$ we have $\partial(e)=u+w$ with vertex $u$ belonging to $K_1$ and vertex $w$ not belonging to $K_1$ then the edge $e$ must be one of the edges $x_1,x_2,\ldots,x_{k}$. Indeed, otherwise $K_1$ cannot be a connected component of the graph $\Gamma$, which is the $G$ without the edges $x_1,x_2,\ldots,x_{k}$. Thus, each one of the $k$ edges $x_1,x_2,\ldots,x_{k}$, when we regard then as edges in $L$, appears as a summand in one and unique chain of edges from the $k$ chains of edges $p_{L,1},p_{L,2},\ldots,p_{L,k}$ in $L$. We assume that each $x_i$ appears in the chain of edges $p_{L,i}$. Since in the graph $G$ we have $\partial(x_i)=u_{i,1}+u_{i,2}$, one of the two vertices $u_{i,1}$ and $u_{i,2}$, for every $i$, must belong to $K_1$ and the other one cannot belong to $K_1$. We can assume that for all $i=1,\ldots,k$, the vertex $u_{i,1}$ belongs to $K_1$ and the vertex $u_{i,2}$ does not belong to $K_1$. Now, in $L$ we get $\partial(x_i)=u_{i,1}+v'_2$ for all $i=1,\ldots,k$.
\\
\includegraphics[scale=0.5]{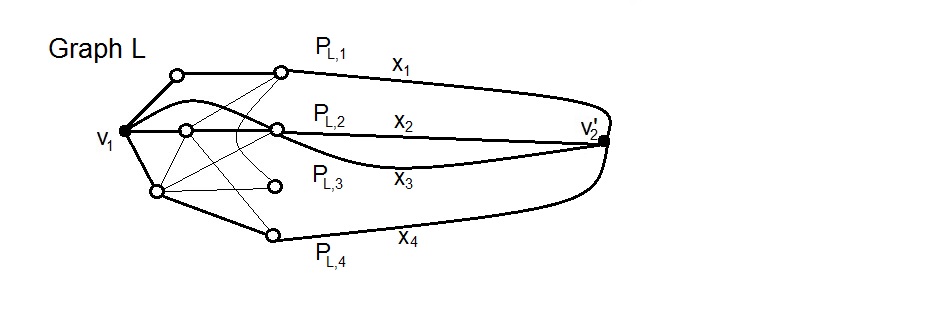}
\\
Let $\Theta$ be the set of all the edges which appear as summands in any of the chains of edges $p_{L,1},\ldots,p_{L,k}$. Here we abuse the notation and regard the edges, appearing as summands in $p_{L,1},\ldots,p_{L,k}$, as edges in $L$ and as edges in $G$. Thus, we regard the set $\Theta$ as a subset of the set of edges of $L$ and as a subset of the set of edges of $G$. We do remind, that the boundary operation in the graph $G$ is the same for all the edges belonging to $L$ as it is the the graph $L$, except for the edges $x_1,x_2,\ldots,x_{k}$ where we have $\partial(x_i)=u_{i,1}+u_{i,2}$ in $G$ and $\partial(x_i)=u_{i,1}+v'_2$ in $L$. Thus, when we speak of edges from the set $\Theta$ as edges in $G$, we mean that the boundary operation is the boundary operation in $G$, and when we speak of edges from the set $\Theta$ as edges in $L$, we mean that the boundary operation is the boundary operation in $L$. Each time we mention the set $\Theta$ or its elements, we will specify if we regard them in $G$ or in $L$.
\\
\includegraphics[scale=0.5]{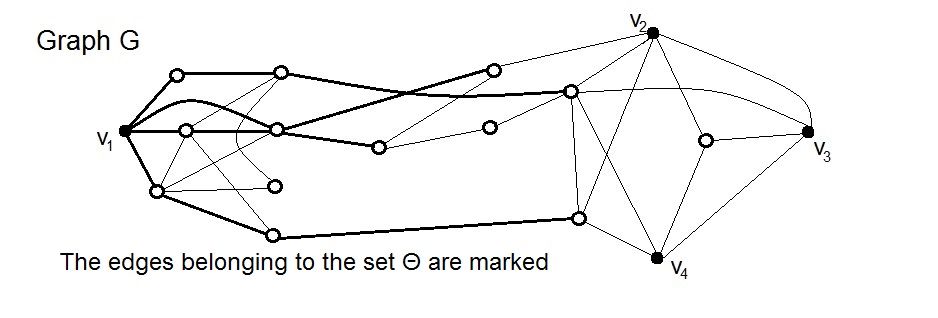}
\\
We are now ready to prove that the graph $G$ cannot be a counter-example to our theorem. Assume that someone constructed a new graph $GN$ by drawing new edges $\epsilon_1,\ldots,\epsilon_n$ in $G$ so that the boundary $\partial(\epsilon_i)$ of each new edge $\epsilon_i$ is a sum of two vertices from the set $\{v_1,v_2,v_3,v_4,y_1,\ldots,y_m\}$ with at least one of these two vertices being one of the vertices $y_1,\ldots,y_m$, and that each vertex $y_1,\ldots,y_m$ is a summand in the boundary of exactly one of the new edges $\epsilon_1,\ldots,\epsilon_n$. We want to show that regardless of how the new graph $GN$ was constructed, there exist $k$ chains of edges $p_1,\ldots,p_k$, such that no two of $p_1,\ldots,p_k$ have any common nontrivial summands and that $\partial(p_1)=\cdots=\partial(p_k)=v_1+v_2+v_{3}+v_{4}$ in $GN$.
\\ \\
We remind, that the vertices $y_1,\ldots,y_r$ in $G$ are all the vertices in $G$, which do not belong to the set $\{v_1,v_2,v_3,v_4\}$ and do not belong to $K_1$, and which have odd degrees in $G$, and that the vertices $y_{r+1},\ldots,y_m$ in $G$ are all the vertices in $G$, which do not belong to the set $\{v_1,v_2,v_3,v_4\}$ and which belong to $K_1$, and which have odd degrees in $G$. Let $\epsilon_1,\ldots,\epsilon_{\tau}$ be all the edges in the set $\{\epsilon_1,\ldots,\epsilon_n\}$ such that their boundaries are either a sum of two elements from the set $\{y_1,\ldots,y_r\}$ or a sum of an element from the set $\{y_1,\ldots,y_r\}$ and an element from the set $\{v_1,v_2,v_3,v_4\}$. Let $\epsilon_{\tau+1},\ldots,\epsilon_{\tau+\lambda}$ be all the edges in the set $\{\epsilon_1,\ldots,\epsilon_n\}$ such that their boundaries are either a sum of two elements from the set $\{y_{r+1},\ldots,y_m\}$ or a sum of an element from the set $\{y_{r+1},\ldots,y_m\}$ and an element from the set $\{v_1,v_2,v_3,v_4\}$. Thus, the boundaries of all the remaining new edges $\epsilon_{\tau+\lambda+1},\ldots,\epsilon_{n}$ must be sums of a vertex from the set $\{y_1,\ldots,y_r\}$ and a vertex from the set $\{y_{r+1},\ldots,y_m\}$.
\\ \\
We are now going to construct a new graph $TN$ from the graph $T$ by drawing new edges $\epsilon_{T,1},\ldots,\epsilon_{T,t}$ in $T$ so that the boundary $\partial(\epsilon_{T,i})$ of each new edge $\epsilon_{T,i}$ is a sum of two vertices from the set $\{v'_1,v_2,v_3,v_4,y_1,\ldots,y_r\}$ with at least one of these two vertices being one of the vertices $y_1,\ldots,y_r$, and that each vertex $y_1,\ldots,y_r$ is a summand in the boundary of exactly one of the new edges $\epsilon_{T,1},\ldots,\epsilon_{T,t}$, in such a way, that if in the graph $TN$ there exist $k$ chains of edges $p_{T,1},\ldots,p_{T,k}$, such that no two of $p_{T,1},\ldots,p_{T,k}$ have any common nontrivial summands and that $\partial(p_{T,1})=\cdots=\partial(p_{T,k})=v'_1+v_2+v_{3}+v_{4}$ in $TN$, then in the graph $GN$ there exist $k$ chains of edges $p_{1},\ldots,p_{k}$, such that no two of $p_{1},\ldots,p_{k}$ have any common nontrivial summands and that $\partial(p_{1})=\cdots=\partial(p_{k})=v_1+v_2+v_{3}+v_{4}$ in $GN$. This will imply that $G$ is not a counter-example.
\\ \\
We define $\Omega=\Theta$, and we regard $\Omega$ as a set of edges in the graph $GN$. At each step of the process, which we describe below, we will modify the set $\Omega$ by adding into it another edge belonging to the graph $GN$. Notice, that at this point, for any vertex $y\ne v_1$ in $K_1$, the number of edges in $\Omega$ which contain $y$ as a summand in their boundaries must be even. After a possible re-indexing, we will have all the vertices $\{y_1,\ldots,y_{\varrho}\}$ appear as summands in the boundaries of the edges $\epsilon_1,\ldots,\epsilon_{\tau}$, and all the vertices $\{y_{\varrho+1},\ldots,y_{r}\}$ appear as summands in the boundaries of the edges $\epsilon_{\tau+\lambda+1},\ldots,\epsilon_{n}$. Moreover, we can assume that $\partial(\epsilon_{\tau+\lambda+1})=y_{\varrho+1}+y_{r+1}$.
\\
\includegraphics[scale=0.5]{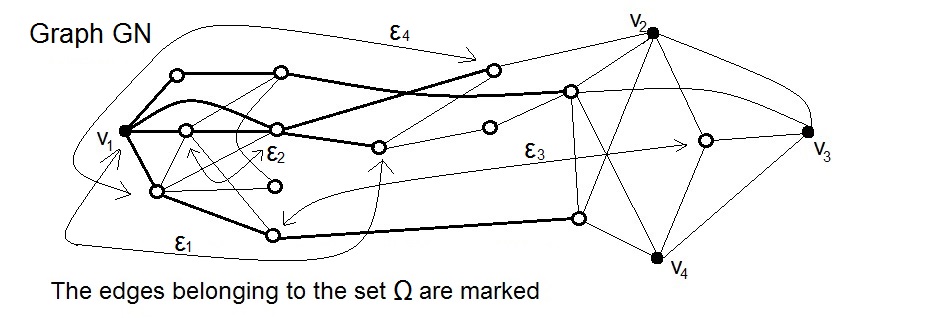}
\\ \\
\textit{The Step 1 of the Round 1} of our process is defining the chain of edges $\gamma_{\tau+1}$ in $GN$ as $\gamma_{\tau+1}=\epsilon_{\tau+\lambda+1}$, and adding the edge $\epsilon_{\tau+\lambda+1}$ as an element into the set $\Omega$. Now the number of edges in $\Omega$ which contain the vertex $y_{r+1}$ as a summand in their boundaries must be odd, while the degree of the vertex $y_{r+1}$ in the graph $GN$ is even. Thus, there exists some edge $\alpha_{(\tau+\lambda+1,1)}$ in $GN$, which does not belong to $\Omega$, and which contains the vertex $y_{r+1}$ as a summand in its boundary. Since the edge $\epsilon_{\tau+\lambda+1}$ and all the edges $x_1,x_2,\ldots,x_{k}$ belong to $\Omega$, the edge $\alpha_{(\tau+\lambda+1,1)}$ must belong to $K_1$. Thus, $\partial(\alpha_{(\tau+\lambda+1,1)})=y_{r+1}+w_1$ with the vertex $w_1$ belonging to $K_1$.
\\
\includegraphics[scale=0.5]{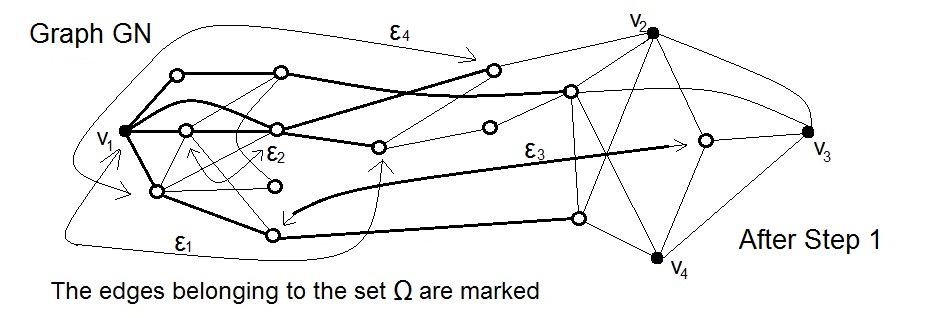}
\\ \\
\textit{The Step 2 of the Round 1} of our process is redefining the chain of edges $\gamma_{\tau+1}$ to be our previous $\gamma_{\tau+1}$ plus the edge $\alpha_{(\tau+\lambda+1,1)}$, and adding the edge $\alpha_{(\tau+\lambda+1,1)}$ as an element into the set $\Omega$. If the vertex $w_1$ is $v_1$, then we end the Round 1 of our process, and we proceed to the Round 2 of our process, starting its Step 1 with the edge $\epsilon_{\tau+\lambda+2}$. Otherwise, the number of edges in $\Omega$ which contain $w_1$ as a summand in their boundaries must be odd, while the degree of the vertex $w_1$ in the graph $GN$ is even. Thus, there exists some edge $\alpha_{(\tau+\lambda+1,2)}$ in $GN$, which does not belong to $\Omega$, and which contains the vertex $w_1$ as a summand in its boundary. Let $\partial(\alpha_{(\tau+\lambda+1,2)})=w_1+w_2$.
\\
\includegraphics[scale=0.5]{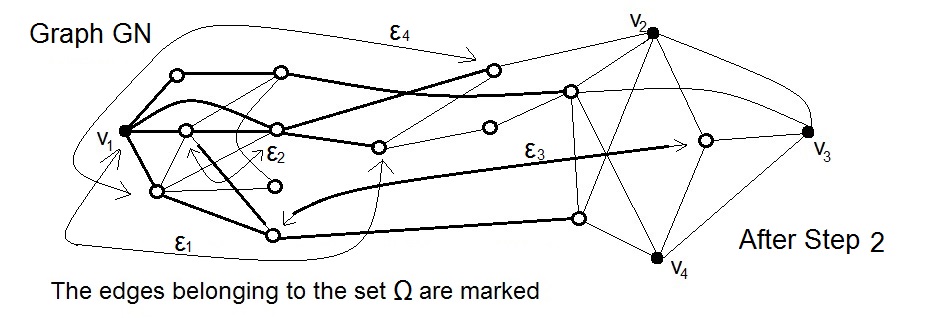}
\\ \\
\textit{The Step 3 of the Round 1} of our process is redefining the chain of edges $\gamma_{\tau+1}$ to be our previous $\gamma_{\tau+1}$ plus the edge $\alpha_{(\tau+\lambda+1,2)}$, and adding the edge $\alpha_{(\tau+\lambda+1,2)}$ as an element into the set $\Omega$.
\\ \\
If the edge $\alpha_{(\tau+\lambda+1,2)}$ is an edge in the graph $G$ then, since all the edges $x_1,x_2,\ldots,x_{k}$ belong to $\Omega$, the vertex $w_2$ must belong to $K_1$. In that case, if we have $w_2=v_1$ then we end the Round 1 of our process, and we proceed to the Round 2 of our process, starting its Step 1 with the edge $\epsilon_{\tau+\lambda+2}$, and if $w_2\ne v_1$ then the number of edges in $\Omega$ which contain $w_2$ as a summand in their boundaries must be odd, while the degree of the vertex $w_2$ in the graph $GN$ is even. Thus, there exists some edge $\alpha_{(\tau+\lambda+1,3)}$ in $GN$ with $\partial(\alpha_{(\tau+\lambda+1,3)})=w_2+w_3$, which does not belong to $\Omega$. At this point we proceed to the Step 4 of the Round 1, which we start by adding the edge $\alpha_{(\tau+\lambda+1,3)}$ as a summand to our $\gamma_{\tau+1}$, and as an element into our $\Omega$.
\\ \\
If the edge $\alpha_{(\tau+\lambda+1,2)}$ is not an edge in the graph $G$ then we have the following two cases:
\\ \\
\textit{Case 1} is if $\alpha_{(\tau+\lambda+1,2)}\in\{\epsilon_{\tau+1},\ldots,\epsilon_{\tau+\lambda}\}$. In this case the vertex $w_2$ is either one of the four vertices $v_1,v_2,v_3,v_4$ or one of the vertices $y_{r+2},\ldots,y_m$. We cannot have $w_2=y_{r+1}$ because $y_{r+1}$ appears as a summand in the boundary of the new edge $\epsilon_{\tau+\lambda+1}$, which prevents $y_{r+1}$ from appearing as a summand in the boundary of any other new edge which was drawn in $G$ to create $GN$. If the vertex $w_2$ is one of the four vertices $v_1,v_2,v_3,v_4$ then we end the Round 1 of our process, and we proceed to the Round 2 of our process, starting its Step 1 with the edge $\epsilon_{\tau+\lambda+2}$. If the vertex $w_2$ is one of the vertices $y_{r+2},\ldots,y_m$ then the number of edges in $\Omega$ which contain $w_2$ as a summand in their boundaries must be odd, while the degree of the vertex $w_2$ in the graph $GN$ is even. Thus, there exists some edge $\alpha_{(\tau+\lambda+1,3)}$ in $GN$ with $\partial(\alpha_{(\tau+\lambda+1,3)})=w_2+w_3$, which does not belong to $\Omega$. At this point we proceed to the Step 4 of the Round 1, which we start by adding the edge $\alpha_{(\tau+\lambda+1,3)}$ as a summand to our $\gamma_{\tau+1}$, and as an element into our $\Omega$.
\\
\includegraphics[scale=0.5]{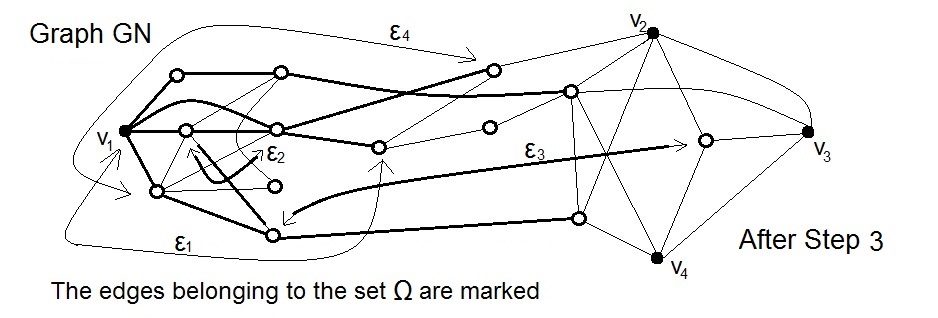}
\\ \\
\textit{Case 2} is if $\alpha_{(\tau+\lambda+1,2)}\in\{\epsilon_{\tau+\lambda+2},\ldots,\epsilon_{n}\}$. In this case the vertex $w_2$ is one of the vertices $\{y_{\varrho+2},\ldots,y_{r}\}$. We cannot have $w_2=y_{\varrho+1}$ because $y_{\varrho+1}$ appears as a summand in the boundary of the new edge $\epsilon_{\tau+\lambda+1}$, which prevents $y_{\varrho+1}$ from appearing as a summand in the boundary of any other new edge which was drawn in $G$ to create $GN$. Thus, we can assume that $w_2=y_{\varrho+2}$ and that $\alpha_{(\tau+\lambda+1,2)}=\epsilon_{\tau+\lambda+2}$. At this point we end the Round 1 of our process, and we proceed to the Round 2 of our process, starting its Step 1 with the edge $\epsilon_{\tau+\lambda+3}$.
\\
\includegraphics[scale=0.5]{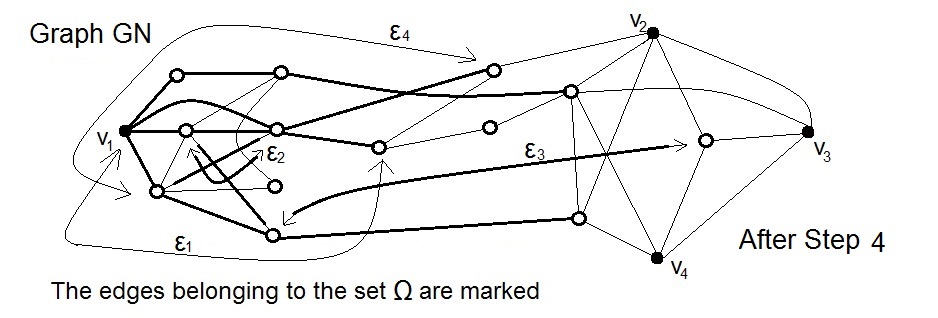}
\\ \\
After a finite number of steps, where Step 4 is similar to Step 3, and so on, we will end our Round 1. If by this point we did not use all the edges $\epsilon_{\tau+\lambda+1},\ldots,\epsilon_n$, we proceed to our Round 2. The Round 2 is identical to the Round 1 in the sense, that we start that round with every vertex $y\ne v_1$ in $K_1$ appearing as a summand in the boundary of an even number of edges belonging to $\Omega$. After a finite number such rounds, all the edges $\epsilon_{\tau+\lambda+1},\ldots,\epsilon_n$ will be used in our process.
\\
\includegraphics[scale=0.5]{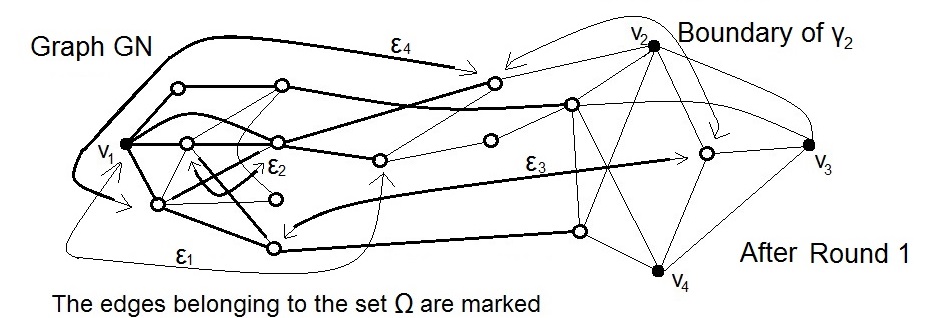}
\\ \\
At this point we have some chains of edges $\gamma_{\tau+1},\gamma_{\tau+2},\ldots,\gamma_{\tau+\varphi}$ in $GN$, such that every vertex $y_{\varrho+1},\ldots,y_{r}$ appears as a summand in the boundary of exactly one of the $\gamma_{\tau+1},\gamma_{\tau+2},\ldots,\gamma_{\tau+\varphi}$, and that the boundary of each one of  $\gamma_{\tau+1},\gamma_{\tau+2},\ldots,\gamma_{\tau+\varphi}$ is either a sum of two of the vertices  $y_{\varrho+1},\ldots,y_{r}$, or of one of the vertices $y_{\varrho+1},\ldots,y_{r}$ and one of the vertices $v_1,v_2,v_3,v_4$. Moreover, no two of the chains of edges $\gamma_{\tau+1},\gamma_{\tau+2},\ldots,\gamma_{\tau+\varphi},p_{L,1},p_{L,2},\ldots,p_{L,k}$ in $GN$ have any common nontrivial summands. Finally, none of the chains of edges $\gamma_{\tau+1},\gamma_{\tau+2},\ldots,\gamma_{\tau+\varphi},p_{L,1},p_{L,2},\ldots,p_{L,k}$ in $GN$ contains any of the edges $\epsilon_1,\ldots,\epsilon_{\tau}$ as a summand in it.
\\ \\
\includegraphics[scale=0.5]{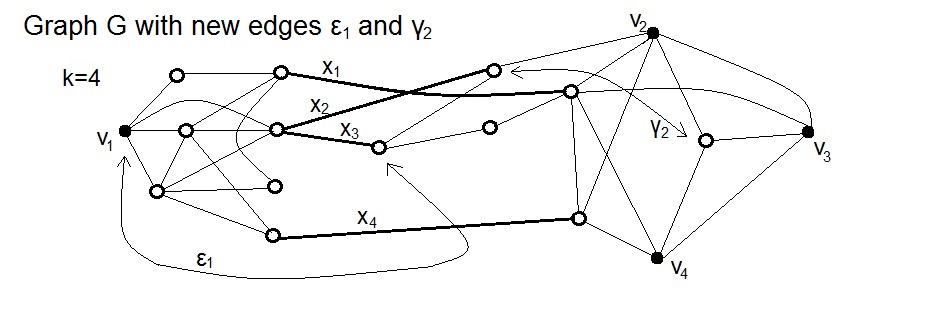}
\\
Now we take the graph $T$, and we construct the graph $TN$ by drawing the new edges $\epsilon_{T,1},\ldots,\epsilon_{T,t}$, where $t=\tau+\varphi$, in $T$ in the following way:
\\ \\
For all $1\le i\le\tau$ we define the boundary $\partial(\epsilon_{T,i})$ of the each edge $\epsilon_{T,i}$ in the graph $TN$ to be equal to the boundary $\partial(\epsilon_i)$ of the edge $\epsilon_i$ in the graph $GN$ with the vertex $v_1$, if the vertex $v_1$ appears as a summand in $\partial(\epsilon_i)$, substituted by the vertex $v'_1$;
\\ \\
For all $\tau+1\le i\le t$ we define the boundary $\partial(\epsilon_{T,i})$ of the each edge $\epsilon_{T,i}$ in the graph $TN$ to be equal to the boundary $\partial(\gamma_i)$ of the chain of edges $\gamma_i$ in the graph $GN$ with the vertex $v_1$, if the vertex $v_1$ appears as a summand in $\partial(\gamma_i)$, substituted by the vertex $v'_1$.
\\
\includegraphics[scale=0.5]{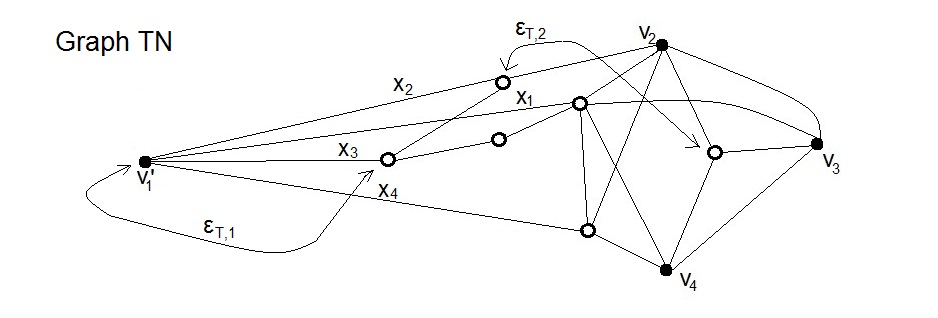}
\\ \\
One can check, that this construction of $TN$ satisfies all the requirements of our theorem. Thus, in the graph $TN$ there exist $k$ chains of edges $p_{T,1},\ldots,p_{T,k}$, such that no two of $p_{T,1},\ldots,p_{T,k}$ have any common nontrivial summands and that $\partial(p_{T,1})=\cdots=\partial(p_{T,k})=v'_1+v_2+v_{3}+v_{4}$ in $TN$.
\\ \\
Finally, we take these $k$ chains of edges $p_{T,1},\ldots,p_{T,k}$, and look at the summands in them. For each summand $x_1,\ldots,x_k$ we substitute the chain of edges $p_{L,1},p_{L,2},\ldots,p_{L,k}$ with the same index. For each summand $\epsilon_{T,1},\ldots,\epsilon_{T,\tau}$ we substitute the edge $\epsilon_{1},\ldots,\epsilon_{\tau}$ with the same index. For each summand $\epsilon_{T,\tau+1},\ldots,\epsilon_{T,t}$ we substitute the chain of edges $\gamma_{\tau+1},\gamma_{\tau+2},\ldots,\gamma_{t}$ with the same index. We obtain $k$ chains of edges $p_{1},\ldots,p_{k}$ in $GN$, and no two of $p_{1},\ldots,p_{k}$ have any common nontrivial summands. This is easy to verify, because the chains of edges $p_{L,1},p_{L,2},\ldots,p_{L,k}$ contain only the edges $x_1,\ldots,x_k$ and some edges from $K_1$, and the chains of edges $\gamma_{\tau+1},\gamma_{\tau+2},\ldots,\gamma_{t}$ contain only the edges $\epsilon_{\tau+\lambda+1},\ldots,\epsilon_n$ and some of the edges $\epsilon_{\tau+1},\ldots,\epsilon_{\tau+\lambda}$ and some of the edges from $K_1$. Finally, the boundary operation in the graph $GN$, when applied to the chains of edges $p_{L,1},p_{L,2},\ldots,p_{L,k}$, produces $\partial(p_{L,i})=v_1+u_{i,2}$ for each $1\le i\le k$. Thus $\partial(p_{1})=\cdots=\partial(p_{k})=v_1+v_2+v_{3}+v_{4}$ in $GN$.
\\ \\
Thus, we proved that $G$ is not a counter-example to our theorem, which means that in the counter-example with the minimal number of edges in it, there cannot be an edge $x_1$ such that its boundary $\partial(x_1)=u_{1,1}+u_{1,2}$ does not contain any one of the four vertices $v_1,v_2,v_3,v_4$ as a summand in it. Thus, for every edge $x$ in $G$ we have $\partial(x)=v_i+y$, where $1\le i\le k$ and $y$ is a vertex in $G$. But now Lemma \ref{minimal.final} asserts that such $G$ cannot be a counter-example. Thus, there cannot exist a counter-example to our theorem.
\end{proof}
\begin{thm}\label{Menger.homol.5} Let $G$ be a graph such that its vertex set $V$ is a union of disjoint sets $\{v_1,v_2,v_3,v_4\}$ and some empty or nonempty set $Y$, such that the degree of every vertex $y\in Y$ in the graph $G$ is even. If after a deletion of any $k-1$ or less edges from $G$ there still exists some chain of edges $p$ such that its boundary $\partial(p)$ in $G$ is $v_1+v_2+v_{3}+v_{4}$ then there exist $k$ chains of edges $p_1,\ldots,p_k$ in the graph $G$ such that no two of $p_1,\ldots,p_k$ have any common nontrivial summands and that $\partial(p_1)=\cdots=\partial(p_k)=v_1+v_2+v_{3}+v_{4}$ in $G$.
\end{thm}
\begin{proof}
Since every vertex $y\in Y$ has an even degree, one cannot draw any new edges $\epsilon_1,\ldots,\epsilon_n$ as described in Theorem \ref{Menger.homol.4}, and the graph $GN$ is just the graph $G$. Thus, Theorem \ref{Menger.homol.5} follows from Theorem \ref{Menger.homol.4}.
\end{proof}
If we take any graph $G$ which satisfies the requirements stated in Theorem \ref{Menger.homol.5}, any draw in $G$ any new edges $ne_1,\ldots,ne_n$ as described in Theorem \ref{Menger.homol.4}, we will create a new graph $GN$ which satisfies the requirements of Theorem \ref{Menger.homol.5}. Thus, Theorem \ref{Menger.homol.4} follows from Theorem \ref{Menger.homol.5}.
\\ \\
The other direction of Theorem \ref{Menger.homol.5} is obviously correct for any positive integer $n$. Namely, if in any graph $G$ there exist $k$ chains of edges $p_1,\ldots,p_k$, such that no two of $p_1,\ldots,p_k$ have any common nontrivial summands and that $\partial(p_1)=\cdots=\partial(p_k)=v_1+v_2+\cdots+v_{2n-1}+v_{2n}$ in $G$, then deleting $k-1$ or less edges from $G$ will effect at most $k-1$ of $p_1,\ldots,p_k$, hence some $p_i$ will not be effected by that deletion.
\\ \\
Our Example \ref{counter} demonstrates, that the requirement that the degrees all the vertices in $Y$ must be even in Theorem \ref{Menger.homol.5} is crucial. Indeed, deleting any one of the seven edges in this counter-example will not break the graph into separate connected components. Thus, after a deletion of any edge, there will remain some chain of edges $p$, satisfying $\partial(p)=v_1+v_2+v_3+v_4$. However, the degree of $y_1$ is $3$ and the degree of $y_2$ is $3$, and there do not exist two chains of edges $p_1$ and $p_2$ in that graph $G$, which do not have any common nontrivial summands, and which satisfy $\partial(p_1)=\partial(p_2)=v_1+v_2+v_3+v_4$.
\\ \\
Now we state the necessary and the sufficient condition for a graph $G$ to have $k$ chains of edges $p_1,\ldots,p_k$, such that no two of $p_1,\ldots,p_k$ have any common nontrivial summands and that $\partial(p_1)=\cdots=\partial(p_k)=v_1+v_2+v_{3}+v_{4}$ in $G$.
\begin{thm}\label{Menger.homol.6} Let $G$ be a graph, and let $v_1,v_2,v_3,v_4$ be four different vertices in $G$. There exist $k$ chains of edges $p_1,\ldots,p_k$ in $G$, such that no two of $p_1,\ldots,p_k$ have any common nontrivial summands and that $\partial(p_1)=\cdots=\partial(p_k)=v_1+v_2+v_{3}+v_{4}$ in $G$, if and only if after some deletion of some non-negative number of edges from $G$, every vertex $y\notin \{v_1,v_2,v_3,v_4\}$ in the modified graph $G$ has an even degree, and after a deletion of any $k-1$ or less edges from the modified graph $G$, there still remains some chain of edges $p$ there, such that its boundary $\partial(p)$ in the modified graph $G$ is $v_1+v_2+v_{3}+v_{4}$.
\end{thm}
\begin{proof}
If we can perform a deletion of some non-negative number of edges from the graph $G$ in such a way that every vertex $y\notin \{v_1,v_2,v_3,v_4\}$ in the modified graph $G$ has an even degree, and that after a deletion of any $k-1$ or less edges from the modified graph $G$, there still remains some chain of edges $p$ there, such that its boundary $\partial(p)$ in the modified graph $G$ is $v_1+v_2+v_{3}+v_{4}$, then, by Theorem \ref{Menger.homol.5}, there exist $k$ chains of edges $p_1,\ldots,p_k$ in the modified graph $G$, such that no two of $p_1,\ldots,p_k$ have any common nontrivial summands and that $\partial(p_1)=\cdots=\partial(p_k)=v_1+v_2+v_{3}+v_{4}$ in the modified graph $G$. These same $p_1,\ldots,p_k$ are also chains of edges in the original graph $G$, which satisfy in it $\partial(p_1)=\cdots=\partial(p_k)=v_1+v_2+v_{3}+v_{4}$.
\\ \\
For the other direction in the theorem, suppose that we found some $k$ chains of edges $p_1,\ldots,p_k$ in $G$, such that no two of $p_1,\ldots,p_k$ have any common nontrivial summands and that $\partial(p_1)=\cdots=\partial(p_k)=v_1+v_2+v_{3}+v_{4}$ in $G$. Delete from $G$ every edge, which does not appear as a summand in one of these $k$ chains of edges. Any vertex $y\notin \{v_1,v_2,v_3,v_4\}$ does not appear as a summand in the boundary of any of these $k$ chains of edges. Hence $y$ must appear as a summand in the boundary of an even number of edges, which appear as summands in any one of these $k$ chains of edges. Thus, the degree of every vertex $y\notin \{v_1,v_2,v_3,v_4\}$ in the modified graph $G$ will be even.
\end{proof}
We conclude our work by stating and proving our extension of the edge version of the classical Menger's Theorem using the language of paths between pairs of vertices.
\begin{thm}\label{Menger.homol.7} Let $G$ be a graph such that its vertex set $V$ is a union of disjoint sets $\{v_1,v_2,v_3,v_4\}$ and some empty or nonempty set $Y$, such that the degree of every vertex $y\in Y$ in the graph $G$ is even. If after a deletion of any $k-1$ or less edges from $G$, the four vertices $v_1,v_2,v_3,v_4$ can be split into two disjoint pairs of vertices $v_a,v_b$ and $v_c,v_d$ in such a way that there exist two paths $P_1$ and $P_2$ which do not have any common edges, such that $P_1$ goes between the vertices $v_a$ and $v_b$ and $P_2$ goes between the vertices $v_c$ and $v_d$, then the four vertices $v_1,v_2,v_3,v_4$ can be split $k$ times into two disjoint pairs of vertices $v_{a(1)},v_{b(1)}$ and $v_{c(1)},v_{d(1)}$, $\ldots$, $v_{a(k)},v_{b(k)}$ and $v_{c(k)},v_{d(k)}$, in such a way that the two vertices in each one of these $2k$ pairs can be connected by a path and no two among these $2k$ paths have any common edges.
\end{thm}
\begin{proof}
If after a deletion of any $k-1$ or less edges from $G$, the four vertices $v_1,v_2,v_3,v_4$ can be split into two disjoint pairs of vertices $v_a,v_b$ and $v_c,v_d$ in such a way that there exist two paths $P_1$ and $P_2$ which do not have any common edges, such that $P_1$ goes between the vertices $v_a$ and $v_b$, and $P_2$ goes between the vertices $v_c$ and $v_d$, then we construct the chain of edges $p$ by adding-up all the edges which appear in the paths $P_1$ and $P_2$. The requirements of Theorem \ref{Menger.homol.5} are satisfied by this $p$, hence there exist $k$ chains of edges $p_1,\ldots,p_k$ in the graph $G$ such that no two of $p_1,\ldots,p_k$ have any common nontrivial summands and that $\partial(p_1)=\cdots=\partial(p_k)=v_1+v_2+v_{3}+v_{4}$ in $G$. Our theorem now follows from applying Lemma \ref{basic} to the chains of edges $p_1,\ldots,p_k$.
\end{proof}
\section{Appendix: The Topological Viewpoint}
In this section we provide a brief introduction to the topological viewpoint of the Graph Theory. We will provide a definition of a topological space which is a geometric realization of a graph, and discuss topological properties of such topological spaces. We will also state and prove the topological versions of Theorem \ref{Menger.homol.2} and Theorem \ref{Menger.homol.5}.
\\ \\
In the literature, the geometric realizations of graphs are usually constructed either as closed subsets of the $n$-dimensional Euclidean space $\mathbb{R}^n$ or by gluing together copies of the closed interval $[0,1]$ and copies of the one-point topological space. Both of these approaches produce topological space with some additional structure.
\\ \\
For example, when one takes a copy of the closed interval $[0,1]$ for each edge of a graph, one implicitly assumes that to each point in the subset of the geometric realization which corresponds to that edge, with respect to that subset, there is an associated real number between $0$ and $1$.
\\ \\
The properties of the geometric realization are then described using these additional structures, but later are attributed to the underlying topological space. In reality, each such attribution requires a justification. For example, when we say that a topological space $X$ is homeomorphic to the closed interval $[0,1]$, we do not assume that any specific homeomorphism is given, but only that such homeomorphisms exist. Thus, for example, we cannot speak of the point in $X$ which corresponds to $0$ in $[0,1]$ or which corresponds to $0.5$ in $[0,1]$.
\\ \\
In this section we construct a geometric realization of a given graph $G$ by starting with a purely topological object, and applying to it a certain quotient map. We obtain a topological space with only one additional structure, which comes from the knowledge of the original topological object and of the quotient map. We then show that the properties which we study are independent of that structure, and as such they are properties of the underlying topological space.
\\ \\
After that we give a purely topological definition of the notion ``topological space which is a geometric realization of a graph'', and we prove that our construction of a geometric realization of a graph $G$ and this definition are compatible.
\\ \\
To make it easier for the readers who do not have a strong background in Topology and Real Analysis, we outline the proofs of several basic results which are used in our study of the geometric realization of a graph. We do assume, however, that the reader is familiar with the construction of the real numbers from the Cauchy sequences, and with the basic concepts in Topology.
\\ \\
In order to discuss various properties of the topological spaces which are geometric realizations of graphs, we first state and prove some facts about intervals of the real numbers. By a finite interval $I$ we mean the closed interval $[a,b]$, the intervals $[a,b)$ and $(a,b]$, and the open interval $(a,b)$. Here $a$ and $b$ are two real numbers, which are called the limiting points of the interval $I$, and $a<b$. By an infinite interval $I$ we mean the closed half-lines $[a,\infty)$ and $(-\infty,b]$, the open half-lines $(a,\infty)$ and $(-\infty,b)$, and the entire real line $\mathbb{R}=(-\infty,\infty)$. Unless we specify, we permit our intervals to be finite or infinite.
\begin{lem} \label{basictopol} Let $I\subset \mathbb{R}$ be an interval. Let $f:I\rightarrow\{0,1\}$ be a function from $I$ to the set $\{0,1\}$. If there exist $x,y\in I$ such that $f(x)=0$ and $f(y)=1$ then for some $\min(x,y)\le z\le\max(x,y)$ in $I$, in every open neighborhood of $z$ in $I$ there are points $a$ and $b$ such that $f(a)=0$ and $f(b)=1$.
\end{lem}
\begin{proof} We apply \textit{the Bolzano$-$Weierstrass (the B$-$W) argument}. We take the midpoint $z_1$ between $x_1=x$ and $y_1=y$ and check if the function $f$ at $z_1$ returns $0$, like at $x$, or $1$, like at $y$. According to how $f$ behaves at $z_1$, we take this $z_1$ as $x_2$ and keep $y_2=y$, or take this $z_1$ as $y_2$ and keep $x_2=x$.
\\ \\
Repeating this process again and again produces a Cauchy sequence of midpoints, which by the construction of the real numbers must converge to some $z$ between $\min(x,y)$ and $\max(x,y)$. If this $z$ is one of the $x_i$ then $f(z)=0$, and every open neighborhood of $z$ in $I$ contains some midpoint $z_j=y_{j+1}$ where $f(y_{j+1})=1$. If this $z$ is one of the $y_j$ then $f(z)=1$, and every open neighborhood of $z$ in $I$ contains some midpoint $z_i=x_{i+1}$ for which $f(x_{i+1})=0$. Otherwise, every open neighborhood of $z$ in $I$ contains some midpoints $z_i=x_{i+1}$ for which $f(x_{i+1})=0$ and contains some midpoints $z_j=y_{j+1}$ for which $f(y_{j+1})=1$.
\end{proof}
We list here several important consequences of Lemma \ref{basictopol}. Any continuous function from an interval $I$ to the set $\{0,1\}$ equipped with the discrete topology must take the entire $I$ either to $0$ or to $1$. Indeed, otherwise the point $z$, described in Lemma \ref{basictopol}, cannot belong to the pre-image of $0$ and cannot belong to the pre-image of $1$, because these two pre-images must be open in $I$. It follows that an interval $I$ is not equal to a disjoint union of two nonempty intervals which are open in $I$, and that a continuous function from $I$ to a disjoint union of two copies of $\mathbb{R}$ must take the entire $I$ to one of these two copies.
\\ \\
Any closed interval $[x,y]$ is compact. Indeed, if we are given a covering of $[x,y]$ by open sets, we can define a function $f:[x,y]\rightarrow\{0,1\}$, where $\{0,1\}$ is equipped with the discrete topology, by $f(z)=0$ if we can cover $[x,z]$ by a finite sub-covering, and $f(z)=1$ if we cannot cover $[x,z]$ by a finite sub-covering. Here $[x,x]$ is just $\{x\}$. Since for every $z\in [x,y]$, the given covering has some open set $U$ containing $z$, this $f$ is continuous. Indeed, if we can cover $[x,z]$ by a finite sub-covering then we can cover by a finite sub-covering $[x,u]$ for all $u\in U$, and if we cannot cover $[x,z]$ by a finite sub-covering then we cannot cover by a finite sub-covering $[x,u]$ for all $u\in U$. Since $f(x)=0$, $f$ takes the entire $[x,y]$ to $0$. Thus, $f(y)=0$ and $[x,y]$ can be covered by a finite sub-covering.
\\ \\
Any infinite sequence $a_1,a_2,\ldots$ of points in any closed interval $[x,y]$ must have a converging subsequence. In other words, there must be a point $z\in [x,y]$ such that any open neighborhood of $z$ in $[x,y]$ contains an infinite number of points from the sequence $a_1,a_2,\ldots$. Indeed, otherwise every point $z\in [x,y]$ will have some open neighborhood $U$ in $[x,y]$ which contains only a finite number of points from $a_1,a_2,\ldots$, and $U$ without all these points is still an open subset of $[x,y]$. Thus, the sequence of open sets $[x,y]-\{a_n,a_{n+1},\ldots\}$, where $n=1,2,\ldots$, is an infinite covering of $[x,y]$, and it does not have a finite sub-covering.
\\ \\
For any continuous function $f:[x,y]\rightarrow \mathbb{R}$ there exist points $m,M\in \mathbb{R}$ such that $m\le f(z)\le M$ for all $z\in [x,y]$. Indeed, otherwise we can construct an infinite sequence $a_1,a_2,\ldots$ of points in $[x,y]$ such that $f(a_1),f(a_2),\ldots$ increases to $\infty$ or decreases to $-\infty$. But $a_1,a_2,\ldots$ must have a convergent subsequence in $[x,y]$, and the function $f$ is not continuous at its limit point. Notice, that we can have $m=M$, in which case $f$ is a constant function.
\\ \\
For any continuous function $f:[x,y]\rightarrow \mathbb{R}$ there exist points $u,w\in [x,y]$ such that $f(u)\le f(z)\le f(w)$ for all $z\in [x,y]$. Indeed, we can find two real numbers $\alpha,\beta\in\mathbb{R}$, such that $\alpha$ is in the image $f([x,y])$ of $[x,y]$ under $f$, and $\beta$ is greater than any number in $f([x,y])$. Applying the B$-$W argument to the interval $[\alpha,\beta]$ produces some real number $\alpha\le M\le\beta$ such that in any open neighborhood of $M$ in $\mathbb{R}$ there is a number in $f([x,y])$ and a number greater than any number in $f([x,y])$. Thus, $f(z)$ cannot be greater than $M$ for all $z\in [x,y]$. On the other hand, we can find an infinite sequence $a_1,a_2,\ldots$ of points in $[x,y]$ such that $f(a_1),f(a_2),\ldots$ converges to $M$. The limit point $w$ of any convergent subsequence of $a_1,a_2,\ldots$ in $[x,y]$ must be taken by $f$ to $M$. Repeating this argument for $\beta$ smaller than any number in $f([x,y])$ produces the required $u$. We may have $u=w$, or we may have $u\ne w$ but $f(u)=f(w)$, and in these cases $f$ is a constant function.
\\ \\
An image of a closed interval $[x,y]$ under a continuous function $f:[x,y]\rightarrow\mathbb{R}$ is either a one-point subset $\{m\}\subset \mathbb{R}$ or a closed interval $[m,M]\subset \mathbb{R}$. Indeed, there must exist points $u,w\in [x,y]$ such that $m=f(u)\le f(z)\le f(w)=M$ for all $z\in [x,y]$. If $m=M$ then $f(z)=m$ for all $z\in[x,y]$. If $m<M$ then for any real number $m<t<M$ we start with the $u$ and $w$, and we take the midpoint $r$ between $u$ and $w$. If $f(r)=t$ we are done. Otherwise, if $f(r)<t$ we use $r$ as our new $u$ and if $f(r)>t$ we use $r$ as our new $w$. Repeating this process again and again generates a Cauchy sequence of midpoints, and in any open neighborhood of the limit point of that Cauchy sequence there are points $\alpha$ for which $f(\alpha)<t$ and points $\beta$ for which $f(\beta)>t$. Thus, $f$ takes the limit point of that Cauchy sequence to $t$.
\\ \\
We are now ready to state and prove the Invariance of Domain Theorem for the one-dimensional Eucledian space $\mathbb{R}^1$. The general case of the Invariance of Domain Theorem for the $n$-dimensional Eucledian space  $\mathbb{R}^n$ is stated and proved in Chapter 2, Paragraph 9, of \cite{Alexandrov}. There is a slight inaccuracy in the proof of Lemma 1 there, because after the small translation of the closed $n$-dimensional simplex $T$, the points which belong to $n+1$ different cubes of the $\epsilon$-covering and which lie on one of the sides of the boundary $\tilde{T}$ of the simplex $T$, may end-up on another side of $\tilde{T}$. This inaccuracy can be easily corrected by imposing an additional requirement on how small the small translation, which is described there, must be.
\begin{lem} \label{InvDom} An injective continuous function $f:\mathbb{R}\rightarrow \mathbb{R}$ always takes $\mathbb{R}$ onto an open interval.
\end{lem}
\begin{proof}
Since we can compose $f$ with the homeomorphism $\arctan:\mathbb{R}\rightarrow (-\frac{\Pi}{2},\frac{\Pi}{2})$, we can assume that $f$ takes $\mathbb{R}$ to the open interval $(-\frac{\Pi}{2},\frac{\Pi}{2})\subset \mathbb{R}$.
\\ \\
For $n=1,2,\ldots$, let the image of the closed interval $[-n,n]$ under $f$ be the closed interval $[m(n),M(n)]$ inside $(-\frac{\Pi}{2},\frac{\Pi}{2})$. Since $f$ is injective, the image of $[-(n+1),-n]$ under $f$ must either be the closed interval $[m(n+1),m(n)]$ or the closed interval $[M(n),M(n+1)]$, and the image of $[n,(n+1)]$ under $f$ must either be the closed interval $[M(n),M(n+1)]$ or the closed interval $[m(n+1),m(n)]$, respectively. We see that $m(1)>m(2)>\cdots$ and that $M(1)<M(2)<\cdots$.
\\ \\
Since the image of $\mathbb{R}$ under $f$ is bounded between $-\frac{\Pi}{2}$ and $\frac{\Pi}{2}$, the sequences $m(1),m(2),\ldots$ and $M(1),M(2),\ldots$ must have some converging subsequences, and we can find some converging subsequence of $m(1),m(2),\ldots$ which converges to some $m$ from above, and some converging subsequence of $M(1),M(2),\ldots$ which converges to some to $M$ from below.
\\ \\
For each open neighborhood $(m-\epsilon,m+\epsilon)$ of $m$ in $\mathbb{R}$ we can find some $m(t)>m$ from the converging subsequence in that neighborhood, and this implies that all the $m(t+1),m(t+2),\ldots$ are also in that neighborhood of $m$. Thus, the entire sequence $m(1),m(2),\ldots$ converges to $m$ from above. By the same argument the sequence $M(1),M(2),\ldots$ converges to $M$ from below. Hence, the open interval $(m,M)$ is contained in the image of $\mathbb{R}$ under $f$.
\\ \\
On the other hand, for any real number $x$ we can find some integer $n$ such that $x\in [-n,n]$. Thus, $m<m(n)<f(x)<M(n)<M$, which means that the image of $\mathbb{R}$ under $f$ is contained inside the open interval $(m,M)$. Hence the image of $\mathbb{R}$ under $f$ is equal to the open interval $(m,M)$.
\end{proof}
It follows from Lemma \ref{InvDom} that any injective continuous function $f$ from an open interval $I$ to $\mathbb{R}$ is a homeomorphism between $I$ and the image of $I$ under $f$. Indeed, we only need to show that the inverse $f^{-1}$ of $f$ is a continuous function from the image of $I$ under $f$ to $I$, and that continuity follows from the fact that the pre-images of the open intervals under $f^{-1}$ are the images of these open intervals under $f$, and these images will be open intervals by Lemma \ref{InvDom}.
\\ \\
It also follows from Lemma \ref{InvDom} that the point $a$ in the interval $I$, which is $[a,b]$, or $[a,b)$, or $(b,a]$, or $[a,\infty)$, or $(-\infty,a]$, does not have any neighborhood in $I$ which is homeomorphic to $\mathbb{R}$. Indeed, such a neighborhood must be some open interval $(m,M)\subset I$, which implies that either $a\le m<M$ or $m<M\le a$, but such an open interval does not contain $a$.
\\ \\
We define \textit{the inner points in a topological space $X$ which is homeomorphic to an interval $I$} as the points in $X$ which have some open neighborhood in $X$ which is homeomorphic to $\mathbb{R}$. It follows from Lemma \ref{InvDom} that all such neighborhoods in $I$ are open intervals $(m,M)\subset I$. Notice, that any homeomorphism between topological spaces $x$ and $Y$ which are homeomorphic to intervals takes the inner points in $X$ to the inner points in $Y$ and the non-inner points in $X$ to the non-inner points in $Y$.
\\ \\
For $\mathbb{R}^{\infty}$, which is the set of all the sequences of real numbers $(x_1,x_2,\ldots)$ with only a finite number of nonzero entries, equipped with the Euclidean metric, these results no longer hold. The injective continuous function $f:\mathbb{R}^{\infty}\rightarrow\mathbb{R}^{\infty}$ which takes each $(x_1,x_2,\ldots)\in \mathbb{R}^{\infty}$ to $(0,x_1,x_2,\ldots)\in \mathbb{R}^{\infty}$ is a homeomorphism between $\mathbb{R}^{\infty}$ and the image of $\mathbb{R}^{\infty}$ under $f$, but the image of $\mathbb{R}^{\infty}$ under $f$ is not open in $\mathbb{R}^{\infty}$.
\\ \\
The injective continuous function $h:\mathbb{R}^{\infty}\rightarrow\mathbb{R}^{\infty}$ which takes each $(x_1,x_2,x_3,\ldots) \in \mathbb{R}^{\infty}$ to $(h_1(x_1),h_2(x_1),x_2,x_3,\ldots)\in \mathbb{R}^{\infty}$, where $h_1(x_1)=x_1$ when $x_1\le 1$, $h_1(x_1)=1$ when $1\le x_1\le 2$, $h_1(x_1)=3-x_1$ when $2\le x_1\le 3$, $h_1(x_1)=0$ when $3\le x_1$, and $h_2(x_1)=0$ when $x_1\le 1$, $h_2(x_1)=x_1-1$ when $1\le x_1\le 2$, $h_2(x_1)=1$ when $2\le x_1\le 3$, $h_2(x_1)=\frac{3}{x_1}$ when $3\le x_1$, is not a homeomorphism between $\mathbb{R}^{\infty}$ and the image of $\mathbb{R}^{\infty}$ under $h$.
\\ \\
From Lemma \ref{InvDom} it also follows, that the image of a disjoint union of two copies $R_1$ and $R_2$ of $\mathbb{R}$ under a continuous injective function to an interval $I$ cannot contain $I$. Indeed, otherwise for each $z\in I$ we can define $f(z)=0$ if $z$ is in the image of $R_1$ and $f(z)=1$ if $z$ is in the image of $R_2$. Since these two images are open in $I$, $f$ is a continuous function from $I$ to the set $\{0,1\}$ equipped with the discrete topology. Thus, either every $z$ in $I$ is in the image of $R_1$ or every $z$ in $I$ is in the image of $R_2$, but an image of a copy of $\mathbb{R}$ under any function cannot be an empty set. Hence there exists some $z\in I$ not in the image of $R_1\cup R_2$.
\\ \\
Let $X$ be a topological space, and let $\sim$ be an equivalence relation between the points in $X$ which breaks $X$ into equivalence classes of points. Let $Y$ be the set of all these equivalence classes, and let the quotient map $q:X\rightarrow Y$ take each $x\in X$ to the equivalence class $y\in Y$ to which $x$ belongs. The topological space $Y$ equipped with the topology in which a set $U\in Y$ is open if and only if its pre-image $q^{-1}(U)$ under $q$ is open in $X$ is called \textit{the quotient space of $X$ under the equivalence relation $\sim$}.
\\ \\
Let $h_1,\ldots,h_t$ be homeomorphisms between a topological space $Y$ and some closed subsets of topological spaces $X_1,\ldots,X_t$, respectively. By \textit{gluing $X_1,\ldots,X_t$ along $h_1(Y),\ldots,h_t(Y)$} we mean taking the disjoint union of $X_1,\ldots,X_t$, introducing in it the equivalence relation under which $h_i(y)$ is equivalent to $x_j(y)$ for all $1\le i,j\le t$ and for all $y\in Y$, and then taking the quotient space of this disjoint union under that equivalence relation. When $Y$ consists of a finite number of points we just say that we glue all the $X_1,\ldots,X_t$ at their such and such points.
\\ \\
The topological space $X$, constructed by taking topological spaces $H_1$ and $H_2$ homeomorphic to the closed half-line $[0,\infty)$, and gluing them at their non-inner points, is homeomorphic to $\mathbb{R}$. Indeed, let $q$ be a quotient map from the disjoint union of $H_1$ and $H_2$ onto $X$, $\phi_1$ be a homeomorphism between $H_1$ and $[0,\infty)$, and $\phi_2$ be a homeomorphism between $H_2$ and $[0,\infty)$. We define the homeomorphism $\phi$ between $X$ and $\mathbb{R}$ as follows: If $q^{-1}(x)$ consists of only one point $z\in H_1$ then $\phi(x)=\phi_1(z)$; If $q^{-1}(x)$ consists of only one point $z\in H_2$ then $\phi(x)=\phi_2(z)$; If $q^{-1}(x)$ consists of $0\in H_1$ and $0\in H_2$ then $\phi(x)=0$. It is easy to verify that $\phi$ is a homeomorphism.
\bde For any nonnegative integer $d$ we define the model space $\Xi_d$ as $d$ copies of the closed half-line $[0,\infty)$ all glued together at their points $0$. For $d=0$ the model space $\Xi_d$ is just the one-point space $\{0\}$. We refer to the common $0$ point of these $d$ copies of $[0,\infty)$ as the point $0$ in $\Xi_d$.
\ede
For each $n=1,3,\ldots$ let $W_{d,n}$ be the subset of $\Xi_d$, consisting of the point $0$ and of all the inner points $z$ inside each one of the $d$ copies of the closed half-line $[0,\infty)$, satisfying $z<\frac{1}{n}$. The local base of topology of $\Xi_d$ at the point $0$ is given by the sets $W_{d,n}$ with $n=2,3,\ldots$. It can be easily verified that for each set $W_{d,n}$ there exists a homeomorphism between $W_{d,n}$ and $\Xi_d$, which takes the common $0$ point in $W_{d,n}$ to the point $0$ in $\Xi_d$.
\\ \\
A direct consequence of Lemma \ref{InvDom} is that there is no injective continuous function $f$ from $\Xi_d$ with $d\ge 3$ into an interval $I$. Indeed, $f$ would have to take two out of the $d$ copies of $[0,\infty)$ glued together at their points $0$ onto some open interval $(m,M)\subset I$ which contains the point $f(0)$ inside it. Thus, the pre-image of that open interval $(m,M)$ under $f$ would not be open in $\Xi_d$.
\begin{lem} \label{topdegree1} For $n>\max(m,2)$ there is no continuous injective function\\ $f:\Xi_n\rightarrow \Xi_m$.
\end{lem}
\begin{proof}
If $m=0$ then the $\Xi_m$ is just the one-point space, and the lemma trivially follows. Thus, assume that $m\ge 1$.
\\ \\
The local base of topology of $\Xi_n$ at the point $0$ is given by the open neighborhoods which are homeomorphic to $\Xi_n$. Thus, if we found a continuous injective function $f:\Xi_n\rightarrow\Xi_m$, then, since there is no continuous injective function from $\Xi_n$ with $n\ge 3$ into an interval, $f$ would have to take the point $0$ in $\Xi_n$ to the the point $0$ in $\Xi_m$.
\\ \\
Hence, the image under $f$ of each one of the $n$ copies of the open half-line $(0,\infty)$ which are inside the $n$ copies of the closed half-line $[0,\infty)$, which were glued together at their points $0$ to form $\Xi_n$, is some open interval of the form $(a,b)$ or $(a,\infty)$ inside one of the $m$ copies of the closed half-line $[0,\infty)$, which were glued together at their points $0$ to form $\Xi_m$.
\\ \\
Since $n>m$, two different copies $A$ and $B$ of $(0,\infty)$, contained inside two different copies of $[0,\infty)$ from among the $n$ copies of the $[0,\infty)$ glued together to create $\Xi_n$, must be taken by $f$ into the same copy of $(0,\infty)$ contained inside one of the $m$ copies of $[0,\infty)$ glued together to create $\Xi_m$. But this implies that the point $0$ in $\Xi_m$ has some open neighborhood $U$ which does not intersect either $f(A)$ or $f(B)$ or both of them. The pre-image of $U$ under $f$ is not open in $\Xi_n$, which means that $f$ is not continuous at $0$.
\end{proof}
\begin{lem} \label{topdegree2} For $m\ne n$, the point $0$ in $\Xi_n$ does not have any open neighborhood $U$ in $\Xi_n$ homeomorphic to $\Xi_m$.
\end{lem}
\begin{proof}
If $m=0$ or $n=0$ then the corresponding topological space is the one-point space, and our lemma is trivial. Thus, we assume that $m\ge 1$ and $n\ge 1$.
\\ \\
The local base of topology of $\Xi_n$ at the point $0$ is given by the open neighborhoods which are homeomorphic to $\Xi_n$. Thus, a homeomorphism $f$ between $\Xi_m$ and any open neighborhood $U$ of the point $0$ in $\Xi_n$ must induce a homeomorphism $h$ between $\Xi_n$ and some open subset $W$ of $\Xi_m$. Indeed, there exists some open neighborhood $V$ of the point $0$ in $\Xi_n$, such that $V$ is homeomorphic to $\Xi_n$ and contained in $U$, and the inverse of $f$, restricted to $V$, is a homeomorphism between $V$ and some open subset $W$ of $\Xi_m$. Since $m\ne n$, it follows from Lemma \ref{topdegree1} that for such $f$ and $h$ to exist, we must have $n\le 2$ and $m\le 2$.
\\ \\
Thus, we only have to consider two cases: $n=2, m=1$ and $n=1, m=2$. For $n=1, m=2$, Lemma \ref{InvDom} asserts that the image under $f$ of $\Xi_2$, which is homeomorphic to $\mathbb{R}$, in $[0,\infty)$ is an open interval, and so it cannot contain the point $0$. For $n=2, m=1$, we notice that $f(0)$ is a non-inner point in the image of $\Xi_1$ under $f$, but any open neighborhood of $f(0)$ in $\Xi_2$ contains an open neighborhood of $f(0)$ which is homeomorphic to $\mathbb{R}$. Thus, the image of $\Xi_1$ under $f$ does not contain any open neighborhood of $f(0)$ in $\Xi_2$, which implies that this image is not open in $\Xi_2$.
\end{proof}
Now we discuss two approaches to constructing a geometric realization of a graph $G$. In both approaches we end-up with an object, which has more structure than just topology. However, the underlying topological space in both approaches is the same up to a homeomorphism.
\\ \\
Let $G$ be a graph. For each edge $e$ in $G$ we define $K(e)$ to be a copy of the closed interval $[0,1]$, and we define $T(e)$ to be a topological space homeomorphic to the closed interval $[0,1]$. For each vertex $v$ in $G$ we define $K(v)$ and $T(v)$ to be copies of the one-point topological space.
\\ \\
We remind that when we say that an object $X$ is a copy of an object $Y$, we mean that together with the objects $X$ and $Y$ we are also given some specific isomorphism between $X$ and $Y$, and under that isomorphism all the relevant properties of the objects $X$ and $Y$ are identical. However, when we just know that there exists a homeomorphism between topological spaces $X$ and $Y$, we are not actually given any specific homeomorphism. Thus, all the statements, which we make, must be independent of a particular choice of a homeomorphism between $X$ and $Y$.
\\ \\
Every homeomorphism $h:T(e)\rightarrow [0,1]$ takes inner points in $T(e)$ to the inner points $0<z<1$ in $[0,1]$, and takes the two non-inner points in $T(e)$ to the non-inner points $0$ and $1$ in $[0,1]$. Regardless of our selection of $h$, the local base of topology of $T(e)$ at any inner point $x\in T(e)$ is given by the open neighborhoods $h^{-1}(U_i)\subset T(e)$ of $x$ in $T(e)$, where $$U_i=(h(x)-\frac{1}{i},h(x)+\frac{1}{i})\subset [0,1]$$ with $i=n,n+1,\ldots$. Here $n$ is an integer satisfying $$n>\max(\frac{1}{h(x)},\frac{1}{1-h(x)})$$
\\
Each open neighborhood $h^{-1}(U_i)$ of $x$ in $T(e)$ is homeomorphic to the open interval $U_i$, and $U_i$ is homeomorphic to a copy of the interval $U_{1,i}=(h(x)-\frac{1}{i},h(x)]$ and a copy of the interval $U_{2,i}=[h(x),h(x)+\frac{1}{i})$ glued together at their points $h(x)$. Thus each $h^{-1}(U_i)$ is homeomorphic to $h^{-1}(U_{1,i})$ and $h^{-1}(U_{2,i})$ glued together at their points $x$.
\\ \\
We define the total topological space $Tot$ as the disjoint union of all the spaces $T(e)$, where $e$ is an edge in $G$, and all the spaces $T(v)$, where $v$ is a vertex in $G$. For each edge $e$ in $G$, the boundary $\partial(e)$ of $e$ in $G$ is a set consisting of two different vertices $v_1$ and $v_2$, and in each topological space $T(e)$ there are two non-inner points. To one of these two non-inner points in $T(e)\subset Tot$ we associate the point $T(v_1)\subset Tot$, and to the other one of them we associate the point $T(v_2)\subset Tot$.
\\ \\
After we perform this association, we introduce an equivalence relation in $Tot$ by saying that $x_1\in T(e_1)$ is equivalent to $x_2\in T(e_2)$ if the same point $T(v)$ is associated to $x_1$ and to $x_2$, and saying that $x\in T(e)$ is equivalent to $T(v)$ if $T(v)$ is associated to $x\in T(e)$. The quotient space $T(G)$ of $Tot$ under this equivalence relation is called \textit{the geometric realization of the graph $G$}. The image of $T(e)$ in $T(G)$ is called \textit{the closed subset of $T(G)$ which corresponds to the edge $e$}. The quotient map restricted to $T(e)$ is a homeomorphism between $T(e)$ and its image in $T(G)$, hence we abuse the notation and refer to the image of $T(e)$ in $T(G)$ as $T(e)$. The image of $T(v)$ in $T(G)$ is called \textit{the point in $T(G)$ which corresponds to the vertex $v$}. Abusing the notation, we refer to it as $T(v)$.
\\ \\
When we associated each one of the two vertices in $\partial(e)$ to one of the two non-inner points in $T(e)$, we had to make choices. If we would choose differently, we would obtain a topological space which is homeomorphic to the one we constructed. Indeed, let $h:T(e)\rightarrow [0,1]$ be any homeomorphism, and define $g(x)=h^{-1}(1-h(x))$. Then $g:T(e)\rightarrow T(e)$ is a homeomorphism which interchanges the two non-inner points in $T(e)$. The map $g$ extends to homeomorphism $g_e:Tot\rightarrow Tot$, defined as the identity map on all $T(e')$, when $e'\ne e$, and on all $T(v)$, and defined as $g$ on $T(e)$. Finally, $g_e$ induces a homeomorphism between $T(G)$ constructed by making one choice of associating two vertices in $\partial(e)$ to the two non-inner points in $T(e)$, and $T(G)$ constructed by making the other choice.
\\ \\
The object $T(G)$, which we constructed, has more structure than just topology. For each point $x\in T(G)$ we know the pre-images of $x$ in $Tot$ under the quotient map, and this, in general, imposes an additional structure on $T(G)$.
\\ \\
For example, if $G$ has three vertices $v_1,v_2,v_3$ and two edges $e_1,e_2$ such that $\partial(e_1)=\{v_1,v_2\}$ and $\partial(e_2)=\{v_2,v_3\}$, then $Tot$ is the disjoint union of the three copies $T(v_1)$, $T(v_2)$ and $T(v_3)$ of the one-point space and of the two topological spaces $T(e_1)$ and $T(e_2)$ homeomorphic to $[0,1]$. Thus, $T(G)$ is obtained as the quotient of $Tot$ by the equivalence relation, which relates one of the non-inner points in $T(e_1)$ to $T(v_1)$, and relates the other non-inner point in $T(e_1)$ to $T(v_2)$ and to one of the non-inner points in $T(e_2)$, and relates the other non-inner point in $T(e_2)$ to $T(v_3)$. However, as a topological space, $T(G)$ is homeomorphic to $[0,1]$, so it has two non-inner points, and all the other points in it are inner points and are topologically indistinguishable from one-another.
\\ \\
Now we introduce a $1$-dimensional simplicial complex $K(G)$ which is also called a geometric realization of the graph $G$.
\\ \\
We define the total complex $Kt$ to be the disjoint union of all the $K(e)$, where $e$ is an edge in $G$, and all the $K(v)$, where $v$ is a vertex in $G$. Next, for each edge $e$ in $G$, the boundary $\partial(e)$ of $e$ in $G$ is a set consisting of two different vertices $v_1$ and $v_2$, and we associate $v_1$ to the point $0$ in $K(e)\subset Kt$ and $v_2$ to the point $1$ in $K(e)\subset Kt$. Next, we introduce an equivalence relation in $Kt$ by saying that $x_1\in K(e_1)$ is equivalent to $x_2\in K(e_2)$ if the same $K(v)$ is associated to $x_1$ and $x_2$, and saying that $x\in K(e)$ is equivalent to $K(v)$ if $K(v)$ is associated to $x\in K(e)$. The quotient $K(G)$ of $Kt$ under that equivalence relation is also called the geometric realization of the graph $G$.
\\ \\
The images in $K(G)$ of $K(e)$ and of $K(v)$ under the quotient map are also called $K(e)$ and $K(v)$. Simplicial complex $K(G)$ comes with more structure than $T(G)$ because we can speak of coordinates of points in each $K(e)\subset K(G)$. However, the same point $x\in K(G)$ can belong to $K(e_1)$ and to $K(e_2)$ with $e_1\ne e_2$, and $x$ can have coordinate $0$ in $K(e_1)$ and coordinate $1$ in $K(e_2)$.
\\ \\
When we speak of a topological space $X$ which is a geometric realization of a graph $G$, we can start with $T(G)$ or $K(G)$, and then speak of a topological space $X$ which is homeomorphic to $T(G)$ or to $K(G)$.
\\ \\
We are now going to associate a non-negative number to each point $x\in T(G)$, and show that this number is a local topological invariant, which does not depend on the additional structure of $T(G)$.
\\ \\
Consider any point $x=T(v)$ in $T(G)$. Let $e_1,\ldots,e_d$, where $d$ is a nonnegative integer, be all the edges in $G$ such that $T(v)\subset Tot$ is associated to one of the two non-inner points in $T(e)\subset Tot$. If $d=0$ then $x$ is an isolated point in $T(G)$, which means that $\{x\}\subset T(G)$ is an open neighborhood of $x$ in $T(G)$. If $d>0$ then let $$h_1:T(e_1)\rightarrow [0,1]$$ $$\ldots$$ $$h_d:T(e_d)\rightarrow [0,1]$$ be any homeomorphisms which take the non-inner points $\mu_i\in T(e_i)$, to which $T(v)$ is associated, to $0\in [0,1]$. For each $n=2,3,\ldots$ let $W_{d,n}$ be the subset of the $d$ copies of the interval $[0,1]$ all glued together at their points $0$, consisting of the common point $0$ and of all the inner points $z$ inside each one of the $d$ copies of the interval $[0,1]$, such that $z<\frac{1}{n}$. The local base of topology of the $d$ copies of the interval $[0,1]$ all glued together at their points $0$, at the common point $0$, is given by the open subsets $W_{d,n}$ with $n=2,3,\ldots$.
\\ \\
The homeomorphisms $h_1^{-1},\ldots,h_d^{-1}$ induce a homeomorphism $g$ between $W_{d,2}$ and the image of $W_{d,2}$ in $T(G)$, which takes the common point $0$, which is a point in $W_{d,2}$, to $x$ in the image of $W_{d,2}$. The images of the sets $W_{d,n}$ with $n=2,3,\ldots$ under $g$ give the local base of topology of $T(G)$ at $x$, and this statement is independent of the choices of the homeomorphisms $h_1,\ldots,h_d$.
\\ \\
Thus, for any vertex $v$ in $G$, we can describe the point $x=T(v)\in T(G)$ as a point $x$ for which there exist an open neighborhood $\Omega$ in $T(G)$ and a homeomorphism $\phi$ between that $\Omega$ and $\Xi_d$ with $d\ge 0$, which satisfies $\phi(x)=0$. From our discussion above it is evident that the same description applies to any inner point $x$ in $T(e)$, where $e$ is any edge in $G$, and for all of them we have $d=2$. From Lemma \ref{topdegree2} it follows, that this description, applied to any $x\in T$, where $T$ is a topological space homeomorphic to $T(G)$, produces a unique nonnegative integer $d$. This integer is called \textit{the degree of $x$ in $T$}.
\\ \\
Now we give a local-topological definition of a topological space which is a geometric realization of a graph. We remind, that a topological space is called Hausdorff if any two different points in that space have disjoint open neighborhoods in it:
\bde\label{defTopolGraph} A compact Hausdorff topological space $T$ is called a geometric realization of a graph if for every point $x\in T$ there exist a nonnegative integer $d$, an open neighborhood $\Omega$ of $x$ in $T$, and a homeomorphism $\phi$ between $\Omega$ and $\Xi_d$ which satisfies $\phi(x)=0$.
\ede
The requirement that every $x\in T$ has an open neighborhood in $T$, which is homeomorphic to $\Xi_d$, is sufficient for $T$ to be a Fech\'{e}t space, which means that every $\{x\}\subset T$ is closed in $T$, but is not sufficient for $T$ to be Hausdorff. Let, for example, $T$ be two copies of the open interval $(-1,1)$, in which every point $x\ne 0$ in the first copy of $(-1,1)$ is glued to that same $x$ in the second copy of $(-1,1)$. The point $0_1\in T$ coming from the first copy has an open neighborhood homeomorphic to $\Xi_2$, the point $0_2\in T$ coming from the second copy has an open neighborhood homeomorphic to $\Xi_2$, and all the other points in $T$ have open neighborhoods homeomorphic to $\Xi_2$. However, any open neighborhood of $0_1$ intersects any open neighborhood of $0_2$.
\\ \\
Intuitively, we took a topological space for which, like for any $T(G)$, its points can be separated by disjoint open sets, and any covering by the open sets has a finite sub-covering, and we imposed a local property on it, requiring it at every point to be locally homeomorphic to some $T(G)$. To justify the name ``a topological space which is a geometric realization of a graph'' we must show that any such space is homeomorphic to some $T(G)$. We will prove that in Theorem \ref{maintopology}, but first we need six technical lemmas.
\begin{lem}\label{topolinter} Let $T$ be a topological space, $V$ be an open subset of $T$, and $f:U\rightarrow I$, where $I$ is a finite interval, be a homeomorphism. If $U\cap V\ne \emptyset$ then $f(U\cap V)$ is a disjoint union of intervals in $I$, which are all open in $I$.
\end{lem}
\begin{proof}
Let $I$ be $(m,M)$ or $[m,M)$ or $(m,M]$ or $[m,M]$, and let $\alpha$ be a point in $f(U\cap V)$. If there is no $\beta_2>\alpha$ in $I$ such that $\beta_2\notin f(U\cap V)$ then $[\alpha,M]$, if $M\in I$, or $[\alpha,M)$, if $M\notin I$, is contained in $f(U\cap V)$. In this case we have two subcases:
\\ \\ 
Subcase 1: There is no $\beta_1<\alpha$ in $I$ such that $\beta_1\notin f(U\cap V)$. Then $[m,\alpha]$, if $m\in I$, or $(m,\alpha]$, if $m\notin I$, is contained in $f(U\cap V)$. In this subcase $f(U\cap V)=I$.
\\ \\
Subcase 2: There is some $\beta_1<\alpha$ in $I$ such that $\beta_1\notin f(U\cap V)$. We define the function $h:[\beta_1,\alpha]\rightarrow \{0,1\}$ as $h(t)=0$ if $[t,\alpha]$ is not a subset of $f(U\cap V)$, and $h(t)=1$ if $[t,\alpha]$ is a subset of $f(U\cap V)$. Here $[\alpha,\alpha]=\{\alpha\}$. By Lemma \ref{basictopol}, there exists some $z_1\in [\beta_1,\alpha]$ such that in any open neighborhood of $z_1$ in $[\beta_1,\alpha]$ there are points $t$ for which $h(t)=0$ and points $t$ for which $h(t)=1$. Since $h(t_2)\ge h(t_1)$ when $t_2>t_1$, $h(t)=0$ for all $t<z_1$ and $h(t)=1$ for all $t>z_1$. Since $U\cap V$ is open in $U$, $f(U\cap V)$ is open in $I$. Thus $h(z_1)=0$. In this subcase $(z_1,M]$ or $(z_1,M)$, depending on if $M\in I$ or $M\notin I$, is contained in $f(U\cap V)$, but $z_1\notin f(U\cap V)$.
\\ \\
If there is some $\beta_2>\alpha$ in $I$ such that $\beta_2\notin f(U\cap V)$ then we define the function $h:[\alpha,\beta_2]\rightarrow \{0,1\}$ by $h(t)=0$ if $[\alpha,t]$ is a subset of $f(U\cap V)$ and $h(t)=1$ if $[\alpha,t]$ is not a subset of $f(U\cap V)$. Repeating the above argument, we obtain some $z_2\in [\alpha,\beta_2]$ such that $h(t)=0$ for all $t<z_2$ and $h(t)=1$ for all $t>z_2$. Since $f(U\cap V)$ is open in $I$, $h(z_2)=1$. Thus, $[\alpha,z_2)\subset f(U\cap V)$, but $z_2\notin f(U\cap V)$. In this case we also have two subcases:
\\ \\
Subcase 1: There is no $\beta_1<\alpha$ in $I$ such that $\beta_1\notin f(U\cap V)$. Then $[m,\alpha]$, if $m\in I$, or $(m,\alpha]$, if $m\notin I$, is contained in $f(U\cap V)$. In this subcase $[m,z_2)$ or $(m,z_2)$, depending on if $m\in I$ or $m\notin I$, is contained in $f(U\cap V)$, but $z_2\notin f(U\cap V)$.
\\ \\
Subcase 2: There is some $\beta_1<\alpha$ in $I$ such that $\beta_1\notin f(U\cap V)$. Repeating the argument from the Subcase 2 above, we produce some $z_1$ satisfying $\beta_1\le z_1<\alpha$, such that $(z_1,z_2)\subset f(U\cap V)$, but $z_1\notin f(U\cap V)$ and $z_2\notin f(U\cap V)$.
\end{proof}
\begin{lem}\label{topolinter2} Let $T$ be a Hausdorff topological space, $U\ne V$ be open subsets of $T$ with neither one of them being a subset of another. Let $f:U\rightarrow I$ and $g:V\rightarrow J$, where $I$ and $J$ are some finite intervals, be homeomorphisms. Then every one of the disjoint intervals, which compose $f(U\cap V)$, is either $(z,t)$ or $(t,z)$, with $z\in I$.
\end{lem}
\begin{proof}
Since $U$ contains points not belonging to $V$, $f(U\cap V)$ must be a proper subset of $I$. Let $D$ be one of the disjoint intervals composing $f(U\cap V)$. Since $f(U\cap V)$ is open in $I$ and not equal to $I$, $D$ cannot be of the form $[m,M]$.
\\ \\
If $D$ is of the form $(z_1,M]$ or of the form $[m,z_2)$, then $z_1$ or $z_2$ must belong $I=f(U)$, because otherwise $I$ itself is $(z_1,M]$ or $[m,z_2)$ and so $I=f(U\cap V)$. The homeomorphism $g\circ f^{-1}$ must take $D$ to $(\zeta_1,M']$ or to $[m',\zeta_2)$, where $M'$ or $m'$ is $g(f^{-1}(M))$ if $D=(z_1,M]$, and $M'$ or $m'$ is $g(f^{-1}(m))$ if $D=[m,z_2)$.
\\ \\
Assume that $D$ is $(z_1,M]$ and that its image $g\circ f^{-1}(D)$ in $g(U\cap V)$ is $(\zeta_1,M']$. Since $(\zeta_1,M']$ must be open in $J$, $J$ must be either some $(m',M']$ or some $[m',M']$. If $\zeta_1\in J$ then any open neighborhood of $g^{-1}(\zeta_1)$ in $V$ intersects any open neighborhood of $f^{-1}(z_1)$ in $U$. Since $T$ is Hausdorff, $g^{-1}(\zeta_1)=f^{-1}(z_1)$ in $T$, which contradicts $z_1\notin f(U\cap V)$. If $\zeta_1\notin J$ then $J=(\zeta_1,M']$, which implies $V=g^{-1}(J)\subset f^{-1}(I)=U$, contradicting the assumptions in the lemma.
\\ \\
The same argument applies if $D$ is $[m,z_2)$ and/or if its image $g\circ f^{-1}(L)$ is $[m',\zeta_2)$. Thus, each one of the disjoint intervals composing $f(U\cap V)$ is of the form $(a,b)$. If $a\notin I$ and $b\notin I$ then $I=(a,b)=f(U\cap V)$, contradicting the assumptions in the lemma.
\end{proof}
Let $T$ be a Hausdorff topological space, $U\ne V$ be open subsets of $T$ with neither one of them being a subset of another. Let $f:U\rightarrow I$ and $g:V\rightarrow J$, where $I$ and $J$ are some finite intervals, be homeomorphisms. Let $D=(z,t)$ or $D=(t,z)$, where $z\in I$, be one of the disjoint intervals composing $f(U\cap V)$. Let $y_1>y_2>\cdots$ if $D=(z,t)$, or $y_1<y_2<\cdots$ if $D=(t,z)$, be a sequence of points in $D$ which converges to $z$ in $I$. Then $g(f^{-1}(y_1)),g(f^{-1}(y_2)),\ldots$ is a sequence of points in the image $g\circ f^{-1}(D)$ of $D$ in $J$, which converges monotonically in $\mathbb{R}$ to one of the two limiting points of the open interval $g\circ f^{-1}(D)$, and that limiting point, which we call $\zeta$, is not contained in $J$.
\\ \\
The argument why the sequence $g(f^{-1}(y_1)),g(f^{-1}(y_2)),\ldots$ converges monotonically in the closure of the finite interval $J$ in $\mathbb{R}$ is identical to the argument in the proof of Lemma \ref{InvDom}, with the intervals $[y_n,t)$ or $(t,y_n]$ here playing the role of the intervals $[-n,n]$ there. If $\zeta\in J$ then $g^{-1}(\zeta)$ must be equal to $f^{-1}(z)$ because $T$ is Hausdorff, and that contradicts $z\notin D$.
\\ \\
Thus, if $(z,t)$ is one of the disjoint intervals composing $f(U\cap V)$, and both $z$ and $t$ belong to $I$, then the image of $(z,t)$ under $g\circ f^{-1}$ must be some $(\zeta,\tau)\subset J$ with both $\zeta$ and $\tau$ not belonging to $J$, and that implies $J=(\zeta,\tau)$ and $V\subset U$. We get the following lemma:
\begin{lem}\label{topolinter3} Let $T$ be a Hausdorff topological space, $U\ne V$ be open subsets of $T$ such that neither one of them is a subset of another. If there exist homeomorphisms $f:U\rightarrow I$ and $g:V\rightarrow J$, where $I$ and $J$ are some finite intervals, then each one of the disjoint intervals composing $f(U\cap V)$ is of the form $(a,b)$ with either $a$ or $b$, but not both, belonging to $I$. Thus, $f(U\cap V)$ is either empty, or it is an open interval, or it is a disjoint union of two open intervals.
\end{lem}
In what follows we assume $I$ to be either $[t_1,t_2)$ or $(t_1,t_2)$, and one of the disjoint intervals composing $f(U\cap V)$ to be $(z_2,t_2)$, where $t_1\le z_2<t_2$ and $t_1=z_2$ can happen only when $I=[t_1,t_2)$. This assumption is done without loss of generality because we can always change $f$ to $h\circ f$ where $h(y)=t_1+t_2-y$. Now, $g\circ f^{-1}$ either takes $(z_2,t_2)$ to some $(\zeta_1,\tau_1)\subset J$, in which case $J$ is either $(\zeta_1,\zeta_2)$ or $(\zeta_1,\zeta_2]$, or takes $(z_2,t_2)$ to some $(\tau_2,\zeta_2)\in J$, in which case $J$ is either $(\zeta_1,\zeta_2)$ or $[\zeta_1,\zeta_2)$.
\\ \\
Since we can change $g$ to $h\circ g$, where $h(y)=\zeta_1+\zeta_2-y$, in each discussion we will consider only the cases when $g\circ f^{-1}$ takes $(z_2,t_2)$ to $(\zeta_1,\tau_1)\subset J$ and $J$ is either $(\zeta_1,\zeta_2)$ or $(\zeta_1,\zeta_2]$, or the cases when $g\circ f^{-1}$ takes $(z_2,t_2)$ to $(\tau_2,\zeta_2)\subset J$ and $J$ is either $(\zeta_1,\zeta_2)$ or $[\zeta_1,\zeta_2)$. We have the following three lemmas:
\begin{lem}\label{topolcomp}
Let $T$ be a Hausdorff topological space and $U\ne V$ be open subsets of $T$ such that neither one of them is a subset of the other one. Let $f:U\rightarrow I$ and $g:V\rightarrow J$, where $I$ and $J$ are some finite intervals, be homeomorphisms. If $f(U\cap V)$ consists of two disjoint intervals then $U\cup V$ is compact.
\end{lem}
\begin{proof}
Without loss of generality, we can assume that $I=(t_1,t_2)$, $J=(\zeta_1,\zeta_2)$, $f(U\cap V)=(t_1,z_1)\cup (z_2,t_2)$ where $t_1<z_1\le z_2<t_2$, and $g\circ f^{-1}$ takes $(t_1,z_1)$ to $(\zeta_1,\tau_1)$ and $(z_2,t_2)$ to $(\tau_2,\zeta_2)$ with $\zeta_1<\tau_1\le\tau_2<\zeta_2$.
\\ \\
Let $b_1\in (t_1,z_1)$, $b_2\in (z_2,t_2)$, $\beta_1=g\circ f^{-1}(b_1)$, and $\beta_2=g\circ f^{-1}(b_1)$. The continuous maps $f^{-1}$ and $g^{-1}$ take compact sets to compact sets, so the image $W_b$ of $[b_1,b_2]$ under $f^{-1}$ and the image $W_{\beta}$ of $[\beta_1,\beta_2]$ under $g^{-1}$ are compact. Hence, the union $W_b\cup W_{\beta}$ is also compact. Since $g\circ f^{-1}$ takes $[b_1,z_1)$ to $(\zeta_1,\beta_1]$, $(t_1,b_1]$ to $[\beta_1,\tau_1)$, $(z_2,b_2]$ to $[\beta_2,\zeta_2)$, and $[b_2,t_2)$ to $(\tau_2,\beta_2]$, $W_b\cup W_{\beta}$ is equal to $U\cup V$. Thus, $U\cup V$ is compact.
\end{proof}
\includegraphics[scale=0.5]{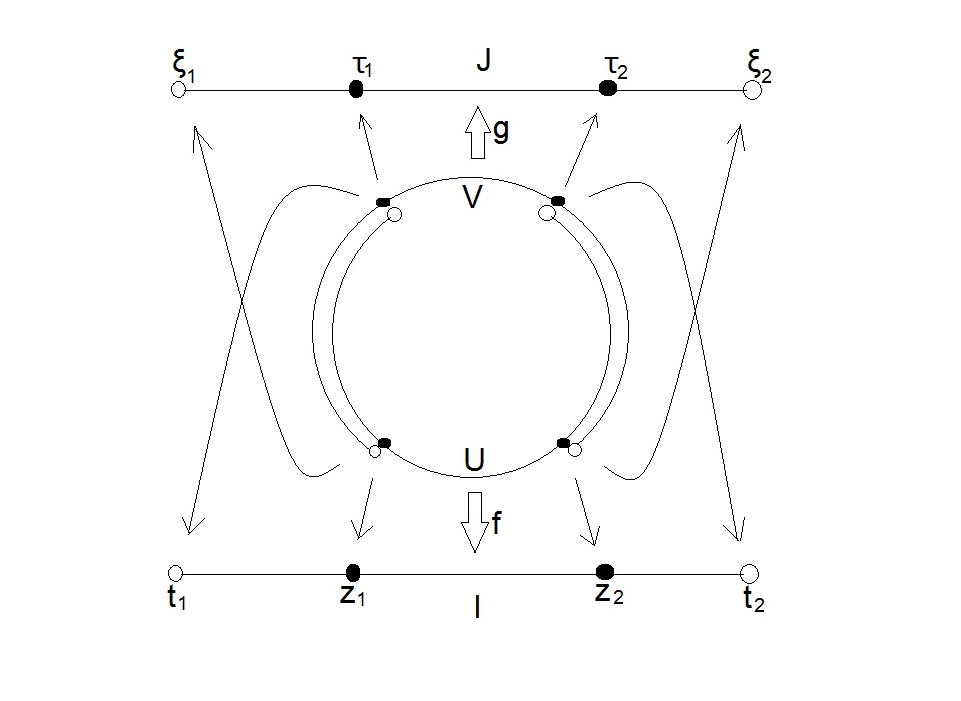}
\begin{lem}\label{topolcicrle}
Let $T$ be a Hausdorff topological space and $U\ne V$ be open subsets of $T$ such that neither one of them is a subset of the other one. Let $W$ be a subset of $T$, and let $f:U\rightarrow I$, $g:V\rightarrow J$, and $h:W\rightarrow \L$, where $I$, $J$ and $L$ are some finite intervals, be homeomorphisms. If $f(U\cap V)$ consists of two disjoint intervals then either $(U\cup V)\cap W=\emptyset$ or $W\subset U\cup V$.
\end{lem}
\begin{proof}
If $W$ intersects $U\cup V$, but is not contained in $U\cup V$, then $W$ contains points which belong to $U\cup V$ and contains points which do not belong to $U\cup V$. Since $U\cup V$ is open in $T$, and since $U\cup V$ is compact and $T$ is Hausdorff, which implies that $U\cup V$ is closed in $T$, the set of points in $W$ which are in $U\cup V$ and the set of points in $W$ which are not in $U\cup V$ are both open in $W$. Thus, the images of these sets under $h$ are open in $L$, but an interval is not equal to the union of two nonempty sets which are open in it and which do not intersect.
\end{proof}
\begin{lem}\label{TopolInterval} Let $T$ be a Hausdorff topological space, $U\ne V$ be open subsets of $T$ such that neither one of them is a subset of the other one. Let $f:U\rightarrow I$ and $g:V\rightarrow J$, where $I$ and $J$ are some finite intervals, be homeomorphisms. If $f(U\cap V)$ consists of one disjoint interval then there exists a homeomorphism $h:(U\cup V)\rightarrow L$, where $L$ is a finite interval in $\mathbb{R}$.
\end{lem}
\begin{proof}
Here $I$ is $(t_1,t_2)$ or $[t_1,t_2)$, $J$ is $(\zeta_1,\zeta_2)$ or $(\zeta_1,\zeta_2]$, and $f(U\cap V)=(z_2,t_2)$ is mapped by $g\circ f^{-1}$ to $(\zeta_1,\tau_1)$.
\\ \\
Select any $b\in (z_2,t_2)$ and let $\beta=g\circ f^{-1}(b)$.
\\ \\
If $I=(t_1,t_2)$ and $J=(\zeta_1,\zeta_2)$ then let $I'=(t_1,b]$, $J'=[\beta,\zeta_2)$, and $L=(0,2)$.
\\
If $I=[t_1,t_2)$ and $J=(\zeta_1,\zeta_2)$ then let $I'=[t_1,b]$, $J'=[\beta,\zeta_2)$, and $L=[0,2)$.
\\
If $I=(t_1,t_2)$ and $J=(\zeta_1,\zeta_2]$ then let $I'=(t_1,b]$, $J'=[\beta,\zeta_2]$, and $L=(0,2]$.
\\
If $I=[t_1,t_2)$ and $J=(\zeta_1,\zeta_2]$ then let $I'=[t_1,b]$, $J'=[\beta,\zeta_2]$, and $L=[0,2]$.
\\ \\
Define $h^{-1}:L\rightarrow (U\cup V)$ by $h^{-1}(x)=f^{-1}(x\cdot (b-t_1)+t_1)$ if $x\le 1$, and $h^{-1}(x)=g^{-1}((x-1)\cdot (\zeta_2-\beta)+\beta)$ if $x\ge 1$. In the first part of this definition we have $h^{-1}(1)=f^{-1}(b)$, and in the second part of it we have $h^{-1}(1)=g^{-1}(\beta)$, which is possible because $f^{-1}(b)=g^{-1}(\beta)$ since $\beta=g\circ f^{-1}(b)$. It is easy to verify that $h^{-1}$ is a homeomorphism between $L$ and $U\cup V$.
\end{proof}
The following theorem justifies the name ``a topological space which is a geometric realization of a graph'' in our Definition \ref{defTopolGraph}.
\begin{thm}\label{maintopology} A topological space $T$ is a geometric realization of a graph if and only if there exists a graph $G$ such that $T(G)$ is homeomorphic to $T$.
\end{thm}
\begin{proof}
Let $T$ be a topological space such that there exists a graph $G$ with $T(G)$ homeomorphic to $T$. The total topological space $Tot$, which we discussed above, is Hausdorff because it is a disjoint union of Hausdorff spaces. Since $G$ has a finite number of vertices and a finite number of edges, $Tot$ is compact. The pre-image of a point in $T(G)$ under the quotient map $q:Tot\rightarrow T(G)$ is a finite set of points in $Tot$. Moreover, every inner point $x$ in every closed subspace $T(e)$ of $Tot$ is the only point in the full pre-image under $q$ of the image $q(x)$ of $x$. In other words, $q^{-1}(\{q(x)\})=\{x\}$.
\\ \\
We now show that $T$ is Hausdorff. Let $x,y\in T(G)$ be two different points. Let $x_1,\ldots,x_i$ be all the pre-images in $Tot$ of $x$, and let $y_1,\ldots,y_j$ be all the pre-images in $Tot$ of $y$. For each $t=1,\ldots,i$ and each $r=1,\ldots,f$ we find in $Tot$ an open neighborhood $U_{t,r}$ of $x_t$ and an open neighborhoods $V_{t,r}$ of $y_r$ so that $U_{t,r}\cap V_{t,r}=\emptyset$. We can select all these open neighborhoods in such a way, that each $U_{t,r}-\{x_t\}$ and each $V_{t,r}-\{y_r\}$ is either an empty set or a nonempty subset of the set of the inner points in some $T(e)$. Thus, if a subset $U$ of the union of all the sets $U_{t,r}$ contains all the points $x_1,\ldots,x_i$ then $U$ is the full pre-image under $q$ of the image $q(U)$ of $U$ under $q$, and if a subset $V$ of the union of all the sets $V_{t,r}$ contains all the points $y_1,\ldots,y_i$ then $V$ is the full pre-image under $q$ of the image $q(V)$ of $V$ under $q$.
\\ \\
For each $t=1,\ldots,i$ let $U_t=U_{t,1}\cap\cdots\cap U_{t,j}$, and for each $r=1,\ldots,j$ let $V_r=V_{1,r}\cap\cdots\cap V_{i,r}$. Clearly, each $U_t$ is an open neighborhood of $x_t$ in $Tot$, and each $V_r$ is an open neighborhood of $y_r$ in $Tot$. It is also clear that $U_{t}\cap V_{r}=\emptyset$ for all $t$ and $r$. Thus, $U=U_1\cup\cdots\cup U_i$ and $V=V_1\cup\cdots\cup V_j$ have an empty intersection. But $U$ is the full pre-image under $q$ of $q(U)$, and $V$ is the full pre-image under $q$ of $q(V)$. Thus, $q(U)$ is an open neighborhood of $x$ in $T(G)$, and $q(V)$ is an open neighborhood of $y$ in $T(G)$, and we have $q(U)\cap q(V)=\emptyset$, which proves that $T(G)$ and $T$, which is homeomorphic to $T(G)$, are Hausdorff.
\\ \\
Since $T(G)$ is the image of $Tot$ under the quotient map, which is continuous, $T(G)$ is compact, and $T$ is compact because $T$ is homeomorphic to $T(G)$.
\\ \\
Finally, select any homeomorphism $f$ between $T$ and $T(G)$ and take any point $x\in T$. For $f(x)\in T(G)$ there exists a nonnegative integer $d$, such that there is a homeomorphism $\phi$ between some open neighborhood $\Omega$ of $f(x)$ in $T(G)$ and $\Xi_d$, which takes $f(x)$ to $0\in \Xi_d$. Then $\phi\circ f$ is a homeomorphism between the open neighborhood $f^{-1}(\Omega)$ of the point $x$ in $T$ and $\Xi_d$, which takes $x$ to $0\in \Xi_d$. The uniqueness of $d$ for $x\in T$ follows from Lemma \ref{topdegree2}. This completes the proof of the fact that if there exists a graph $G$ such that $T(G)$ is homeomorphic to $T$, then $T$ satisfies our definition of a topological space which is a geometric realization of a graph.
\\ \\
Now assume that we are given a compact Hausdorff topological space $T$ such that for every point $x\in T$ there are a nonnegative integer $d$, an open neighborhood $\Omega$ of $x$ in $T$, and a homeomorphism $\phi$ between $\Omega$ and $\Xi_d$ which satisfies $\phi(x)=0\in \Xi_d$. We will construct a graph $G$ together with a homeomorphism $h:K(G)\rightarrow T$. A simplicial complex $K$ together with a homeomorphism $h:K\rightarrow T$ is called \textit{a triangulation of the topological space $T$}.
\\ \\
First we prove that there is only a finite set of points $x$ in $T$ for which the degree $\deg(x)$ in $T$ is different from $2$. Assume that we found some infinite set of points $\{x_1,x_2,x_3,\ldots\}$ in $T$ such that their degrees in $T$ are all different from $2$. Since $T$ is compact there must exist a point $x\in T$ such that in every open neighborhood $U$ of $x$ in $T$ there are two different points $x_i$ and $x_j$ from the set $\{x_1,x_2,x_3,\ldots\}$. Indeed, otherwise taking for every $x\in T$ an open neighborhood $U$ of $x$ which contains one or zero points from $\{x_1,x_2,x_3,\ldots\}$ produces an infinite covering of $T$ by open sets, for which there is no finite sub-covering.
\\ \\
But for this point $x$ there exists a unique nonnegative integer $d$, such that there are an open neighborhood $\Omega$ of $x$ in $T$ and a homeomorphism $\phi:\Omega\rightarrow \Xi_d$ which satisfies $\phi(x)=0\in \Xi_d$. The homeomorphism $\phi$ must take at least one of the two points $x_i$ and $x_j$ in $\Omega$ to an inner point inside one of the $d$ copies of a closed half-lines $[0,\infty)$ glued together at their points $0$, which is impossible by Lemma \ref{topdegree2}. Thus, the set of all the points $x\in T$ such that $\deg(x)\ne 2$ is finite. Let $x_1,\ldots,x_n$ be all $x\in T$ such that $\deg(x)=1$ or $\deg(x)>2$, $z_1,\ldots,z_{\eta}$ be all $x\in T$ such that $\deg(x)=0$, and $T_1$ be $T-\{z_1,\ldots,z_{\eta}\}$.
\\ \\
For $i=1,\ldots,n$ and $j=1,\ldots,n$, for all $i\ne j$ let $U_{i,j}$ be an open neighborhood of $x_i$ and $U_{j,i}$ be an open neighborhood of $x_j$ such that $U_{i,j}\cap U_{j,i}=\emptyset$, and let $U_{i,i}=T$. For $i=1,\ldots,n$ we define $U_i=U_{i,1}\cap\cdots\cap U_{i,n}$. Clearly, each $U_i$ is an open neighborhood of $x_i$ in $T$, and $U_i\cap U_j=\emptyset$ when $i\ne j$. We can assume that each $U_i$ is homeomorphic to $\Xi_{d(i)}$, where $d(i)=\deg(x_i)$, by a homeomorphism $\phi_i$ satisfying $\phi_i(x)=0$, because the local base of topology of $T$ at each $x_i$ is given by such open neighborhoods.
\\ \\
For $i=1,\ldots,n$ we denote by $\Gamma_i\subset T$ the image under $\phi^{-1}_i$ of the subset of $\Xi_{d(i)}$ consisting of the point $0$ and of all the inner points $z<1$ inside each one of the $d(i)$ copies of the closed half-line $[0,\infty)$ glued together to create $\Xi_{d(i)}$, and we denote by $v_{i,1},\ldots,v_{i,d(i)}\in T$ the images under $\phi^{-1}_i$ of the $d(i)$ points $1$ inside these $d(i)$ copies of the closed half-line $[0,\infty)$.
\\ \\
Each $\Gamma_i$ is open in $T$, so the complement $T'$ of $\Gamma_1\cup\cdots\cup\Gamma_n$ in $T_1$ is a closed subset of $T$. Since $T$ is compact, $T'$ is compact. For each $j=1,\ldots,d(i)$ we restrict $\phi^{-1}_i$ to a homeomorphism between the interval $[1,2)$ inside the $j^{th}$ copy of the closed half-line $[0,\infty)$ and some open neighborhood $\Upsilon(v_{i,j})$ of $v_{i,j}$ in $T'$, and that homeomorphism takes $1$ to $v_{i,j}$. For each $y\in T'$ which is not one of the points $v_{i,j}$ there is a homeomorphism between the interval $(0,1)$ and some open neighborhood $\Upsilon(y)$ of $y$ in $T'$.
\\ \\
The collection of all the sets $\Upsilon(v_{i,j})$ and $\Upsilon(y)$ is a covering of $T'$ by subsets which are open in $T'$, and we can find some finite sub-covering. That finite sub-covering consists of all the sets $\Upsilon(v_{i,j})$ with $1\le i\le n$ and $1\le j\le d(i)$, because each point $v_{i,j}$ is not contained in any other open set in the infinite covering, and of some sets $\Upsilon(y_1),\ldots,\Upsilon(y_t)$. Next, we perform the following steps:
\\ \\
Step 1: We remove from the finite sub-covering every open set which is a proper subset of another open set in that finite sub-covering. Next, if there are several equal open sets in the finite sub-covering, we keep one of them and remove all the others from the sub-covering.
\\ \\
Step 2: If for some sets $U$ and $V$ in the finite sub-covering $U\cup V$ is homeomorphic to a disjoint union of two open intervals then we remove from our finite sub-covering all the sets $W$ which intersect with $U\cup V$. By Lemma \ref{topolcicrle}, we still have a sub-covering of $T'$. We repeat Step 2 until for every such remaining $U$ and $V$ in the finite sub-covering, there is no $W$ in the finite sub-covering which intersects $U\cup V$.
\\ \\
Step 3: If for some sets $U$ and $V$ in the finite sub-covering $U\cup V$ is homeomorphic to one open interval then we add to the sub-covering the set $U\cup V$ and, by Lemma \ref{TopolInterval}, find a homeomorphism between $U\cup V$ and either $(1,2)$ or $[1,2)$ or $[1,2]$. After doing that we return to Step 1.
\\ \\
After a finite number of repetitions of these steps, the finite sub-covering of $T'$ will consist of the sets $V_1,\ldots,V_m$ and the sets $W_1,\ldots,W_{2\mu}$, such that each $V_i$ does not intersect any other set in the sub-covering, and each pair of sets $W_{2r-1},W_{2r}$ has an intersection homeomorphic to a disjoint union of two open intervals with $W_{2r-1}\cup W_{2r}$ not intersecting any other set in the sub-covering.
\\ \\
Since $T'$ is compact $V_1,\ldots,V_m$ must be compact. Hence, each $V_t$ is homeomorphic to the closed interval $[0,1]$ and, thus, must contain exactly two different points $v_{i,j}$ and $v_{i',j'}$ which are mapped by any homeomorphism between $V_t$ and $[0,1]$ to the points $0$ and $1$. For each pair $W_{2r-1},W_{2r}$ we select any two points $w_{2r-1},w_{2r}\in W_{2r-1}\cap W_{2r}$ in such a way that $w_{2r-1}$ and $w_{2r}$ are mapped to the different copies of $(0,1)$ by any homeomorphism between $W_{2r-1}\cap W_{2r}$ and two disjoint copies of $(0,1)$.
\\ \\
We are now ready to construct the graph $G$. The vertex set of $G$ is the union $$\{z_1,\ldots,z_{\eta}\}\cup\{x_1,\ldots,x_n\}\cup\{v_{1,1},\ldots,v_{1,d(1)},\ldots,v_{n,1},\ldots,v_{n,d(n)}\}\cup\{w_1,\ldots,w_{2\mu}\}$$ The edge set of $G$ consists of the edges $e_{1,1},\ldots,e_{1,d(1)},\ldots,e_{n,1},\ldots,e_{n,d(n)}$ for which $\partial(e_{i,j})=\{x_i,v_{i,j}\}$, the edges $\alpha_1,\ldots,\alpha_m$ for which $\partial(\alpha_t)=\{v_{i,j},v_{i',j'}\}$ with $v_{i,j}$ and $v_{i',j'}$ being the two non-inner points in $V_t$, and of the edges $\beta_1,\ldots,\beta_{2\mu}$ for which $\partial(\beta_{2r-1})=\partial(\beta_{2r})=\{w_{2r-1},w_{2r}\}$.
\\ \\
We construct a homeomorphism $h:K(G)\rightarrow T$ as follows:
\\ \\
Step 1. The points $$z_1,\ldots,z_{\eta}\,x_1,\ldots,x_n,v_{1,1},\ldots,v_{1,d(1)},\ldots,v_{n,1},\ldots,v_{n,d(n)},w_1,\ldots,w_{2\mu}\in K(G)$$ are mapped by $h$ to the points with the same label and the same index in $T$;
\\ \\
Step 2. For all the points $x\in K(e_{i,j})$ for all $i$ and $j$ we define $h(x)=\phi^{-1}_i(x)$, where $\phi_i$ is the above-mentioned homeomorphism between $U_i\subset T$ and $\Xi_{d(i)}$;
\\ \\
Step 3. For all the points $x\in K(\alpha_t)$ for all $t$ we define $h$ to be any homeomorphism between $[0,1]$ and $V_t\subset T$ which takes $0\in [0,1]$ to the point $v_{i,j}\in T$ such that $0\in K(\alpha_t)$ corresponds to the vertex $v_{i,j}$ in $G$. This $h$ will take $1\in [0,1]$ to the point $v_{i',j'}\in T$ such that $1\in K(\alpha_t)$ corresponds to the vertex $v_{i',j'}$ in $G$;
\\ \\
Step 4. For all the $K(\beta_{2r-1})$ we first select any homeomorphism $\varsigma_{2r-1}$ between $(-1,2)$ and $W_{2r-1}\subset T$ such that $\varsigma_{2r-1}$ takes $0$ to the point $\varpi_0\in K(\beta_{2r-1})$ which has coordinate $0$ in $K(\beta_{2r-1})$ and takes $1$ to the point $\varpi_1\in K(\beta_{2r-1})$ which has coordinate $1$ in $K(\beta_{2r-1})$. Here $\varpi_0$ is either $w_{2r-1}$ or $w_{2r}$, and $\varpi_1$ is the other one of them. Such $\varsigma_{2r-1}$ can be constructed by taking any homeomorphism $\kappa:(-1,2)\rightarrow W_{2r-1}\subset T$ and composing it with the homeomorphism from $(-1,2)$ to itself which takes $0$ to $\kappa^{-1}(\varpi_0)$ and takes $1$ to $\kappa^{-1}(\varpi_1)$. For all $x\in K(\beta_{2r-1})$ we define $h(x)=\varsigma_{2r-1}(x)$;
\\ \\
Step 5. For all the $K(\beta_{2r})$ we first select any homeomorphism $\varsigma_{2r}$ between $(-1,2)$ and $W_{2r}\subset T$ such that $\varsigma_{2r}$ takes $0$ to the point $\varpi'_0\in K(\beta_{2r})$ which has coordinate $0$ in $K(\beta_{2r})$ and takes $1$ to the point $\varpi'_1\in K(\beta_{2r})$ which has coordinate $1$ in $K(\beta_{2r})$. Here $\varpi'_0$ is either $w_{2r-1}$ or $w_{2r}$, and $\varpi'_1$ is the other one of them. For all $x\in K(\beta_{2r})$ we define $h(x)=\varsigma_{2r}(x)$.
\\ \\
It is an easy to verify, that $h$ is a homeomorphism between $K(G)$ and $T$.
\end{proof}
Before proceeding to the topological version of Menger's Edge Theorem we demonstrate how a slight relaxation of some of the requirements in our definition of the topological space which is a geometric realization of a graph leads to exotic topological spaces, which are not homeomorphic to $T(G)$ for any graph $G$.
\\ \\
We start by relaxing the requirement on $T$ to be a Hausdorff space. First, we construct each one of the topological spaces $X_1$ and $X_2$ by taking two disjoint copies of the closed interval $[0,1]$, and gluing $0$ in the first copy of $[0,1]$ with $0$ in the second copy of $[0,1]$ and $1$ in the first copy of $[0,1]$ with $1$ in the second copy of $[0,1]$.
\\ \\
Next, we take the disjoint union $X$ of $X_1$ and $X_2$ and introduce in it an equivalence relation $\mu_1$ under which every point $0<x<1$ inside the first copy of $[0,1]$ in $X_1$ is equivalent to the point $x$ inside the first copy of $[0,1]$ in $X_2$. The quotient of $X$ by $\mu_1$ is a topological space $T_1$ which satisfies almost all the requirements in Definition \ref{defTopolGraph}, except that the points $0\in X_1$ and $0\in X_2$ under that quotient become two different points $0_1,0_2\in T_1$, and the points $1\in X_1$ and $1\in X_2$ become two different points $1_1,1_2\in T_1$, and $0_1$ and $0_2$ do not have disjoint open neighborhoods in $T_1$, and $1_1$ and $1_2$ do not have disjoint open neighborhoods in $T_1$.
\\ \\
Finally, we introduce in $X$ a new equivalence relation $\mu_2$ under which every point $0<x\le 1$ inside the first copy of $[0,1]$ in $X_1$ is equivalent to the point $x$ inside the first copy of $[0,1]$ in $X_2$, and every point $0<y\le 1$ inside the second copy of $[0,1]$ in $X_1$ is equivalent to the point $y$ inside the second copy of $[0,1]$ in $X_2$. The quotient of $X$ by $\mu_2$ is a topological space $T_2$ which satisfies almost all the requirements in Definition \ref{defTopolGraph}, except that the points $0_1,0_2\in T_2$ which come from $X_1$ and from $X_2$, respectively, do not have disjoint open neighborhoods in $T_2$.
\\ \\
In both $T_1$ and $T_2$, every point has in that space an open neighborhood which is homeomorphic to $\mathbb{R}$. Intuitively, $T_1$ is almost the geometric realization of the graph, consisting of two vertices and three edges each one of which goes between these two vertices, and $T_2$ is almost the geometric realization of the graph, consisting of two vertices and two edges each one of which goes between these two vertices.
\\ \\
Next we relax the requirement that every point $x\in T$ has an open neighborhood in $T$ homeomorphic to some $\Xi_d$. We will construct a compact Haudorff topological space $T$ in which the set of all the points which do have an open neighborhood homeomorphic to $\Xi_2$ is a dense open subset, and in which every point $x\in T$ has an open neighborhood $\Omega_x$ in $T$ and a continuous injective map $\phi_x:\Omega_x\rightarrow\Xi_3$ satisfying $\phi_x(x)=0$, such that the image of $\Omega_x$ under $\phi_x$ contains two out of the three copies of $[0,\infty)$ which compose $\Xi_3$. This  topological space $T$ is constructed as an exotic compactification of $\mathbb{R}$. Before performing that construction we briefly review compactifications of topological spaces.
\\ \\
The Stone-\v{C}ech compactification of a topological space $X$ is a compact Hausdorff topological space $\beta{X}$ together with a continuous map $\iota_X:X\rightarrow \beta{X}$ such that for any continuous map $f$ from $X$ to a compact Hausdorff topological space $Y$ there exists a unique continuous map $\beta{f}:\beta{X}\rightarrow Y$ such that $f=\beta{f}\circ \iota_X$. The standard construction of the Stone-\v{C}ech compactification is done by taking the set $C$ of all the continuous functions from $X$ to the closed interval $[0,1]$ and then considering the set $[0,1]^C$ of all the functions from the set $C$ to the set $[0,1]$. The next step is to equip $[0,1]^C$ with the product topology, the base of which is given by all the sets of functions from $C$ to $[0,1]$ which take some $f\in C$ to a real number inside some open interval $(a,b)\subset [0,1]$.
\\ \\
One can verify that $[0,1]^C$ with the product topology is a compact Hausdorff topological space and the map $\iota$ from $X$ to $[0,1]^C$, which for each $x\in X$ produces a function from $C$ to $[0,1]$ which takes $f\in C$ to $f(x)\in [0,1]$, is a continuous map. One can show that the closure in $[0,1]^C$ of the image of $X$ under $\iota$ together with the map $\iota$, regarded as a function from $X$ to $\beta{X}$, satisfy all the requirements of $\beta{X}$ and $\iota_X$.
\\ \\
If $X$ is Hausdorff and locally compact, which means that every $x\in X$ has an open neighborhood $U$ in $X$ such that the closure of $U$ in $X$ is compact, then one can show that $X$ is also a completely regular space, which means that for any closed subset $D$ of $X$ and for any $x\in X$ not belonging to $D$ there exists some continuous function $f:X\rightarrow [0,1]$ such that $f(x)=1$ and the image of $D$ under $f$ is $\{0\}\subset [0,1]$. In that and only in that case the image of $X$ in $\beta{X}$ under $\iota_X$ is an open dense subset of $\beta{X}$ and $\iota_X$ is a homeomorphism between $X$ and that open dense subset of $\beta{X}$.
\\ \\
In general, a compactification of a locally compact Hausdorff topological space $X$ is a compact Hausdorff topological space $Y$ and an injective continuous map $f:X\rightarrow Y$ such that the image $f(X)$ of $X$ in $Y$ is an open dense subset of $Y$, and that $f$ is a homeomorphism between $X$ and $f(X)$. Intuitively, the Stone-\v{C}ech compactification of $X$ is the largest compactification of $X$. If $X$ is compact then any compactification $Y$ and $f$ of $X$ must satisfy $f(X)=Y$. However, if $X$ is not compact then $Y$ must have at least one more point than $f(X)$.
\\ \\
The Alexandroff compactification of a locally compact Hausdorff topological space $X$ which is not compact is a compact Hausdorff topological space $\dot{X}$, which is the space $X$ with one point $\infty$ added to it, and the map $f$ which takes each $x\in X$ to the same $x\in \dot{X}$. The open subsets of $\dot{X}$ are $\dot{X}$ itself, all the open subsets of $X$, and all the subsets of $\dot{X}$ such that their complement in $\dot{X}$ is a compact subset of $X\subset \dot{X}$. Intuitively, this is the smallest possible compactification of a locally compact Hausdorff topological space which is not compact. Indeed, if $Y$ and $g:X\rightarrow Y$ is any other compactification of $X$ then the continuous surjective map $h:Y\rightarrow \dot{X}$ defined as $h(y)=g^{-1}(y)\in X\subset\dot{X}$ if $y\in g(X)$ and $h(y)=\infty\in \dot{X}$ if $y\in Y-g(X)$ satisfies $f=h\circ g$.
\\ \\
Since every finite closed interval in $\mathbb{R}$ is compact, $\mathbb{R}$ is locally compact. It is trivial to verify that $\mathbb{R}$ is Hausdorff and is not compact. What we have shown here above is that if for a compactification $Y$ and $g:\mathbb{R}\rightarrow Y$ of $\mathbb{R}$ every $y\in Y$ has an open neighborhood $\Omega_y$ in $Y$ and a homeomorphic $\phi_y:\Omega_y\rightarrow \Xi_d$ satisfying $\phi_y(y)=0\in \Xi_d$ then either this is the Alexandroff compactification and $Y=\dot{\mathbb{R}}$ is homeomorphic to a circle or this is the standard two-point compactification and $Y=\mathbb{R}\cup\{-\infty,\infty\}$ is homeomorphic to $[-1,1]$.
\\ \\
Let now $Y$ be a disjoint union of the sets $\mathbb{R}$, $C^{-}$ and $C^{+}$, where $C^{-}$ and $C^{+}$ as sets are copies of the interval $[0,1)$. We equip $Y$ with the following topology:
\\ \\
The local base of topology at any $x\in \mathbb{R}\subset Y$ is given by the sets $(x-\epsilon,x+\epsilon)\subset \mathbb{R}$ where $\epsilon$ is any positive real number;
\\ \\
The local base of topology at any $x\ne 0$ in $C^{+}$ is given by the sets which are the union of any set $(x-\epsilon,x+\epsilon)\subset C^{+}$ with $0<\epsilon<\min(x,1-x)$ and of all the sets $(n+x-\epsilon,n+x+\epsilon)\subset \mathbb{R}$ with $n=N,N+1,N+2,\ldots$, where $N$ is some integer. The local base of topology at $0\in C^{+}$ is given by the sets which are the union of any two sets $[0,\epsilon)\subset C^{+}$ and $(1-\epsilon,1)\subset C^{+}$ with $0<\epsilon<0.5$ and of all the sets $(n-\epsilon,n+\epsilon)\subset \mathbb{R}$ with $n=N,N+1,N+2\ldots$, where $N$ is some integer;
\\ \\
The local base of topology at any $x\ne 0$ in $C^{-}$ is given by the sets which are the union of any set $(x-\epsilon,x+\epsilon)\subset C^{-}$ with $0<\epsilon<\min(x,1-x)$ and of all the sets $(n+x-\epsilon,n+x+\epsilon)\subset \mathbb{R}$ with $n=N,N-1,N-2,\ldots$, where $N$ is some integer. The local base of topology at $0\in C^{-}$ is given by the sets which are the union of any two sets $[0,\epsilon)\subset C^{-}$ and $(1-\epsilon,1)\subset C^{-}$ with $0<\epsilon<0.5$ and of all the sets $(n-\epsilon,n+\epsilon)\subset \mathbb{R}$ with $n=N,N-1,N-2,\ldots$, where $N$ is some integer.
\\ \\
The function $f:\mathbb{R}\rightarrow Y$ in this compactification takes each $x\in \mathbb{R}$ to $x\in \mathbb{R}\subset Y$. Intuitively, $C^{+}$ is the circle in the positive infinity and $C^{-}$ is the circle in the negative infinity, and the line $\mathbb{R}$ is an infinite slinky which converges to $C^{-}$ in the left direction and converges to $C^{+}$ in the right direction. One can to embed this compactification in $\mathbb{R}^2$ as two circles $x=\cos(t),y=\sin(t)$ and $x=3\cos(t),y=3\sin(t)$ and one infinite spiral $x=\frac{2\arctan(t)+2\Pi}{\Pi}\cos(t),y=\frac{2\arctan(t)+2\Pi}{\Pi}\sin(t)$.
\\ \\
Every $y\in C^{-}\cup C^{+}$ has an open neighborhood $\Omega_y$ in $Y$ and a continuous injective map $\phi_y:\Omega_y\rightarrow\Xi_3$ satisfying $\phi_y(y)=0$, such that the image of $\Omega_y$ under $\phi_y$ contains two out of the three copies of $[0,\infty)$ composing $\Xi_3$.
\\ \\
Any continuous function $h:\mathbb{R}\rightarrow [0,1]$ which converges asymptotically to any two continuous periodic functions with periods $\frac{1}{p}$ and $\frac{1}{q}$ for any positive integers $p$ and $q$, one as $x\rightarrow -\infty$ and the other one as $x\rightarrow +\infty$, can be extended as a continuous function $h:Y\rightarrow [0,1]$.
\\
\includegraphics[scale=0.5]{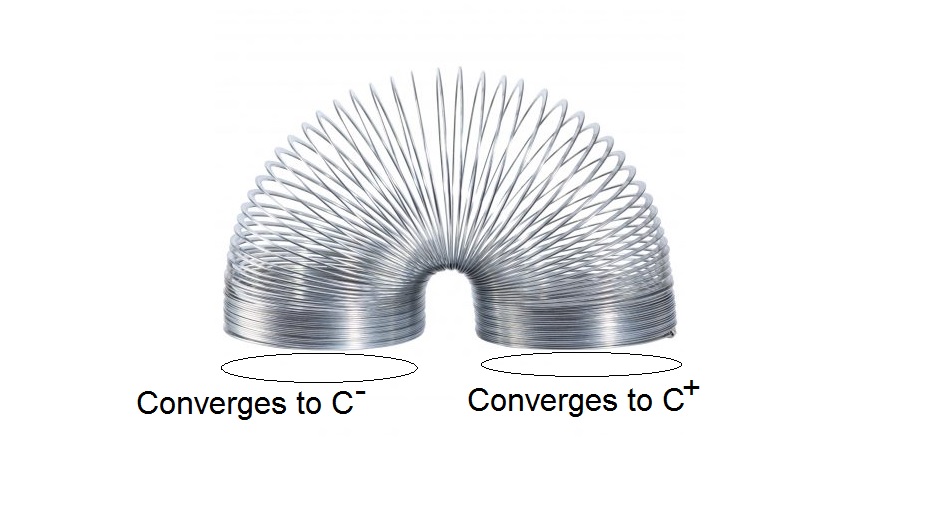}
\\
We now return to our discussion of topological spaces which are geometric realizations of graphs. Let $T$ be a topological space which is a geometric realization of a graph. For a graph $G$, for which $T$ is homeomorphic to $T(G)$, we take the formal sum of all the edges in $G$ and apply the boundary map $\partial$ to that formal sum. Regardless of the particular choice of the graph $G$, the boundary map will return the formal sum of all the vertices in $G$ which correspond to the points in $T$ which have odd degrees in $T$. Thus, we define \textit{the boundary of $T$}, which we denote by $\partial(T)$, as the set of all the points in $T$ which have odd degrees in $T$. If $T$ does not contain any points with odd degrees then $\partial(T)=\emptyset$. Now we state and prove the topological version of Menger's Edge Theorem:
\begin{thm}\label{Menger2Topol} Let $T$ be a topological space which is a geometric realization of a graph, and let $D=\{v_1,v_2\}$ be a subset of $T$. The following two statements are equivalent:
\\
I: For any $k-1$ or less pairwise disjoint nonempty open subsets $\Gamma_1,\ldots,\Gamma_{\vartheta}$ of $T$ we can find some nonempty open subsets $\Gamma'_1\subset\Gamma_1,\ldots,\Gamma'_{\vartheta}\subset\Gamma_{\vartheta}$ such that after the deletion of all the subsets $\Gamma'_1,\ldots,\Gamma'_{\vartheta}$ from $T$ there still exists a closed subspace $P$ of $T$ such that $P$ is a geometric realization of a graph and $\partial(P)=D$;
\\
II: There exist $k$ closed subspaces $P_1,\ldots,P_k$ in $T$ each one of which is a geometric realization of a graph, such that for all $j\ne i$ the subset $P_i\cap P_j$ of $T$ does not contain a subset homeomorphic to $\mathbb{R}$, for which $\partial(P_1)=\cdots=\partial(P_k)=D$.
\end{thm}
\begin{proof}
For an isolated point $x$ in $T$ which is contained in any closed subspace $P$ of $T$ which is a geometric realization of a graph $L$, the degree of $x$ in $P$ is $0$, and $x\in P$ corresponds to a vertex of degree $0$ in $L$. Thus, the isolated points in $T$ have no effect on both statements in our theorem, and we can assume that there are no isolated points in $T$.
\\ \\
Let $G$ be any graph such that $T$ is homeomorphic to $K(G)$ by a homeomorphism $\phi$ which takes the points $v_1$ and $v_2$ in $T$ to the points in $K(G)$ which underly two vertices in $G$, and we call these vertices $v_1$ and $v_2$. We fix any such $\phi$, and speak of $T$ and of $K(G)$ interchangeably.
\\ \\
First we prove that II follows from I:
\\ \\
Suppose that someone selected some $k-1$ or less edges $e_1,\ldots,e_{\vartheta}$ in $G$ for a deletion. Let each $\partial(e_i)$ in $G$ be $v_{i,1}+v_{i,2}$ such that in $K(e_i)$ the point $K(v_{i,1})$ has coordinate $0$ and the point $K(v_{i,2})$ has coordinate $1$. Let graph $\tilde{G}$ be $G$ without the edges $e_1,\ldots,e_{\vartheta}$. We want to find a chain of edges $p$ in $\tilde{G}$ such that its boundary $\partial(p)$ in $\tilde{G}$ is $v_1+v_2$.
\\ \\
To the edges $e_1,\ldots,e_{\vartheta}$ we correspond $\vartheta$ pairwise disjoint open subsets $\Gamma_1,\ldots,\Gamma_{\vartheta}$ in $T=K(G)$, where each $\Gamma_i$ is the set of all the inner points in the closed subspace $K(e_i)$ of $K(G)$. Thus, each $\Gamma_i$ is a copy of $(0,1)$.
\\ \\
As per part I, we can find some nonempty open subsets $\Gamma'_1\subset\Gamma_1,\ldots,\Gamma'_{\vartheta}\subset\Gamma_{\vartheta}$ such that after the deletion of all the $\Gamma'_1,\ldots,\Gamma'_{\vartheta}$ from $T$ we can still find a closed subspace $P$ of $T$ such that $P$ is a geometric realization of a graph and that the boundary $\partial(P)$ of the topological space $P$ is $D$. If $P$ contains any isolated points then $P$ without all these isolated points is also a closed subspace of $T$, and is a geometric realization of a graph, and its boundary is $D$. Thus, we assume that we found $P$ which has no isolated points.
\\ \\
We can make the open sets $\Gamma'_1,\ldots,\Gamma'_{\vartheta}$ smaller without effecting $P\subset T$. Thus, we assume that each subset $\Gamma'_1$ is some $(\alpha_i,\beta_i)\subset \Gamma_i$ with $0<\alpha_i<\beta_i<1$. Let $G'$ be the graph $\tilde{G}$ with the new vertices, which we, abusing the notation, denote by $\alpha_1,\ldots, \alpha_{\vartheta},\beta_1,\ldots, \beta_{\vartheta}$, and the new edges $e_{1,1},\ldots,e_{1,\vartheta},e_{2,1},\ldots,e_{2,\vartheta}$, such that $\partial(e_{1,i})=v_{i,1}+\alpha_i$ and $\partial(e_{2,i})=v_{i,2}+\beta_i$ for $i=1,\ldots,\vartheta$, drawn in it.
\\ \\
The topological space $T'$, which is obtained from $T=K(G)$ by deleting from it all the open sets $\Gamma'_1,\ldots,\Gamma'_{\vartheta}$, is homeomorphic to $K(G')$ by the homeomorphism $\psi$ which is defined as the identity map on $T-(\Gamma_1\cup\cdots\cup\Gamma_{\vartheta})$, and on each $K(e_{i})-(\alpha_i,\beta_i)$ is defined by $\psi(x)=\frac{\phi(x)}{\alpha_i}\in K(e_{1,i})$ if $x\le\alpha_i$ and $\psi(x)=\frac{\phi(x)-\beta_i}{1-\beta_i}\in K(e_{2,i})$ if $\beta_i\le x$. We fix $\psi$ and speak of $T'$ and of $K(G')$ interchangeably.
\\ \\
Now we will show that in $K(G')$ the subset $P$ is underlying $K(L')$ for some subgraph $L'$ of $G'$. Then we will show that this $L'$ is also a subgraph of $\tilde{G}$, which implies that in $K(\tilde{G})$ the subset $P$ is underlying $K(L')$.
\\ \\
Let $\Lambda$ be the set of all the edges $\lambda$ in $G'$ for which the intersection of the closed subsets $K(\lambda)$ and $P$ of $K(G')$ contains a subset homeomorphic to $\mathbb{R}$. Let $L$ be any graph such that $P$ is homeomorphic to $K(L)$.
\\ \\
The union of all the closed subsets $K(\epsilon)$ of $K(G')$, as $\epsilon$ runs over all the edges in $G'$, contains $K(e)$ for any edge $e$ in $L$, but only a finite number of points in $K(G')$ correspond to the vertices in $G'$. Thus, the disjoint union $\Delta(G')$ of all the open subsets $\Delta(\epsilon)$ of $K(G')$, where each $\Delta(\epsilon)$ is the set of all the inner points in $K(\epsilon)$, as $\epsilon$ runs over all the edges in $G'$, contains $\Delta(e)$, which is the set of all the inner points in $K(e)$, without a finite number of points. Since each $\Delta(e)$ is a copy of $(0,1)$, $\Delta(G')$ contains each $\Delta(e)$ without some points $0<\delta(e,1)<\cdots<\delta(e,\chi(e))<1$.
\\ \\
Each one of the open intervals $(0,\delta(e,1)), (\delta(e,1),\delta(e,2)),\ldots, (\delta(e,\chi(e)),1)$ is an open subset of $\Delta(\epsilon)$ for some edge $\epsilon$ in $G'$, and is homeomorphic to $\mathbb{R}$. Thus, the union $\Upsilon$ of all the closed subsets $K(\lambda)$ of $K(G')$, as $\lambda$ runs over all the edges in $\Lambda$, contains $K(e)$. Since this is true for every edge $e$ in $L$, and since $P$ does not contain any isolated points, $\Upsilon$ contains $P$.
\\ \\
Now assume that some $u\in \Upsilon$ is not contained in $P$. We can find an edge $\lambda$ in $G'$ such that the closed subset $K(\lambda)$ of $K(G')$ contains $u$, and that the intersection of $K(\lambda)$ with $P$ contains an open interval $(a,b)\subset K(\lambda)$. Thus, for any $a<w<b$ in $K(\lambda)$ we have $w\in P$, and for some $0\le u\le 1$ in $K(\lambda)$ we have $u\notin P$.
\\ \\
By the B$-$W argument there exists some $\min(w,u)\le\theta\le\max(w,u)$ in $K(\lambda)$, such that any open neighborhood of $\theta$ in $K(\lambda)$ contains points which belong to $P$ and points which do not belong to $P$. Thus, any open neighborhood of $\theta$ in $K(G')$ contains points which belong to $P$ and points which do not belong to $P$. Since $P$ is closed in $K(G')$, $\theta$ belongs to $P$. Thus, $\theta\ne u$, which implies that $0<\theta<1$, because $0<w<1$. We get that $\deg(\theta)$ in $P$ cannot be $2$, and thus must be $1$, which means that either $\theta=v_1$ or $\theta=v_2$, and that is impossible because $0<\theta<1$ in $K(\lambda)$. So $P$ contains $\Upsilon$.
\\ \\
Thus, $\Upsilon=P$, which implies that in $K(G')$ the subset $P$ is underlying $K(L')$, where $L'$ is the subgraph of $G'$, whose edge set is $\Lambda$ and whose vertex set is the set of all the vertices in $G'$ which belong to the boundaries of the edges in $\Lambda$. Since $\Lambda$ cannot contain any edges $e_{1,1},\ldots,e_{1,\vartheta},e_{2,1},\ldots,e_{2,\vartheta}$, which were drawn in the graph $\tilde{G}$ to create the graph $G'$, because each one of these edges contains a vertex of degree $1$ in $G'$ different from $v_1$ and $v_2$, all the edges in $\Lambda$ belong to the graph $\tilde{G}$. Thus, $L'$ is a subgraph of $\tilde{G}$, and the subset $P$ of $K(\tilde{G})$ is underlying $K(L')$ in $K(\tilde{G})$.
\\ \\
Let the chain of edges $p$ in $\tilde{G}$ be the sum of all the edges in $\Lambda$. The boundary $\partial(p)$ of $p$ in $\tilde{G}$ is $v_1+v_2$. Since this is true for $\tilde{G}$ obtained by deleting any $k-1$ or less edges in $G$, by Theorem \ref{Menger.homol.2} there exist $k$ chains of edges $p_1,\ldots,p_k$ in $G$, such that no two of $p_1,\ldots,p_k$ have any common nontrivial summands and that $\partial(p_1)=\cdots=\partial(p_k)=v_1+v_2$ in $G$.
\\ \\
Finally, let the closed subspace $P_i$ of $K(G)$, for each $i=1,\ldots,k$, be the union of all $K(e)$ as $e$ runs over all the edges in $G$ which appear as summands in $p_i$. It is easy to verify that the subspaces $P_1,\ldots,P_k$ satisfy the requirements in the part II of the theorem.
\\ \\
The proof that I follows from II:
\\ \\
We can remove the isolated points from $P_1,\ldots,P_k$, so we assume that none of them has isolated points. For each $P_i$, let $\Lambda_i$ be the set of all the edges $\lambda$ in $G$ for which the intersection of the subsets $K(\lambda)$ and $P_i$ of $K(G)$ contains a subset homeomorphic to $\mathbb{R}$, and let $L_i$ be the graph with the set of edges $\Lambda_i$ and with the set of vertices consisting of all the vertices with appear in the boundaries of the edges in $\Lambda_i$. As was established above, each $P_i$ underlies its $K(L_i)$ in $K(G)$. Since any intersection $P_i\cap P_j$ when $j\ne i$ does not contain a subset homeomorphic to $\mathbb{R}$, there is no edge in $G$ which belongs to $\Lambda_i$ and to $\Lambda_j$. Let each $\Lambda_i$ be $\{e_{i,1},\ldots, e_{i,\varrho(i)}\}$. Let each $\Delta(e_{i,j})$ be the set of all the inner points in $K(e_{i,j})$.
\\ \\
Assume that someone selected any $k-1$ or less pairwise disjoint open subsets $\Gamma_1,\ldots,\Gamma_{\vartheta}$ of $K(G)$. For each $\Gamma_r$ which contains a point not belonging to the union $X$ of all the closed subsets $K(e_{i,j})$ of $K(G)$ we find some open neighborhood $\Gamma'_r$ of that point in $\Gamma_r$ which does not intersect $X$. For each $\Gamma_p$ contained in $X$ we find a nonempty open subset $\Gamma'_p$ of $\Gamma_p$ which is contained in one of the open subsets $\Delta(e_{i,j})$ of $K(G)$. Thus, we obtain some nonempty open subsets $\Gamma'_1\subset\Gamma_1,\ldots,\Gamma'_{\vartheta}\subset\Gamma_{\vartheta}$ such that the deletion of all of them from $K(G)$ effects at most $k-1$ of the $k$ closed subspaces $P_1,\ldots,P_k$ of $K(G)$. Part I of the theorem follows.
\end{proof}
Similarly, the topological version of our Theorem \ref{Menger.homol.5} states:
\begin{thm}\label{Menger4Topol} 
Let $T$ be a topological space which is a geometric realization of a graph, and let $D=\{v_1,v_2,v_3,v_4\}$ be a subset of $T$. If $\partial(T)=D$ or $\partial(T)\subset D$ then the following two statements are equivalent:
\\
I: For any $k-1$ or less pairwise disjoint nonempty open subsets $\Gamma_1,\ldots,\Gamma_{\vartheta}$ of $T$ we can find some nonempty open subsets $\Gamma'_1\subset\Gamma_1,\ldots,\Gamma'_{\vartheta}\subset\Gamma_{\vartheta}$ such that after the deletion of all the subsets $\Gamma'_1,\ldots,\Gamma'_{\vartheta}$ from $T$ there still exists a closed subspace $P$ of $T$ such that $P$ is a geometric realization of a graph and $\partial(P)=D$;
\\
II: There exist $k$ closed subspaces $P_1,\ldots,P_k$ in $T$ each one of which is a geometric realization of a graph, such that for all $j\ne i$ the subset $P_i\cap P_j$ of $T$ does not contain a subset homeomorphic to $\mathbb{R}$, for which $\partial(P_1)=\cdots=\partial(P_k)=D$.
\end{thm}
\begin{proof} Let $G$ be any graph such that $T$ is homeomorphic to $K(G)$ by a homeomorphism $\phi$ which takes the points $v_1,v_2,v_3,v_4$ in $T$ to the points in $K(G)$ which underly four vertices in $G$, and we call these vertices $v_1,v_2,v_3,v_4$. We fix any such $\phi$, and speak of $T$ and of $K(G)$ interchangeably. Since either $K(G)$ has no boundary or the boundary of $K(G)$ is contained in $D$, all the other vertices in $G$ have even degrees in $G$. The rest of the proof that II follows from I, and the proof that I follows from II can be copied from the proof of Theorem \label{Menger2Topol}.
\end{proof}
In our triangulation of the topological space which is a geometric realization of a graph we used the approach which appears in the Appendix in \cite{Milnor}. For an exposition on the Stone-\v{C}ech compactification and the related notions in Topology we refer to \cite{Bell}. For a treatment of simplicial complexes, chains of simplices, the boundary operation, cycles, boundaries, homology, topological spaces which are geometric realizations of simplicial complexes, and their topological properties, we refer to \cite{HW}, \cite{Munkres}, \cite{Hatcher}. Homological paths, which appear in this work, were introduced in \cite{Chern}.

\end{document}